\documentclass{article}
\usepackage[utf8]{inputenc}

\usepackage[style=alphabetic,natbib,giveninits=true,doi=false,isbn=false,url=false,eprint=false]{biblatex}
\AtEveryBibitem{\clearlist{language}}
\AtEveryBibitem{\clearfield{pages}} 
\usepackage{soul}
\usepackage{amsmath}
\usepackage{amssymb}
\usepackage{amsthm}
\usepackage{mathtools}

\usepackage{todonotes}

\usepackage{pgf,tikz,pgfplots}

\usepackage{tikz-cd}
\usetikzlibrary{patterns}
\usetikzlibrary{spy}
\usetikzlibrary{arrows,decorations.markings}
\usetikzlibrary{math}

\usepackage{dirtytalk}

\usepackage[mathscr]{euscript}
\usepackage{enumerate}
\usepackage[shortlabels]{enumitem}

\usepackage{subfigure}

\usepackage[margin=1in]{geometry}

\usepackage{xcolor}

\usepackage{setspace}
\theoremstyle{plain}
\newtheorem{thm}{Theorem}

\theoremstyle{plain}
\newtheorem{lemma}{Lemma}

\theoremstyle{plain}
\newtheorem{prop}{Proposition}

\theoremstyle{plain}

\theoremstyle{definition}
\newtheorem{definition}{Definition}

\theoremstyle{definition}
\newtheorem{remark}{Remark}

\usepackage[T2A,T1]{fontenc}

\usepackage[symbol]{footmisc}

\usepackage{array}
\usepackage{makecell}

\renewcommand{\d}{\,\mathrm{d}}
\newcommand{\limn}{\lim_{n\rightarrow \infty}}

\newcommand{\tor}{\mathbb{T}^2}

\usepackage{subfiles}

\usepackage{multicol}
\providecommand{\Acknowledgements}[1]
{
  \small	
  \textbf{\textit{Acknowledgements---}} #1
}

\title{Loss of Exponential Mixing in a Non-Monotonic Toral Map}

\author{J. Myers Hill$^{1\,\star}$, R. Sturman$^{1}$, M. C. T. Wilson$^{2}$  \\
        \small $^{1}$School of Mathematics, University of Leeds, Leeds LS2 9JT, United Kingdom \\
        \small $^{2}$School of Mechanical Engineering, University of Leeds, Leeds LS2 9JT, United Kingdom \\
        \small $^\star$ E: j.d.myershill@leeds.ac.uk
}

\pgfplotsset{compat=1.17}

\begin{document}

\date{}

\maketitle

\begin{abstract}

We consider a Lebesgue measure preserving map of the 2-torus, given by the composition of orthogonal tent shaped shears. We establish strong mixing properties with respect to the invariant measure and polynomial decay of correlations for H\"older observables, making use of results from the chaotic billiards literature. The system serves as a prototype example of piecewise linear maps which sit on the boundary of ergodicity, possessing null measure sets around which mixing is slowed and which birth elliptic islands under certain perturbations. 
 
\end{abstract}

\Acknowledgements{JMH supported by EPSRC under Grant Refs. EP/L01615X/1 and EP/W524372/1.}


\section{Introduction}
The statistics of chaotic dynamics driven by an area-preserving map are often described by its mixing properties. Given such a map $f:X \to X$, preserving a measure $\mu$, we say that $f$ is \emph{mixing} if
its correlations $C_n(\varphi,\psi)$ decay to 0 for $L^2$ observables $\varphi,\psi:X \to \mathbb{R}$, where
\[  C_n(\varphi,\psi) =  \int \left(\varphi \circ f^n \right) \psi \d \mu - \int \varphi \d \mu \int \psi \d \mu \]
denotes the \emph{correlation function}.

Given $C_n \to 0$, the speed at which these correlations decay (the \emph{mixing rate}) further characterises a map's dynamics. We say that $f$ enjoys exponential decay of correlations if there exists constants $0<\theta<1$ and $c(\varphi,\psi)>0$ such that
\begin{equation}
    \label{eq:expMixing}
    |C_n(\varphi,\psi)| \leq c \, \theta^n.
\end{equation}
Similarly we say that $f$ enjoys polynomial decay of correlations if there exists $\alpha>0$ and $c(\varphi,\psi)>0$ such that
\begin{equation}
    \label{eq:polyMixing}
    |C_n(\varphi,\psi)| \leq c \, n^{-\alpha}.
\end{equation}
Some regularity on the observables is typically assumed; we assume H\"older continuity throughout this article. Other statistical properties (e.g. the central limit theorem) are also intimately linked this rate of decay.

Building on \cite{young_statistical_1998,young_recurrence_1999,chernov_decay_1999,markarian_billiards_2004}, \cite{chernov_billiards_2005} gives conditions under which a uniformly hyperbolic map with singularities satisfies (\ref{eq:expMixing}). These include mild restrictions on the nature of the singularities and regularity of local manifolds, alongside a \emph{one-step expansion} estimate which ensures expansion by hyperbolicity dominates the cutting by singularities:
\begin{equation}
    \label{eq:oldOneStepExpansion}
     \liminf_{\delta \to 0} \sup_{W: |W|<\delta} \sum_i \frac{|W_i|}{|f(W_i)|} < 1 
\end{equation}
where the supremum is taken over unstable manifolds $W$, of length $|W|$, split into components $W_i$ by the singularities. Key to this analysis is construction of a \emph{Young tower} \cite{young_statistical_1998}. Given a subset $A \subset X$ and $x \in A$, define the 
\emph{return time} of $x$ to $A$ under $f$ as $R(x;f,A) = \inf \{ i>0 \, | \, f^i(x) \in A  \}$. Young considers returns\footnote{In particular `good' returns which satisfy additional technical constrains, see \cite{young_statistical_1998}. A precise definition of \emph{hyperbolic product structure} is also found therein.} to some subset $\Delta_0$ (the tower base) with hyperbolic product structure, showing that if returns satisfy an exponential tail bound:
\begin{equation}
    \label{eq:expRtailYoung}
    \mu( \{ x \in \Delta_0  \, | \, R(x;f,\Delta_0) > n\}) < C \theta^n,
\end{equation}
then (\ref{eq:expMixing}) holds. Explicitly constructing $\Delta_0$ is challenging in many systems, as is estimating its recurrence, requiring all the iterates of $f$ to be considered. The scheme of \cite{chernov_billiards_2005} both avoids the explicit construction of $\Delta_0$ and reduces the analysis down to conditions such as (\ref{eq:oldOneStepExpansion}) concerning a single iterate of the map $f$.

The scheme has utility beyond uniformly hyperbolic examples. Following \cite{young_recurrence_1999}, if there exists $C>0$ and $\alpha>0$ such that
\begin{equation}
    \label{eq:polyRtailYoung}
    \mu( \{ x \in \Delta_0  \, | \, R(x;f,\Delta_0) > n\}) < C n^{-\alpha},
\end{equation}
then $f$ satisfies (\ref{eq:polyMixing}). Suppose $f:X \to X$ has suspected polynomial decay of correlations, non-uniformly hyperbolic and possessing some region $N$ where $f$ is non-hyperbolic with escape times $E(x;f,N) = \inf \{i>0 \, | \, f^i(x) \notin N \}$ satisfying
\begin{equation}
\label{eq:escapeDist}
     \mu( \{ x \in N  \, | \, E(x;f,N) > n\}) < C n^{-\alpha}.
\end{equation}
By non-uniform hyperbolicity, a.e. $x \in N$ eventually escapes and hits some region of `strong' hyperbolicity, precisely a subset $M \subset X$ with uniformly hyperbolic \emph{return map} $f_M(x) = f^R(x)$, $R=R(x;f,M)$. Using its strong hyperbolic properties to satisfy the conditions of \cite{chernov_billiards_2005}, $f_M$ then admits a Young tower with base $\Delta_0 \subset M$, satisfying
\begin{equation}
\label{eq:expTailReturnMap}
    \mu( \{ x \in M  \, | \, R(x;f_M,\Delta_0) > n\}) < C \theta^n.
\end{equation}
Extending the domain of $f_M$ to $X$ in the obvious fashion, the bound (\ref{eq:escapeDist}) suggests
\begin{equation}
\label{eq:E1}
     \mu( \{ x \in X  \, | \, R(x;f,M) > n\}) < C n^{-\alpha}, 
\end{equation}
which can be extended, making use of (\ref{eq:expTailReturnMap}), to give (\ref{eq:polyRtailYoung}). This final step is non-trivial and typically relies on utilising precise mapping behaviour of $f_M$. The above scheme has been used to establish polynomial decay of correlations for various billiards maps including certain stadia and tables with cusps \cite{chernov_billiards_2005,chernov_improved_2008}. Beyond billiards, in \cite{springham_polynomial_2014} $\mathcal{O}(1/n)$ correlation decay was shown for a family of \emph{linked twist maps} (hereafter LTMs). These are Lebesgue measure preserving continuous maps on the 2-torus $\tor$, composing monotonic shears restricted to horizontal and vertical annuli $P,Q \subsetneq \tor$. Here, mixing is slowed by orbits remaining trapped in $P \triangle Q$\footnote{The symmetric difference $P \triangle Q = (P \setminus Q) \cup (Q \setminus P)$.} for arbitrarily long periods, with recurrence to $M = P \cap Q$ satisfying the tail bound (\ref{eq:E1}). Monotonicity of the shears was important in the analysis, allowing for a straightforward proof of the mixing property.

More recently the scheme was directly applied to a family of non-monotonic toral maps \cite{myers_hill_exponential_2022}. Parameterising $\tor$ by $(x,y) \in \mathbb{R}^2 / \mathbb{Z}^2$, these maps similarly compose horizontal and vertical shears $H_{(\xi,\eta)} = G \circ F$ where
\[ F(x,y) =
\begin{cases}
\left(  x + \frac{y}{1-\eta}, y  \right) \text{ mod 1 } & \text{ for } y \leq 1-\eta, \\

\left(  x + \frac{1-y}{\eta}, y  \right) \text{ mod 1 } & \text{ for } y \geq 1-\eta, \\
\end{cases}
\quad
G(x,y) =
\begin{cases}
\left(  x  , y + \frac{x}{1-\xi}  \right) \text{ mod 1 } & \text{ for } x \leq 1-\xi, \\

\left(  x  , y + \frac{1-x}{\xi}  \right) \text{ mod 1 } & \text{ for } x \geq 1-\xi, \\
\end{cases}\]
and $0<\xi,\eta<1$. Exponential mixing rates were established over a wide neighbourhood of $(\xi,\eta) = (0,0)$, with boundary determined by (\ref{eq:oldOneStepExpansion}), and are expected over $1-\frac{1}{4 \xi} < \eta < \frac{1}{4-4 \xi}$ where $H_{(\xi,\eta)}$ is uniformly hyperbolic (with singularities). This includes the parameter subspace $\eta=\xi$ corresponding to matching $F$ and $G$, excluding the cusp $\xi = \eta = 1/2$ where $F$ and $G$ are symmetric tent maps. This cusp is notable in the transverse subspace $\eta = 1- \xi$ also, being the only parameters for which elliptic islands do not form. Following \cite{myers_hill_exponential_2022}, we refer to $H_{(\xi,\eta)}$ at these precise parameters as the \emph{orthogonal tents map} (OTM) and denote it simply by $H$. 

\begin{figure}
    \centering
    \begin{tikzpicture}
    \definecolor{white}{RGB}{255,255,255}
    \definecolor{tomato}{RGB}{255, 99, 71}
    \definecolor{teal}{RGB}{95, 158, 160} 
    \node at (-5,0) {
    \begin{tikzpicture}[scale=0.4]

    \draw[ultra thick, tomato] (0,7.5) -- (2.5,10);
    \node[text=tomato] at (1.5,8.2) {$l_1$};
    \draw[ultra thick, gray] (2.5,5) -- (5,7.5);
    \node[text=gray] at (3.5,6.8) {$l_2$};
    \draw[ultra thick, teal] (5,2.5) -- (7.5,0);
    \node[text=teal] at (7,1.5) {$l_3$};
    \draw[ultra thick, black] (7.5,5) -- (10,2.5);
    \node at (8,3.5) {$l_4$};
    
    \draw (0,0) rectangle (10,10);
    \end{tikzpicture}
    };
    \node at (0,0) {
    \begin{tikzpicture}[scale=0.4]
    
    \draw[ultra thick, tomato] (5,7.5) -- (2.5,10);
    \draw[ultra thick, gray] (2.5,5) -- (0,7.5);
    \draw[ultra thick, teal] (10,2.5) -- (7.5,0);
    \draw[ultra thick, black] (7.5,5) -- (5,2.5);

    \draw (0,0) rectangle (10,10);

    \end{tikzpicture}
    };
    \node at (5,0) {
    \begin{tikzpicture}[scale=0.4]
    \draw[ultra thick, gray] (0,7.5) -- (2.5,10);
    \draw[ultra thick, tomato] (2.5,5) -- (5,7.5);
    \draw[ultra thick, black] (5,2.5) -- (7.5,0);
    \draw[ultra thick, teal] (7.5,5) -- (10,2.5);

    \draw (0,0) rectangle (10,10);
    \end{tikzpicture}
    };
    
\draw[->, thick] (-2.7,0) -- (-2.2,0); 
    \node[scale=1.5] at (-2.45,0.4) {$F$};
    \draw[->, thick] (2.2,0) -- (2.7,0); 
    \node[scale=1.5] at (2.45,0.4) {$G$};
    
    \end{tikzpicture}
    \caption{Line segments $l_j$ satisfying $H:l_1 \leftrightarrow l_2$, $l_3  \leftrightarrow l_4$. Each are periodic with period 2, their union is invariant under $H$.}
    \label{fig:periodic}
\end{figure}
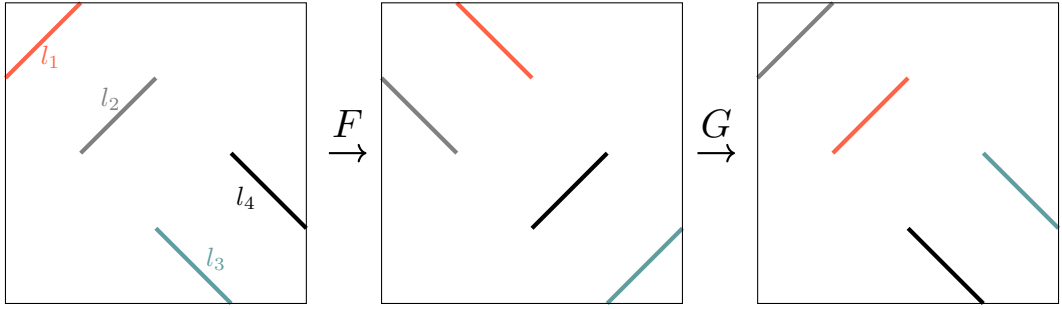

Other authors \cite{cheng_numerical_2023} have recognised the interest of this map, including it as part of a wider fundamental class of alternating wedge flows. We claim in particular it serves as a prototype example of piecewise linear maps which sit on the boundary of ergodicity. Limiting onto $H$ from its non-ergodic perturbations, the nature of the periodic orbits seeding the islands changes from elliptic to parabolic. Provided such an orbit does not limit onto a singularity line, its surrounding islands shrink, leaving behind periodic line segments ($H$ possesses four, sketched in Figure \ref{fig:periodic}) of null measure. This permits mixing with respect to Lebesgue but, as observed in \cite{cheng_numerical_2023}, only at a reduced polynomial rate for we can find orbits which `stick' to the segment for arbitrarily long periods. Here we show that $H$ mixes no slower than this. Our main theorem is the following:

\begin{thm}
    \label{thm:mainTheoremMixingRate}
    Correlations for $H$ decay as $|C_n(\varphi,\psi) | = \mathcal{O}(1/n)$ for H\"older observables $\varphi,\psi$.
\end{thm}

We expect a similar law to hold for piecewise linear systems obeying the limiting behaviour described above, for example the pointwise limit of $H_{(\xi,\eta)}$ as $\xi \to 0$ at $\eta=1/3$ (see \citealp{myers_hill_family_2022}). We focus on $H$ in particular for two key reasons. Firstly, as a fundamental piecewise linear model of alternating shear flows where no-slip boundary conditions force non-monotonic shear profiles, it is of interest to (laminar) fluid mixing applications \cite{cerbelli_continuous_2005}. Indeed, it is the logical extension to \citeauthor{cerbelli_continuous_2005}'s map, incorporating non-monotonicity into both the horizontal and vertical shears. Questions surrounding the mechanism by which $H$ is mixing, but at a reduced rate, are natural in this setting and are answered conclusively by a proof of Theorem \ref{thm:mainTheoremMixingRate}. Secondly, the map possesses certain properties which speed up its analysis. Since both $1/\eta$ and $1/(1-\eta)$ are integer valued over $0<\eta<1$ if and only if $\eta=1/2$, $H$ is the only map in the $0<\xi,\eta<1$ parameter space with all integer valued Jacobians and can be expressed as $H(x,y) = DH \cdot (x,y)^T$ mod 1\footnote{This also implies that periodic orbits are dense on $\tor$, as the cardinality of any orbit containing a rational point $(s/q,p/q) \in \tor$ with $s,p,q \in \mathbb{N}$ is bounded above by $q^2$ \cite{cerbelli_continuous_2005}}. This will prove useful for tracking the orbits of certain points under large powers of $H$. In addition $H$ can be related to its inverse by a conjugacy and behaves symmetrically on certain regions, reducing the calculations required to establish growth conditions by a factor of four.

The following sections are organised as follows. In section \ref{sec:outline} we state two theorems from the billiards literature that we rely upon to establish Theorem \ref{thm:mainTheoremMixingRate}. We next prove hyperbolicity for $H$ in section \ref{sec:OTMhyp} and the mixing property in section \ref{sec:Hmixing}. Central to this analysis is recurrence to a set $\sigma$ with the return map $H_\sigma$ exhibiting strong hyperbolic properties. We establish more formal properties of the map $H_\sigma$ in section \ref{sec:Hsigma}, sufficient to establish exponential decay of correlations. We use this to infer a polynomial bound on correlations for $H$ in section \ref{sec:polyMixingRate}, proving Theorem \ref{thm:mainTheoremMixingRate}. Finally in section \ref{sec:discussion} we comment on the relevance of our work to similar systems and suggest possible extensions.

\section{Some results from the billiards literature}
\label{sec:outline}
A necessary prerequisite for applying the machinery of \cite{chernov_billiards_2005} and similar is establishing mixing with respect to the invariant measure. In hyperbolic systems possessing singularities, the following scheme of \cite{katok_invariant_1986} is useful, giving conditions for the (stronger) Bernoulli property. We paraphrase from \cite{sturman_mathematical_2006}:

\begin{thm}[\citeauthor{katok_invariant_1986}]
\label{thm:katok-strelcyn}
Let $(X,\mathcal{F},\mu,f)$ be a measure preserving dynamical system such that $f$ is $C^2$ smooth outside of a singularity set $S$. Suppose that the Katok-Strelcyn conditions hold:
\begin{enumerate}[label={\bfseries (KS\arabic*):}]
    \item There exist $a,C_1>0$ such that for all $\epsilon>0$, $\mu(B_\varepsilon(S)) \leq C_1 \varepsilon^a$.
    \item There exist $b,C_2>0$ such that for all $x \in X \setminus S$, $||D^2_xf|| \leq C_2 \,d(x,S)^{-b}$.
    \item Lyapunov exponents exist and are non-zero almost everywhere.
\end{enumerate}
Then at almost every $x$ we can define local unstable and stable manifolds $\gamma_u(x)$ and $\gamma_s(x)$. Suppose that the manifold intersection property holds:
\begin{enumerate}[label={\bfseries (M):}]
    \item For almost any $x,x'\in X$, there exist $m,n$ such that $f^m(\gamma_u(x)) \cap f^{-n}(\gamma_s(x')) \neq \varnothing$.
\end{enumerate}
Then $f$ is ergodic. Provided the repeated manifold intersection property holds:
\begin{enumerate}[label={\bfseries (MR):}]
    \item For almost any $x,x'\in X$, there exist $M,N$ such that for all $m>M$ and $n>N$, $f^m(\gamma_u(x)) \cap f^{-n}(\gamma_s(x')) \neq \varnothing$,
\end{enumerate}
the Bernoulli property follows.
\end{thm}

The nature of the constant $a$ giving \textbf{(KS1)} plays an important role in showing expansion conditions such as (\ref{eq:oldOneStepExpansion}). In systems possessing a finite number of singularity curves, see for example \cite{przytycki_ergodicity_1983,myers_hill_exponential_2022,myers_hill_family_2022}, a covering by $\varepsilon$-balls immediately gives \textbf{(KS1)} with $a=1$. Showing (\ref{eq:oldOneStepExpansion}) is then quite straightforward; the singularity set splits an unstable manifold $W$ of vanishing length $|W| \to 0$ into at most $K$ components $W_k$, where $K$ is the maximum number of singularity curves which meet at a given point. This reduces (\ref{eq:oldOneStepExpansion}) to calculating the expansion factors $\lambda_k = |f(W_k)|/|W_k|$ and verifying the finite summation $\sum_k \lambda_k^{-1} <1 $. In many systems, in particular those driven by a return map where recurrence follows a law such as (\ref{eq:E1}), singularity curves instead form a \emph{countable} family. Expansion factors $\lambda_k \sim c\,k$ are typical so that bounding the above sum is challenging, indeed it may even diverge. Such systems satisfy \textbf{(KS1)}, but only with some $a<1$. In certain $a<1$ scenarios, precise mapping behaviour may reduce (\ref{eq:oldOneStepExpansion}) to a finite summation; see for example the return map considered in \cite{springham_polynomial_2014}. Such a scenario is not typical however, with (\ref{eq:oldOneStepExpansion}) failing in many examples \cite{chernov_billiards_2005}. More recent schemes for bounds on correlations have revised (\ref{eq:oldOneStepExpansion}) to suit these more general $a < 1$ systems. We quote the first of these, given in \cite{chernov_statistical_2009}, which is sufficient for our purposes.

Let $\Omega$ denote a two dimensional connected compact Riemannian manifold, $f: \Omega \to \Omega$ preserving a measure $\mu$. Let $d$ denote the distance in $\Omega$ induced by the Riemannian metric $\rho$. For any smooth curve $W$ in $\Omega$, denote by $|W|$ its length, and by $m_W$ the Lebesgue measure on $W$ induced by the Riemannian metric $\rho_W$ restricted to $W$. Also let $\nu_W$ = $m_W /|W|$ be the normalised (probability) measure on W.

\textbf{(H1):} Hyperbolicity of $f$ (with uniform expansion and contraction). There exist two families of cones $C_x^u$ (unstable) and $C_x^s$ (stable) in the tangent spaces $\mathcal{T}_x\Omega$, for all $x \in \Omega$, and there exists a constant $\Lambda$ > 1, with the following properties:
\begin{enumerate}
    \item $Df(C_x^u) \subset C_{fx}^u$ and $Df(C_x^s) \supset C_{fx}^s$ whenever $Df$ exists.
    \item $\| D_xf(v) \| \geq \Lambda \| v \|$ for all  $v \in C_x^u$ and $\| D_xf^{-1}(v) \| \geq \Lambda \| v \|$ for all $v \in C_x^s$.
    \item These families of cones are continuous on $\Omega$ and the angle between $C_x^u$ and $C_x^s$ is uniformly bounded away from zero.
\end{enumerate}

We say that a smooth curve $W \subset \Omega$ is an unstable (stable) curve if at every point $x \in W$
the tangent line $\mathcal{T}_x W$ belongs in the unstable (stable) cone $C_x^u$ ($C_x^s$).

\textbf{(H2):} Singularities and smoothness. Let $\mathcal{S}_0$ be a closed subset in $\Omega$, such that $M := \Omega \setminus \mathcal{S}_0$ is a dense set in $\Omega$. We put $\mathcal{S}_{\pm 1} = f^\mp \mathcal{S}_0$.
\begin{enumerate}
    \item $f:M\setminus \mathcal{S}_1 \to M \setminus \mathcal{S}_{-1}$ is a $C^2$ diffeomorphism.
    \item $\mathcal{S}_0 \cup \mathcal{S}_1$ is a finite or countable union of smooth, compact curves in $\Omega$.
    \item Curves in $\mathcal{S}_0$ are transversal to stable and unstable cones. Every smooth curve in $\mathcal{S}_1$ (resp. $\mathcal{S}_{-1}$) is a stable (resp. unstable) curve. Every curve in $\mathcal{S}_1$ terminates either inside another curve of $\mathcal{S}_1$ or on $\mathcal{S}_0$.
    \item There exists $b \in (0,1)$ and $c>0$ such that for any $x\in M\setminus \mathcal{S}_1$
    \begin{equation}
        \label{eq:DxbCondition}
        \| D_xf  \| \leq c\, d(x,\mathcal{S}_1)^{-b}.
    \end{equation}
\end{enumerate}

\textbf{(H3):} Regularity of smooth unstable curves. We assume that there is a $f$-invariant class of unstable curves $W \subset M$ that are \emph{regular} (see \citealp{chernov_statistical_2009}).

\textbf{(H4):} SRB measure. $\mu$ is a Sinai-Ruelle-Bowen (SRB) measure which is mixing.

\textbf{(H5):} One-step expansion. There exists $q \in (0,1]$ such that
\begin{equation}
    \label{eq:oneStep}
    \liminf_{\delta \to 0} \sup_{W: |W|< \delta} \sum_i \left( \frac{|W|}{|f(W_i)|}\right)^q \frac{|W_i|}{|W|} < 1,
\end{equation}
where the supremum is taken over all unstable curves, $W_i$ are the components of $W$ split by the singularity set for $f$.

\begin{thm}[\citeauthor{chernov_statistical_2009}]
\label{thm:chernovZhang}
Under the conditions \textbf{(H1)}–\textbf{(H5)}, the system $(f, \mu)$ enjoys exponential decay of correlations.
\end{thm}

Note that the new one-step expansion condition (\ref{eq:oneStep}) may be reduced to the old (\ref{eq:oldOneStepExpansion}) by taking $q=1$. The new condition ensures that the images of unstable curves grow `on average'. Choosing a $q<1$ essentially permits summing over countably many components, broadening the potential applications of the scheme to a wider class of $a<1$ systems. The image coupling methods (\cite{young_recurrence_1999}, see also \cite{chernov_chaotic_2006} and the references therein) used to establish Theorem \ref{thm:chernovZhang} differ substantially from those employed in \cite{chernov_billiards_2005}. The key `magnet' construction \cite{chernov_advanced_2006,chernov_chaotic_2006}, however, further serves as the base $\Delta_0$ of a Young tower satisfying the exponential tail bound (\ref{eq:expRtailYoung}) \cite{chernov_statistical_2009}. As such the scheme may similarly be applied to some return map $f_M$ as a step towards proving polynomial decay of correlations for $f$. We conclude this section with two technical adjustments we will refer back to later in section \ref{sec:Hsigma}.

\begin{remark}
\label{remark:firstRemark}
Condition \textbf{(H1.3)} has been relaxed in subsequent schemes \cite{demers_spectral_2014,wang_decay_2021} and can be replaced by
\begin{enumerate}[label=3'.]
    \item These families of cones are continuous on components of $\Omega \setminus \mathcal{S}_0$ and the angle between $C_x^u$ and $C_x^s$ is uniformly bounded away from zero.
\end{enumerate}
Theorem \ref{thm:chernovZhang} still follows under this relaxed assumption by applying (for example) Theorem 1 of \cite{wang_decay_2021}. Despite the improvement over older growth conditions, condition \textbf{(H5)} still fails for many systems over one iterate. See, for example, the modified stadia considered in \cite{chernov_statistical_2009}. It can be replaced by a multi-step expansion condition, establishing \textbf{(H5)} for some higher power $f^n$ of the map and its enlarged singularity set.
\end{remark}

\section{Hyperbolicity}
\label{sec:OTMhyp}
\begin{prop}
\label{prop:tentMap}
$H$ is non-uniformly hyperbolic. That is, Lyapunov exponents
\[\chi (z,v) = \limn \frac{1}{n} \log||DH^n_z v|| \]
are non-zero for almost every $z \in \tor$ and tangent vector $v \neq 0$.
\end{prop}

The key ingredients of the proof were sketched out in \cite{myers_hill_exponential_2022}. We provide a more detailed treatment here as certain constructions are central to the analysis of later sections. We begin with a description of the Jacobian $DH$ and recall a decomposition of the cocycle $DH_z^n$ into blocks (Lemma \ref{lemma:itineraries}) which share an invariant expanding cone (Lemma \ref{lemma:cone}). Associating these blocks with recurrence to a region $\sigma$ then allows us to deduce non-zero Lyapunov exponents on a full measure set.

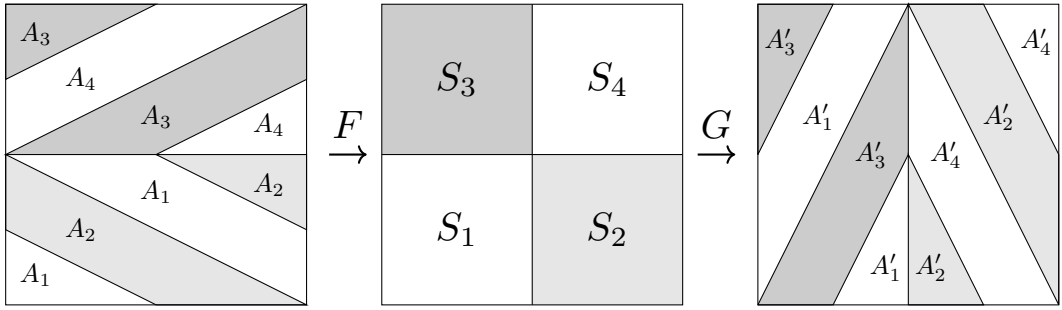
\begin{figure}
    \centering
 \begin{tikzpicture}
     \tikzmath{\e = 0.5;} 
    \node[scale=1] at (0,0) {
    \begin{tikzpicture}
    
    \fill[gray!40] (0,4-4*\e) rectangle (4-4*\e,4);
    
    \fill[gray!20] (4-4*\e,0) rectangle (4,4-4*\e);
    
    \draw (0,0) rectangle (4,4);
    \draw (4-4*\e,0) -- (4-4*\e,4);
    \draw (0,4-4*\e) -- (4,4-4*\e);

    
    \node[scale=1.5] at (2-0.5*4*\e,2-0.5*4*\e) {$S_1$};
    \node[scale=1.5] at (4-0.5*4*\e,2-0.5*4*\e) {$S_2$};
    \node[scale=1.5] at (2-0.5*4*\e,4-0.5*4*\e) {$S_3$};
    \node[scale=1.5] at (4-0.5*4*\e,4-0.5*4*\e) {$S_4$};
    
    \end{tikzpicture}
    };
    \node[scale=1] at (-5,0) {
    \begin{tikzpicture}
    
    \filldraw[fill=gray!40] (0, {4*(1-\e+\e*\e)}) -- ({4*(1-\e)},4) -- (0,4) --  (0, {4*(1-\e+\e*\e)});
    
    \filldraw[fill=gray!40] (4,4) -- (0,{4*(1-\e)}) -- ({4*(1-\e)},{4*(1-\e)}) -- (4, {4*(1-\e+\e*\e)}) -- (4,4);
    
    \filldraw[fill=gray!20] ({4*(1-\e)},{4*(1-\e)}) -- (4,{4*(1-\e)}) -- (4, {4*(1-\e)^2}) -- ({4*(1-\e)},{4*(1-\e)});
    
    \filldraw[fill=gray!20] (0,{4*(1-\e)}) -- (4,0) -- ({4*(1-\e)},0) -- (0, {4*(1-\e)^2}) -- (0,{4*(1-\e)});
    
    \draw (0,0) rectangle (4,4);

    
        
    \node at (2-0.5*4*\e,2-0.5*4*\e) {$A_2$};
    \node at (3.5,1.6) {$A_2$};
    \node at (3.5,2.4) {$A_4$};
    
    \node at (2,1.5) {$A_1$};
    \node at (2-0.5*4*\e,2+0.5*4*\e) {$A_4$};

    \node at (2,2.5) {$A_3$};
    \node at (0.4,3.6) {$A_3$};
    \node at (0.4,0.4) {$A_1$};
    \end{tikzpicture}
    };
    
    \node[scale=1] at (5,0) {
    \begin{tikzpicture}
    
    \filldraw[fill = gray!40] (0,4-4*\e) -- ({4*\e*(1-\e)},4) -- (0,4) -- (0,4-4*\e);
    \filldraw[fill = gray!40] (0,0) -- (4-4*\e,4) -- (4-4*\e,4-4*\e) -- ({4*\e*(1-\e)},0) -- (0,0);
    
    \filldraw[fill=gray!20] (4-4*\e,0) -- (4-4*\e,4-4*\e) -- (4-4*\e*\e,0) -- (4-4*\e,0);
    \filldraw[fill=gray!20] (4,0) -- (4-4*\e,4) -- (4-4*\e*\e,4)-- (4,4-4*\e) -- (4,0);
    
    \draw (0,0) rectangle (4,4);
    
    

    \node at (1.5,2) {$A_3'$};
    \node at (1.7,0.5) {$A_1'$};
    \node at ({0.4*4*(1-\e)},2.5) {$A_1'$};  
    \node at (0.3,3.5) {$A_3'$};

    \node at (4-1.5,2) {$A_4'$};
    \node at (4-1.7,0.5) {$A_2'$};
    \node at ({4-0.4*4*(1-\e)},2.5) {$A_2'$};  
    \node at (4-0.3,3.5) {$A_4'$};

    \end{tikzpicture}
    };
    
\draw[->, thick] (-2.7,0) -- (-2.2,0); 
    \node[scale=1.5] at (-2.45,0.4) {$F$};
    \draw[->, thick] (2.2,0) -- (2.7,0); 
    \node[scale=1.5] at (2.45,0.4) {$G$};
    
    \end{tikzpicture}
    \caption{A partition of the torus into four rectangles $S_j$, and their preimages $A_j$, $A_j'$ under $F,G^{-1}$.}
        \label{fig:firstPartitions}
\end{figure}

Partition the torus into the four squares $S_j$ shown in Figure \ref{fig:firstPartitions}. The Jacobian $DH$ is then constant on the preimages $A_j=F^{-1}(S_j)$, given by the matrix $M_j$ where
\[ M_1 = \begin{pmatrix} 1 & 2 \\ 2 & 5\end{pmatrix}, \quad M_2 = \begin{pmatrix} 1 & 2 \\ -2 & -3\end{pmatrix}, \quad M_3 = \begin{pmatrix} 1 & -2 \\ 2 & -3\end{pmatrix}, \quad M_4 = \begin{pmatrix} 1 & -2 \\ -2 & 5\end{pmatrix}, \]
undefined on the singularity set $\mathcal{D} = \cup_j \partial A_j$. Letting $X'$ denote the full measure set $\tor \setminus \cup_{i \geq 0} H^{-i}(\mathcal{D})$, the $n$-step itinerary
\[  A_{j_1}, A_{j_2}, A_{j_3}, \dots, A_{j_n},\]
is well defined for any $z \in X'$. The related cocycle $DH_z^n$ given by
\[ DH_z^n = M_{j_n} \dots M_{j_3} M_{j_2}  M_{j_1}\]
with each $j_k \in \{ 1,2,3,4 \}$. Our aim is to decompose any cocycle into hyperbolic matrices which share an invariant expanding cone.
Note that while $M_1$ and $M_4$ are hyperbolic, $M_2$ and $M_3$ are not. Hence when $M_2$ or $M_3$ appear in a cocycle at $M_{j_k}$, we must combine them with its neighbouring matrices $M_{j_{k+l}},\dots,M_{j_{k+2}}, M_{j_{k+1}}$ for some $l \in \mathbb{N}$.

Consider the countable family of matrices 
\[\mathcal{M} = \{ M_1, M_4,  M_1M_2^n, M_3M_2^n, M_4M_2^n, M_1M_3^n, M_2M_3^n, M_4M_3^n\} \]
with $n\in\mathbb{N}$. Similarly define 
\[\mathcal{M}' = \{ M_1^{-1}, M_4^{-1},  M_1^{-1}M_2^{-n}, \dots, M_4^{-1}M_3^{-n}\}.\]
It was shown in \cite{myers_hill_exponential_2022} that:

\begin{lemma}
\label{lemma:itineraries}
At almost every $z$, the cocycle $DH_z^n$ can be decomposed into blocks from $\mathcal{M}$.
\end{lemma}

The result essentially follows from the fact that essentially no orbits get trapped in $A_3$,
\begin{equation}
    \label{eq:escapeA3statement}
    \limn \mu\left( \{ z \in A_3 \, | \, H^i(z) \in A_3 \text{ for all } 0 \leq i \leq n-1  \}\right) = 0,
\end{equation}
and the equivalent statement for $A_2$. An entirely analogous argument, considering escapes from $A_2'$ and $A_3'$ under $H^{-1}$, gives that at a.e. $z$ the cocycle $DH_z^{-n}$ can be decomposed into blocks from $\mathcal{M}'$. 

\begin{lemma}
\label{lemma:cone}
The matrices in $\mathcal{M}$ admit an invariant expanding cone $\mathcal{C}$.
\end{lemma}

\begin{proof}
Parameterise the tangent space by $(v_1,v_2)^T\in \mathbb{R}^2$. The lemma was shown in \cite{myers_hill_exponential_2022} using the cone $C = \{(v_1,v_2)\neq 0 \,|\, |v_2| \geq \phi \, |v_1|\}$ where $\phi$ is the golden ratio $(1+\sqrt{5})/2$. Here we define a slightly wider cone $\mathcal{C} = \{(v_1,v_2)\neq 0 \,|\, |v_2| \geq  \varphi \, |v_1|\}$, $\varphi = 21/13$, which still contains all the unstable eigenvectors of matrices in $\mathcal{M}$ and none of the stable eigenvectors. Hence $\mathcal{C}$ is invariant and one can verify that it is also expanding (minimum expansion factors are calculated later in Table \ref{tab:tab2}, in particular the minimum expansion of a matrix $M$ over $\mathcal{C}$ under the $\| \cdot\|_\infty$ norm is given by $\min_\pm K_\pm(M)$). 
\end{proof}

\subsubsection*{Recurrence to $\sigma$}

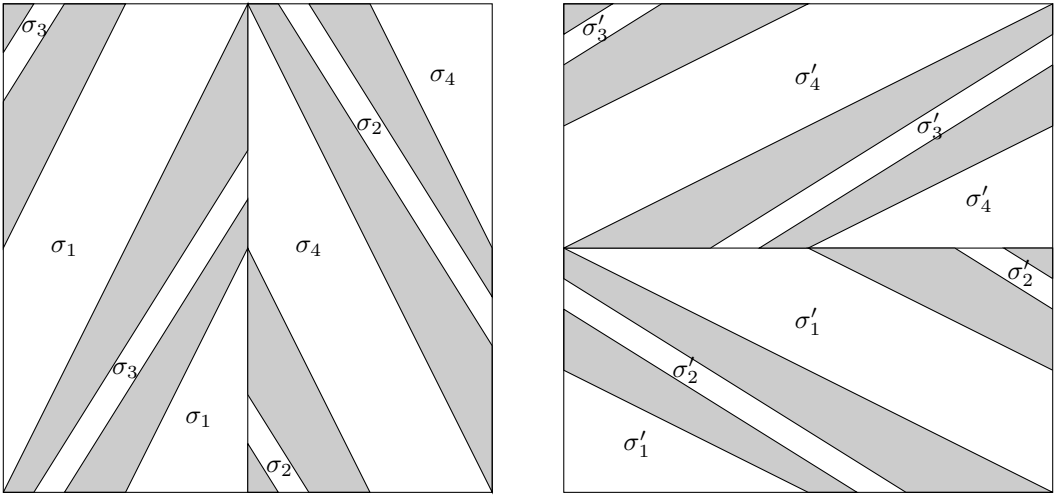
\begin{figure}
    \centering
\begin{tikzpicture}[scale=0.65]
    
    \draw (5,0) -- (5,10);
    \filldraw[fill=gray!40] (0,8) -- (10/8,10) -- (2.5,10) -- (0,5) -- (0,8);
      \filldraw[fill=gray!40] (0,9) -- (0,10) -- (10/16,10) -- (0,9);
      \filldraw[fill=gray!40] (0,0) -- (5,10) -- (5,7) -- (10/16,0) -- (0,0);
      \filldraw[fill=gray!40] (10/8,0) -- (2.5,0) -- (5,5) -- (5,6) -- (10/8,0) ;
      \filldraw[fill=gray!40] (5,2) -- (50/8,0) -- (7.5,0) -- (5,5) -- (5,2);
      \filldraw[fill=gray!40] (5,1) -- (5,0) -- (90/16,0) -- (5,1);
      \filldraw[fill=gray!40] (5,10) -- (10,0) -- (10,3) -- (90/16,10) -- (5,10);
      \filldraw[fill=gray!40] (50/8,10) -- (7.5,10) -- (10,5) -- (10,4) -- (50/8,10);
 
     \node at (0.65,9.5) {$\sigma_3$};
     \node at (1.25,5) {$\sigma_1$};
     \node at (2.5,2.5) {$\sigma_3$};
     \node at (4,1.5) {$\sigma_1$};
     
     \node at (0.65+5,10-9.5) {$\sigma_2$};
     \node at (1.25+5,5) {$\sigma_4$};
     \node at (2.5+5,10-2.5) {$\sigma_2$};
    \node  at (4+5,10-1.5) {$\sigma_4$};
     \draw(0,0) rectangle (10,10);
    \end{tikzpicture}
    \hspace{2em}
    \begin{tikzpicture}[scale=0.65]
    
    \draw (0,5) -- (10,5);
    \filldraw[fill=gray!40] (0,150/16) -- (0,10) -- (1,10) -- (0,150/16);
      \filldraw[fill=gray!40] (0,70/8) -- (2,10) -- (5,10) -- (0,7.5);
      \filldraw[fill=gray!40] (0,5) -- (10,10) -- (10,150/16) -- (3,5) -- (0,5);
      \filldraw[fill=gray!40] (10,70/8) -- (10,7.5) -- (5,5) -- (4,5) -- (10,70/8) ;
   
     \filldraw[fill=gray!40] (10,150/16-5) -- (10,10-5) -- (10-1,10-5) -- (10,150/16-5);
      \filldraw[fill=gray!40] (10,70/8-5) -- (10-2,10-5) -- (10-5,10-5) -- (10,7.5-5);
      \filldraw[fill=gray!40] (10,5-5) -- (0,10-5) -- (0,150/16-5) -- (10-3,5-5) -- (10,5-5);
      \filldraw[fill=gray!40] (0,70/8-5) -- (0,7.5-5) -- (10-5,5-5) -- (10-4,5-5) -- (0,70/8-5) ;   
 
     \node at (0.65,9.5) {$\sigma_3'$};
     \node at (5,8.5) {$\sigma_4'$};
     \node at (7.5,7.5) {$\sigma_3'$};
     \node at (8.5,6) {$\sigma_4'$};
     
     \node at (10-0.65,9.5-5) {$\sigma_2'$};
     \node at (10-5,8.5-5) {$\sigma_1'$};
     \node at (10-7.5,7.5-5) {$\sigma_2'$};
     \node at (10-8.5,6-5) {$\sigma_1'$};

     \draw(0,0) rectangle (10,10);
    \end{tikzpicture}
   
     \caption{Partitions of the return sets $\sigma,\sigma'$ (white) for $H,H^{-1}$ into four sets $\sigma_j \subset H(A_j)$, $\sigma_j' \subset H^{-1}(A_j')$.  \label{fig:sigmaPartition}}
\end{figure}

Define $\sigma$ as the union of the sets $\sigma_1 = H(A_1)$, $\sigma_2 = H(A_2 \cap H(A_3))$, $\sigma_3 = H(A_3 \cap H(A_2)$, $\sigma_4 = H(A_4)$. By construction, any orbit escaping $A_2,A_3$ or passing through $A_1,A_4$ must pass through $\sigma$. The return map $H_\sigma : \sigma \to \sigma$, $z \mapsto H^R(z)$ where $R = R(z;H,\sigma)$, is well defined at $\mu$-almost every $z \in \sigma$ by (\ref{eq:escapeA3statement}) and the equivalent statement for $A_2$. We similarly define $\sigma' = \cup_j \sigma_j'$ using the $A_j'$ and the return map $H^{-1}_{\sigma'}:\sigma' \rightarrow \sigma'$ for $H^{-1}$. The sets $\sigma,\sigma'$ are shown as the unshaded regions in Figure \ref{fig:sigmaPartition}.

We begin by identifying the points in $\sigma$ with return time 1, i.e. $H^{-1}(\sigma) \cap \sigma$. The preimages of $\sigma_1, \sigma_4$ are simply $A_1$, $A_4$ and by definition we have $H^{-1}(\sigma_2) = A_2 \cap H(A_3)$ so that $H^{-1}(\sigma_2) \cap \sigma = A_2 \cap \sigma_3 := \varsigma_3$ and similarly $H^{-1}(\sigma_3) \cap \sigma = A_3 \cap \sigma_2 := \varsigma_2$. See Figure \ref{fig:returnTime1} for an illustration. 

\begin{figure}
    \centering
    \subfigure[][]{
    \label{fig:returnTime1}
    
    \definecolor{tomato}{RGB}{255, 99, 71}
    \definecolor{teal}{RGB}{95, 158, 160} 
    
     \begin{tikzpicture}[scale=0.7]
    
    \clip (0,0) rectangle (10,10);
    
    \draw (0,5) -- (10,5);
    
    \draw (0,10) -- (10,10);
    \draw (5,0) -- (10,0);

    \filldraw[fill=tomato!80] (0,5) -- (5,5) -- (10,2.5) -- (10,0) -- (0,5);
    \filldraw[fill=tomato!80] (0,2.5) -- (5,0) -- (0,0) -- (0,2.5);
    
    \filldraw[fill=tomato!80] (0,5) -- (10,10) -- (5,10) -- (0,7.5) -- (0,5);
    \filldraw[fill=tomato!80] (5,5) -- (10,5) -- (10,7.5) -- (5,5);

    \begin{scope}
    \clip (0,0) -- (5,10) -- (5,5) -- (2.5,0) -- (0,0);
    \filldraw[fill=teal!80] (0,5) -- (10,0) -- (5,0) -- (0,2.5) -- (0,5);
    \end{scope}
    
    \begin{scope}
    \clip (10,0) -- (5,10) -- (7.5,10) -- (10,5) -- (10,0);
    \filldraw[fill=teal!80] (0,5) -- (10,10) -- (10,7.5) -- (5,5) -- (0,5);
    \end{scope}

      \filldraw[fill=gray!80] (0,8) -- (10/8,10) -- (2.5,10) -- (0,5) -- (0,8);
      \filldraw[fill=gray!80] (0,9) -- (0,10) -- (10/16,10) -- (0,9);
      \filldraw[fill=gray!80] (0,0) -- (5,10) -- (5,7) -- (10/16,0) -- (0,0);
      \filldraw[fill=gray!80] (10/8,0) -- (2.5,0) -- (5,5) -- (5,6) -- (10/8,0) ;
      
      \filldraw[fill=gray!80] (5,2) -- (50/8,0) -- (7.5,0) -- (5,5) -- (5,2);
      \filldraw[fill=gray!80] (5,1) -- (5,0) -- (90/16,0) -- (5,1);
      \filldraw[fill=gray!80] (5,10) -- (10,0) -- (10,3) -- (90/16,10) -- (5,10);
      \filldraw[fill=gray!80] (50/8,10) -- (7.5,10) -- (10,5) -- (10,4) -- (50/8,10);

    \draw (0,0) -- (0,3);
    \draw (0,4) -- (0,5);
    
    \draw (10,5) -- (10,8);
    \draw (10,9) -- (10,10);
    \node at (2.5,2.5) {$\varsigma_3$};
    \node at (7.5,7.5) {$\varsigma_2$};
    \end{tikzpicture}

    }
         \subfigure[][]{
    \label{fig:HSigma}
    
    \begin{tikzpicture}[scale=0.7]
    
    \clip (0,0) rectangle (10,10);
    
    \draw (0,5) -- (10,5);
    
    \draw (0,10) -- (10,10);
    \draw (5,0) -- (10,0);

        \foreach \k in {1,...,25}{
        \draw ({ (2*\k+2)*(7.5-7.5)/(2*\k+1)  }, { 7.5 })--({ (2*\k+2)*(10-7.5)/(2*\k+1)  }, { 10 });
        \draw ({ 5+ (2*\k+2)*(5-7.5)/(2*\k+1)  }, { 5 })--({ 5+ (2*\k+2)*(7.5-7.5)/(2*\k+1) }, { 7.5 });
        \draw({ 2.5 + (2*\k+1)*(5-5)/(2*\k)  }, { 5 })--({ 2.5 + (2*\k+1)*(10-5)/(2*\k)  }, { 10 });
        \draw ({ 12.5 + (2*\k+1)*(5-10)/(2*\k)  }, { 5 }) -- ({ 12.5 + (2*\k+1)*((30*\k+20)/(4*\k+2)-10)/(2*\k)  }, { (30*\k+20)/(4*\k+2) });
        }
        
        \draw (30/4,70/8) -- (4,6);
        \draw (70/8,150/16) -- (7,8);
        \draw (15-30/4,15-70/8) -- (15-5,15-190/28);
        \draw (15-70/8,15-150/16) -- (15-7,15-8);
        \draw (2.5,10) -- (0 , 50/6);
        \draw (1,9) -- (0,15-190/28);

        \draw (5,7)--(90/14,230/28);
        \draw (15-5,15-7)--(15-90/14,15-230/28);
        
        \filldraw[fill=white] (0,5) -- (10,10) -- (5,10) -- (0,7.5) -- (0,5);
        \filldraw[fill=white] (5,5) -- (10,7.5) -- (10,5) -- (5,5);

        \draw (0,150/16) -- (1,10);
        \draw (0,70/8) -- (2,10);
        \draw (10,150/16) -- (3,5);
        \draw (10,70/8) -- (4,5);
      
       \foreach \k in {1,...,25}{
        \draw ({10- (2*\k+2)*(7.5-7.5)/(2*\k+1)  }, { 7.5-5 })--({10- (2*\k+2)*(10-7.5)/(2*\k+1)  }, { 10-5 });
        \draw ({ 5 - (2*\k+2)*(5-7.5)/(2*\k+1)  }, { 5-5 })--({ 5 - (2*\k+2)*(7.5-7.5)/(2*\k+1)  }, { 7.5-5 });
        \draw ({ 7.5 - (2*\k+1)*(5-5)/(2*\k)  }, { 5-5 })--({ 7.5- (2*\k+1)*(10-5)/(2*\k)  }, { 10-5 });
        \draw ({ -2.5 - (2*\k+1)*(5-10)/(2*\k)  }, { 5-5 }) -- ({ -2.5 - (2*\k+1)*((30*\k+20)/(4*\k+2)-10)/(2*\k)  }, { (30*\k+20)/(4*\k+2) -5 });
        }
        
        \draw (10-30/4,70/8-5) -- (10-4,6-5);
        \draw (10-70/8,150/16-5) -- (10-7,8-5);
        \draw (-5+30/4,15-70/8-5) -- (0,15-190/28-5);
        \draw (-5+70/8,15-150/16-5) -- (-5+7,15-8-5);
        \draw (10-2.5,10-5) -- (10, 50/6-5);
        \draw (10-1,9-5) -- (10,15-190/28-5);

        \draw (10-5,7-5)--(10-90/14,230/28-5);
        \draw (0,15-7-5)--(-5+90/14,15-230/28-5);
        
        \filldraw[fill=white] (10,0) -- (0,5) -- (5,5) -- (10,2.5) -- (10,0);
        \filldraw[fill=white] (0,0) -- (5,0) -- (0,2.5) -- (0,0);

        \draw (10,150/16-5) -- (10-1,10-5);
        \draw (10,70/8-5) -- (10-2,10-5);
        \draw (0,150/16-5) -- (10-3,5-5);
        \draw (0,70/8-5) -- (10-4,5-5);
       
   \begin{scope}
    \clip (0,0) -- (5,10) -- (5,5) -- (2.5,0) -- (0,0);
    \filldraw[fill=white] (0,5) -- (10,0) -- (5,0) -- (0,2.5) -- (0,5);
    \end{scope}
    
    \begin{scope}
    \clip (10,0) -- (5,10) -- (7.5,10) -- (10,5) -- (10,0);
    \filldraw[fill=white] (0,5) -- (10,10) -- (10,7.5) -- (5,5) -- (0,5);
    \end{scope}

      \filldraw[fill=gray!80] (0,8) -- (10/8,10) -- (2.5,10) -- (0,5) -- (0,8);
      \filldraw[fill=gray!80] (0,9) -- (0,10) -- (10/16,10) -- (0,9);
      \filldraw[fill=gray!80] (0,0) -- (5,10) -- (5,7) -- (10/16,0) -- (0,0);
      \filldraw[fill=gray!80] (10/8,0) -- (2.5,0) -- (5,5) -- (5,6) -- (10/8,0) ;
      
      \filldraw[fill=gray!80] (5,2) -- (50/8,0) -- (7.5,0) -- (5,5) -- (5,2);
      \filldraw[fill=gray!80] (5,1) -- (5,0) -- (90/16,0) -- (5,1);
      \filldraw[fill=gray!80] (5,10) -- (10,0) -- (10,3) -- (90/16,10) -- (5,10);
      \filldraw[fill=gray!80] (50/8,10) -- (7.5,10) -- (10,5) -- (10,4) -- (50/8,10);

    \draw (0,0) -- (0,3);
    \draw (0,4) -- (0,5);
    
    \draw (10,5) -- (10,8);
    \draw (10,9) -- (10,10);

    \draw[red, dashed, very thick] (0,3)--(0,4);
    \draw[red, dashed, very thick] (0,8)--(0,9);
    \draw[red, dashed, very thick] (10,3)--(10,4);
    \draw[red, dashed, very thick] (10,8)--(10,9);

    \draw[red, dashed, thick] (5,1) -- (5,2);
    \draw[red, dashed, thick] (5,6) -- (5,7);

    \end{tikzpicture}
    }
    
    \caption{Part (a) shows the portions of $\sigma$ (red, blue) with return time 1 to $\sigma$. Points in the white region have return times of 2 or more. Part (b) shows the singularity set $\mathcal{S}$ for the return map $H_\sigma$. Red dashed lines denote the shared boundaries of the $\sigma_j$.}
    \label{fig:sigmaSingSet}
\end{figure}
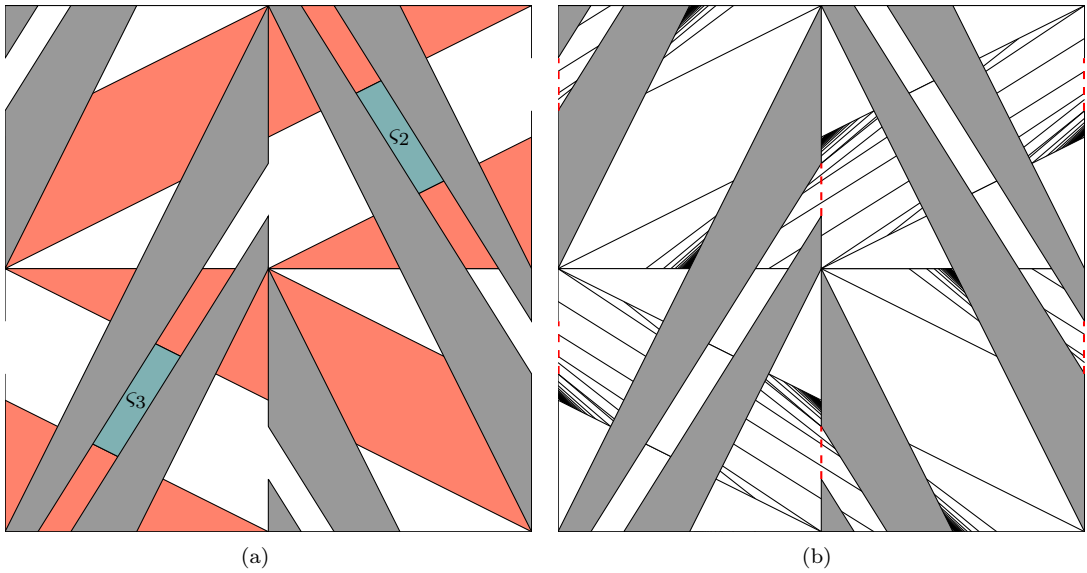

\begin{figure}
    \centering
    
         \subfigure[][]{
    \label{fig:escapeTimes}

    \begin{tikzpicture}[scale=1.4]
    
    \clip (0,5) rectangle (10,10);
    
    
        
        \foreach \k in {1,...,45}{
        \draw[very thin] ({ (2*\k+2)*(7.5-7.5)/(2*\k+1)  }, { 7.5 })--({ (2*\k+2)*(10-7.5)/(2*\k+1)  }, { 10 });
        \draw[very thin] ({ 5+ (2*\k+2)*(5-7.5)/(2*\k+1)  }, { 5 })--({ 5+ (2*\k+2)*(7.5-7.5)/(2*\k+1) }, { 7.5 });
        \draw[very thin] ({ 2.5 + (2*\k+1)*(5-5)/(2*\k)  }, { 5 })--({ 2.5 + (2*\k+1)*(10-5)/(2*\k)  }, { 10 });
        
         \draw[very thin] ({ 2.5  }, { 10 })--({ 2.5 - (2*\k+1)*(10-5)/(2*\k)  }, { 5 });
        
        \draw[very thin] ({ 12.5 + (2*\k+1)*(5-10)/(2*\k)  }, { 5 }) -- ({ 12.5 + (2*\k+1)*((30*\k+20)/(4*\k+2)-10)/(2*\k)  }, { (30*\k+20)/(4*\k+2) });
        }

        \fill[fill=white] (0,5) -- (10,10) -- (5,10) -- (0,7.5);
        \fill[fill=white] (10,5) -- (5,5) -- (10,7.5);
        
        \draw[thick] (0,10) -- (5,10) -- (0,7.5);
        \draw[thick] (0,5) -- (10,10);
        \draw[thick] (0,5) -- (5,5) -- (10,7.5);

        \node at (7.5,7.5) {$A^1$};
        \node at (6.6,8) {$A^2$};
        \node[scale=0.7] at (6,7.9) {$3$};
        \node[scale=0.5] at (5.7,7.8) {$4$};
        
         \node at (0.5,9.5) {$A^1$};
        \node at (0.12,8.3) {$A^2$};
        \node[scale=0.7] at (0.1,8) {$3$};
        
        \node at (15-6.6,15-8) {$A^2$};
        \node[scale=0.7] at (15-6,15-7.9) {$3$};
        \node[scale=0.5] at (15-5.7,15-7.8) {$4$};

        \node at (3.5,9.6) {$A^1$};
        \node at (2.9,9.8) {$A^2$};
        \node[scale=0.7] at (2.8,9.9) {$3$};

        \node at (5-3.5,15-9.6) {$A^1$};
        \node at (5-2.95,15-9.8) {$A^2$};
        \node[scale=0.7] at (5-2.8,15-9.9) {$3$};
        
        \draw[very thick] (0,7.5) -- (2.5,10);
        \draw[very thick] (2.5,5) -- (5,7.5);

    \end{tikzpicture}
    }
         \subfigure[][]{
    \label{fig:escapeBehaviour}
    
    \begin{tikzpicture}[scale=1.4]
    
    \clip (0,5) rectangle (10,10);
        
        \foreach \k in {1,...,45}{
                \draw[very thin] ({ (2*\k+2)*(7.5-7.5)/(2*\k+1)  }, { 7.5 })--({ (2*\k+2)*(10-7.5)/(2*\k+1)  }, { 10 });
        \draw[very thin] ({ 5+ (2*\k+2)*(5-7.5)/(2*\k+1)  }, { 5 })--({ 5+ (2*\k+2)*(7.5-7.5)/(2*\k+1) }, { 7.5 });
        \draw[very thin] ({ 2.5 + (2*\k+1)*(5-5)/(2*\k)  }, { 5 })--({ 2.5 + (2*\k+1)*(10-5)/(2*\k)  }, { 10 });
        
         \draw[very thin] ({ 2.5  }, { 10 })--({ 2.5 - (2*\k+1)*(10-5)/(2*\k)  }, { 5 });
        
        \draw[very thin] ({ 12.5 + (2*\k+1)*(5-10)/(2*\k)  }, { 5 }) -- ({ 12.5 + (2*\k+1)*((30*\k+20)/(4*\k+2)-10)/(2*\k)  }, { (30*\k+20)/(4*\k+2) });
        }
        
        \draw[red] (30/4,70/8) -- (4,6);
        \draw[red] (70/8,150/16) -- (7,8);
        \draw[red] (15-30/4,15-70/8) -- (15-5,15-190/28);
        \draw[red] (15-70/8,15-150/16) -- (15-7,15-8);

        \draw[red] (1,9) -- (0,15-190/28);

        \draw[red] (5,7)--(90/14,230/28);
        \draw[red] (15-5,15-7)--(15-90/14,15-230/28);

        \fill[fill=white] (0,5) -- (10,10) -- (5,10) -- (0,7.5);
        \fill[fill=white] (10,5) -- (5,5) -- (10,7.5);
        
        \draw[thick] (0,10) -- (5,10) -- (0,7.5);
        \draw[thick] (0,5) -- (10,10);
        \draw[thick] (0,5) -- (5,5) -- (10,7.5);

        \draw[red] (0,150/16) -- (1,10);
        \draw[red] (0,70/8) -- (2,10);
        \draw[red] (10,150/16) -- (3,5);
        \draw[red] (10,70/8) -- (4,5);
      
        \node at (7.5,7.5) {$A_{2,3}^1$};
        \node at (7,7.75) {$A_{1,3}^1$};
        \node at (8,7.25) {$A_{1,3}^1$};

        \node at (3.5,9.6) {$A_{4,3}^1$};
        
        \node at (5-3.5,15-9.6) {$A_{4,3}^1$};

        \draw[very thick] (0,7.5) -- (2.5,10);
        \draw[very thick] (2.5,5) -- (5,7.5);
      
    \end{tikzpicture}
    }

    \caption{Partitions of the region $A_3$. Part (a) shows a partition into sets $A^k$ where $k$ is the escape time. Part (b) shows a subdivision into sets $A_{j,3}^k \subset A^k$ where $j$ is such that $H^k(A_{j,3}^k) \subset A_j$. Red lines in each $A^k$ are the preimages of the $A_1A_2$ boundary under $H^k$.}
    \label{fig:singSetConstruction}
\end{figure}
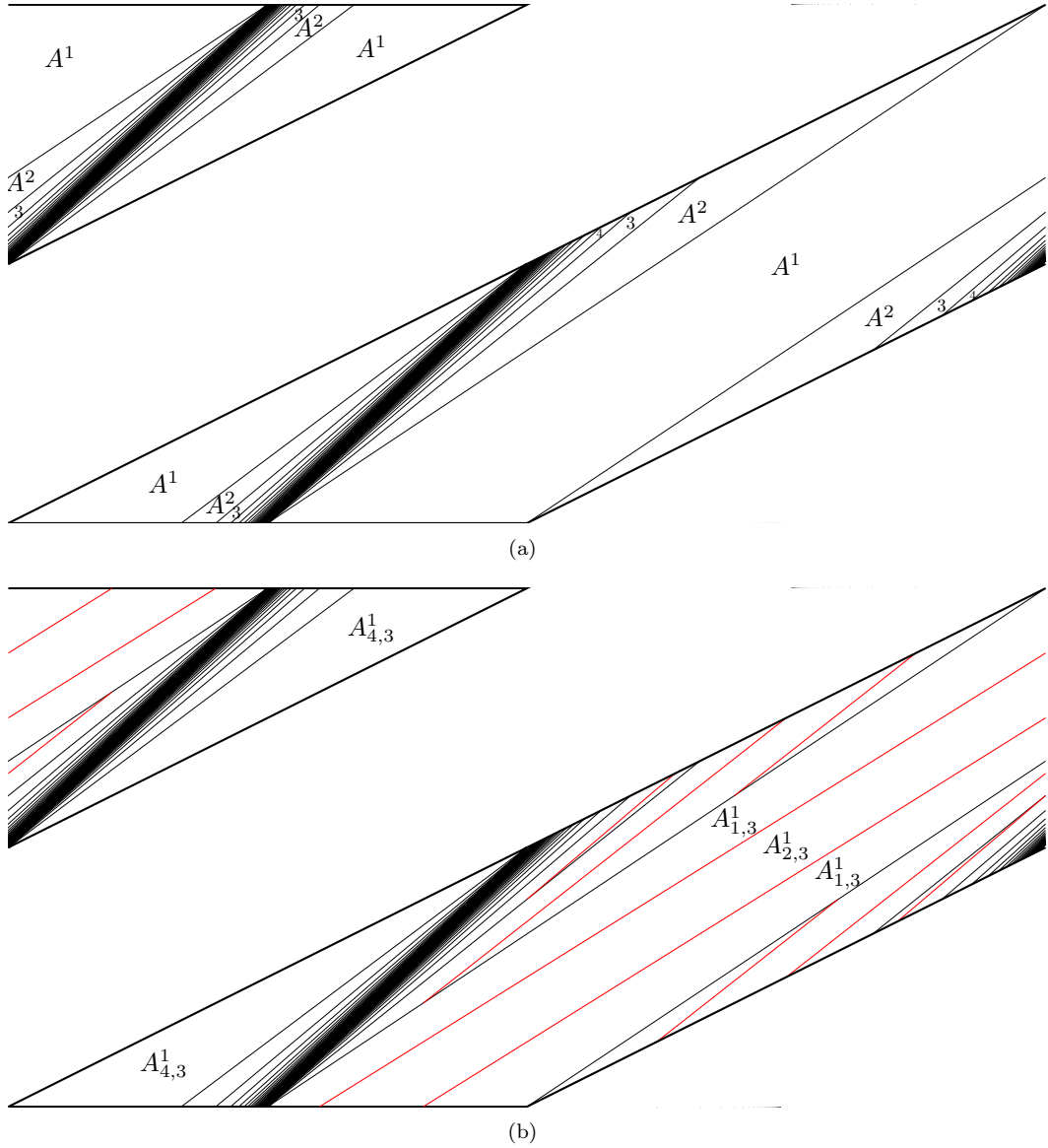

Now consider recurrence to $\sigma$ with return times greater than 1, the white regions of Figure \ref{fig:returnTime1}. Starting with $z \in A_3$, by the definition of $\sigma$, the return time $R(z;H,\sigma) = k+1$ where $k$ is the escape time $E(z;H,A_3)$. Figure \ref{fig:escapeTimes} shows a partition of $A_3$ into sets $A^k$ of constant escape time, bounded by the boundary preimages $H^{-k}(\partial A_3)$. Points in $A^k$ spend $k$ iterates in $A_3$ then escape via $A_1$, $A_2$, or $A_4$ and consequently return to $\sigma$. We partition each $A^k$ based on this escape path, shown as the red lines in Figure \ref{fig:escapeBehaviour}. The labelling $A_{j,i}^k$ is such that $A_{j,i}^k \subset A^k \subset A_i$ and $H^k(A_{j,i}^k) \subset A_j$. It transpires that when points escape after spending 4 or more iterates in $A_3$, they can only do so via $A_1$ or $A_4$. Similarly partitioning $A_2$ and combining with Figure \ref{fig:returnTime1} gives a partition of $\sigma$ into sets on which $DH_\sigma$ is constant. The boundaries of these partition elements are shown in Figure \ref{fig:HSigma} and constitutes, together with $\partial \sigma$, the singularity set $\mathcal{S}$ for $H_\sigma$. We remark that outside of the sets $\varsigma_2$, $\varsigma_3$ the Jacobian $DH_\sigma$ takes values in $\mathcal{M}$. Noting $H_\sigma(\varsigma_2) = H(\varsigma_2) \subset \sigma_3$ and $H_\sigma(\varsigma_3) = H(\varsigma_3) \subset \sigma_2$ we have that within $\varsigma_2$ the Jacobian of $H_\sigma^2$ is given by $M M_3$ for some $M \in \mathcal{M} \cup \{  M_2  \}$ and within $\varsigma_3$ it is given by $M M_2$ for some $M \in \mathcal{M} \cup \{  M_3  \}$. Hence, at almost every $z \in \sigma$ the Jacobian of $H_\sigma$ or $H_\sigma^2$ is some matrix from $\mathcal{M}$. We are now ready to establish non-uniform hyperbolicity.

\begin{proof}[Proof of Proposition \ref{prop:tentMap}]

The proof of Lemma \ref{lemma:itineraries} shows that almost every orbit $H^n(z)$ hits $\sigma$. Similar to LTMs, we can show that almost all of those then continue to return to $\sigma$ with some positive frequency $\alpha_z$. This follows straightforwardly from the fact that $H$ preserves the Lebesgue measure on $\tor$, a compact metric space, and $\sigma$ is measurable. A proof is given in Lemma 6.3.3 of \cite{sturman_mathematical_2006}, originally from \cite{burton_ergodicity_1980}. For large $n$ and a.e. $z$ the cardinality of $\{ 0 \leq i  \leq n - 1 \, | \, H^i(z) \in \sigma \}$ is roughly $\alpha_z n$, certainly
bounded below by $\alpha zn/2$\footnote{By a combinatorial argument, see \cite{sturman_mathematical_2006}.}. The cocycle $DH_z^n$ then contains as many applications of $DH_\sigma$. By the above, applying $DH_\sigma$ either completes a block from $\mathcal{M}$ or does so over the next iterate (the case where we land in $\varsigma_2,\varsigma_3$). At worst, then, we have roughly half as many blocks from $\mathcal{M}$ in $DH_z^n$ as we have returns to $\sigma$. Certainly this proportion is greater than a quarter, so $DH_z^n$ contains at least $\alpha_z n/8$ blocks from $\mathcal{M}$. Defining
\[ K = \inf_{\substack{M \in \mathcal{M}\\ v \in \mathcal{C}}} \frac{\| M v \|}{\| v \|}, \]
Lemma \ref{lemma:cone} gives $K>1$. Noting cone invariance, for any $v_0 \in \mathcal{C}$,
\begin{equation*}
\begin{split}
\frac{1}{n} \log \| DH_z^n v_0  \| & \geq \frac{1}{n} \log \left( K^{\frac{1}{8} \alpha_z n} \| v_0 \|  \right) \\
& = \frac{\alpha_z}{8} \log ( K) + \frac{1}{n} \log \| v_0 \|
\end{split}
\end{equation*}
so that $\chi(z,v_0) \geq \alpha_z \log (K) /8> 0$. We may then extend to non-zero Lyapunov exponents for general $v \neq 0$ using a particular form of Oseledets' theorem in two dimensions (Theorem 3.14 of \citealp{viana_lectures_2014}, see \citealp{myers_hill_exponential_2022}).
\end{proof}

\section{The mixing property}
\label{sec:Hmixing}
In this section we build on hyperbolicity, establishing mixing properties using Theorem \ref{thm:katok-strelcyn}.

\begin{thm}
\label{thm:mixingOTM}
The map $H: \tor \to \tor$ is Bernoulli with respect to the Lebesgue measure.
\end{thm}

\subsection{Nature of local manifolds}

Noting that \textbf{(KS1-2)} were shown in \cite{myers_hill_exponential_2022}, by Theorem \ref{thm:katok-strelcyn}, local unstable and stable manifolds $\gamma_u(z), \gamma_s(z)$ exist at a.e. $z$. By definition, for any $\zeta,\zeta' \in \gamma_u(z)$
\begin{equation}
    \label{eq:unstableManifold}
    \mathrm{dist}(H^{-n}(\zeta),H^{-n}(\zeta')) \rightarrow 0
\end{equation} 
as $n \rightarrow \infty$. Similarly for any $\zeta,\zeta' \in \gamma_s(z)$
\begin{equation}
    \label{eq:stableManifold}
    \mathrm{dist}(H^{n}(\zeta),H^{n}(\zeta')) \rightarrow 0
\end{equation}
as $n\rightarrow \infty$.
Piecewise linearity of $H$ ensures that these local manifolds are line segments containing $z$, aligned with some vector $v = (v_1,v_2)^T$ of \emph{gradient} $v_2/v_1$. The following two lemmas establish bounds on their gradients when mapped under $H$ and its inverse.

\begin{lemma}
\label{lemma:OTMalignment}
For almost every $z$, there exists $m,n \in \mathbb{N}$ such that $H^{m}(\gamma_u(z))$ contains a line segment in $\sigma$ aligned with some $v\in \mathcal{C}$, and $H^{-n}(\gamma_s(z))$ contains a line segment in $\sigma'$ aligned with some $v' \in \mathcal{C}'$.
\end{lemma}

\begin{proof}
By definition of the $\sigma_j, \sigma_j'$, we have that $\sigma_j = H(\sigma_j')$ for $j=1,4$, and $\sigma_j = H^2(\sigma_j')$ for $j=2,3$. For almost every $z$ the number $m = \min \{ k \geq 1 \,| \, H^k(z) \in \sigma \}$ is well defined, as is the cocycle $DH_z^m$. On some portion of $\gamma_u(z)$ around $z$, the cocycle $DH_z^m$ will be constant so that the portion maps to some line segment $\Gamma$ under $H^m$. Hence $H^{m}(\gamma_u(z))$ contains a segment $\Gamma$ in $\sigma$, aligned with some vector $v$. Now if $\Gamma$ lies in $\sigma_1$, its preimage is a segment in $\sigma_1'$ aligned with the vector $M_1^{-1}v$. Now to satisfy (\ref{eq:unstableManifold}), $M_1^{-1}v$ must lie in some stable cone $\mathcal{C}_s'$ which contains all the stable eigenvectors of matrices in $\mathcal{M}'$ and none of the unstable eigenvectors. Hence $v \in M_1 \mathcal{C}_s'$. Similarly if $z\in \sigma_4$ then $v \in M_4 \mathcal{C}_s'$, if $z\in \sigma_2$ then $v \in M_2 M_3 \mathcal{C}_s'$, and if $z\in \sigma_3$ then $v \in M_3M_2 \mathcal{C}_s'$. Such a stable cone $\mathcal{C}_s'$ is given by $\{(v_1,v_2)\neq 0 \,|\, |v_2| \geq |v_1|\}$; one can verify that $M\mathcal{C}_s' \subset \mathcal{C}$ for each $M \in \{ M_1, M_4, M_2M_3, M_3M_2 \}$, verifying $v\in \mathcal{C}$. The argument for $v'\in \mathcal{C}'$ is entirely analogous.
\end{proof}

The expanding and invariance properties of the cone $\mathcal{C}$ formed from $\mathcal{M}$ will be key to growing the images of unstable manifolds. We can ensure stronger expansion by refining the cone, defining $\mathcal{C}_+,\mathcal{C}_-\subset \mathcal{C}$ by
\begin{enumerate}[label={($\mathcal{C}_{+}$)}]
    \item $ 3\,|v_1| \geq |v_2| \geq \varphi \,|v_1|$, $v_1v_2>0$,
\end{enumerate}
\begin{enumerate}[label={($\mathcal{C}_{-}$)}]
    \item $ 3\,|v_1| \geq |v_2| \geq \varphi \,|v_1|$, $v_1v_2<0$.
\end{enumerate}

\begin{lemma}
\label{lemma:mappedOTMalignment}
Let $\Gamma$ be a line segment in $\sigma$, aligned with some $v\in \mathcal{C}$. It follows that $H_\sigma(\Gamma)$ or $H_\sigma^2(\Gamma)$ contains a line segment:
\begin{enumerate}[label={(A\arabic*)}]
    \item Contained within $\sigma_1 \cup \sigma_3$, aligned with some vector in $\mathcal{C}_+$, or
    \item Contained within $\sigma_2 \cup \sigma_4$, aligned with some vector in $\mathcal{C}_-$.
\end{enumerate}
\end{lemma}

\begin{proof}
Suppose first that $\Gamma$ does not lie entirely within $\varsigma_2$ or $\varsigma_3$. Then $\Gamma$ contains a component $\tilde{\Gamma}$ (possibly the whole of $\Gamma$) on which $DH_\sigma$ is a matrix from $\mathcal{M}$. If $H_\sigma\left(\tilde{\Gamma} \right)$ lands in $\sigma_1 \cup \sigma_3$, then this Jacobian is in the subset $\{M_1, M_1M_2^n, M_3M_2^n, M_1M_3^n\} \subset \mathcal{M}$. Case (A1) then follows from verifying that $M\mathcal{C} \subset \mathcal{C}_+$ for each $M$ in this subset. Case (A2) can be argued similarly. If $\Gamma \subset \varsigma_2 \cup \varsigma_3$ then it contains a component on which the Jacobian of $H_\sigma^2$ is in $\mathcal{M}$ and we can follow a similar argument.
\end{proof}

\subsection{Growth lemma}

The recall some useful properties of line segments from \cite{myers_hill_exponential_2022}.

\begin{definition}
Let $\Gamma$ be a line segment. We define the \emph{height} of $\Gamma$ as $\ell_v(\Gamma) = \nu \left( \{y \,|\, (x,y) \in \Gamma \} \right)$, the \emph{width} of $\Gamma$ as $\ell_h(\Gamma) = \nu \left( \{x \,|\, (x,y) \in \Gamma \} \right)$, where $\nu$ is the Lebesgue measure on $\mathbb{R}$.

Given a partition element $A$, we say that $\Gamma$ has \emph{simple intersection} with $A$ if its restriction to $A$ is empty or a single line segment. Conversely we say that $\Gamma$ has \emph{non-simple intersection} with $A$ if its restriction to $A$ contains more than one connected component.
\end{definition}

\begin{lemma}
\label{lemma:unstableGrowth}
Let $\Gamma \subset \sigma$ be a line segment which satisfies either (A1) or (A2) and has simple intersection with each of the $A_j$. Then at least one of the following consequences hold:
\begin{enumerate}[label={(C\arabic*)}]
    \item There exists $k$ such that $H^k(\Gamma)$ contains a line segment having non-simple intersection with some $A_j$,
    \item There exists $k$ such that $H^k(\Gamma)$ contains a line segment $\Lambda$ satisfying (A1) or (A2) with $\ell_v(\Lambda) \geq (1+\delta) \, \ell_v(\Gamma)$ for some $\delta>0$, independent of $\Gamma$.
\end{enumerate}
\end{lemma}

The proof involves splitting into several cases based on the specific location of $\Gamma$ in $\sigma$. The analysis of the first case (roughly up to equation (\ref{eq:sigma1aSimple})) gives a complete exposition of our method, reducing the lemma to checking bounds on growth factors and lengths of partition elements. The other cases are then argued similarly, either by exploiting symmetries or by recalculating bounds on different partition elements. This geometric information, i.e. the equations of the lines which make up $\mathcal{S}$, is vital to our mixing rate analysis in sections \ref{sec:Hsigma}, \ref{sec:polyMixingRate} so we present the full analysis here.

\begin{proof}
Figure \ref{fig:singSetFT} shows the singularity set for the return map $H_\sigma$ over $\sigma_1 \cup \sigma_3 \setminus \varsigma_3$, and the singularity set of $H_\sigma^2$ over $\varsigma_3$. The singularity lines partition $\sigma_1 \cup \sigma_3$ into sets $A_{j,i}^k$ with the same labelling scheme as Figure \ref{fig:escapeBehaviour}.

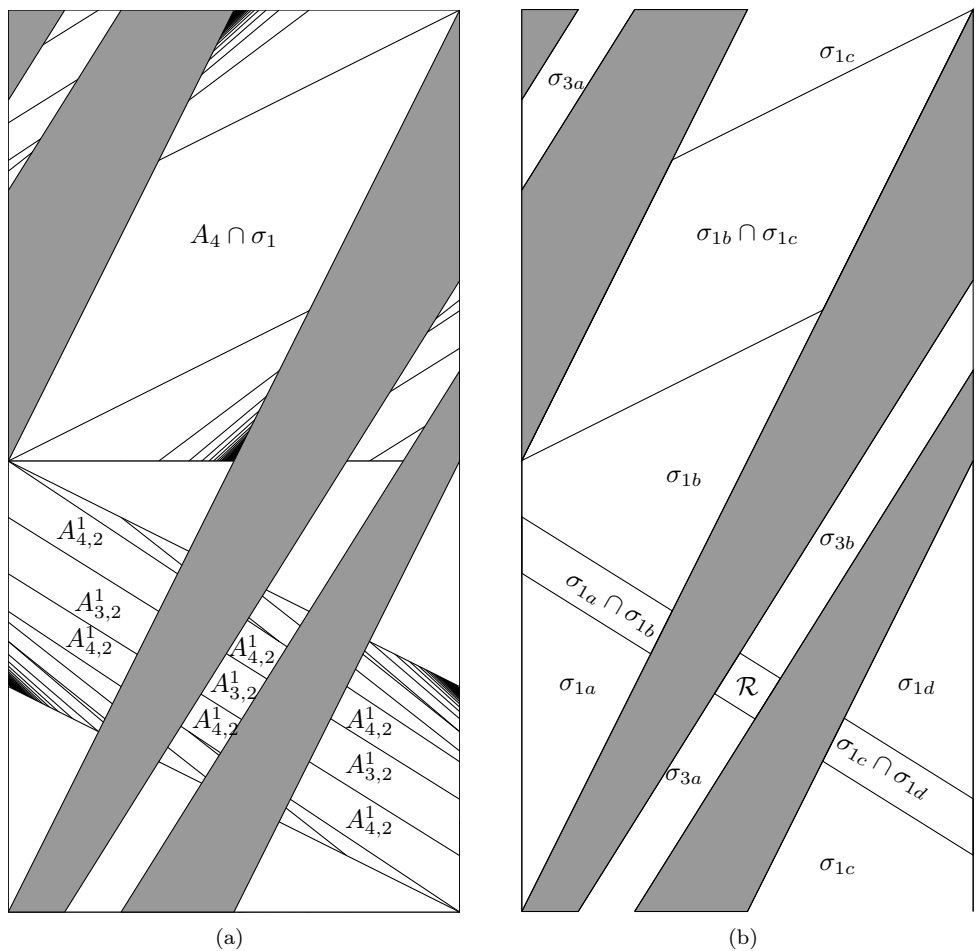
\begin{figure}
    \centering
    
     \subfigure[][]{
    \label{fig:singSetFT}
    \begin{tikzpicture}[scale=1.2]
    
    \clip (0,0) rectangle (5,10);
    
    \draw (0,5) -- (10,5);
    \draw (5,10) -- (5,0);
    
    \draw (0,0) rectangle (10,10);
        
        \foreach \k in {1,...,45}{
        \draw ({ (2*\k+2)*(7.5-7.5)/(2*\k+1)  }, { 7.5 })--({ (2*\k+2)*(10-7.5)/(2*\k+1)  }, { 10 });
        \draw ({ 5+ (2*\k+2)*(5-7.5)/(2*\k+1)  }, { 5 })--({ 5+ (2*\k+2)*(7.5-7.5)/(2*\k+1) }, { 7.5 });
        \draw ({ 2.5 + (2*\k+1)*(5-5)/(2*\k)  }, { 5 })--({ 2.5 + (2*\k+1)*(10-5)/(2*\k)  }, { 10 });
        \draw ({ 12.5 + (2*\k+1)*(5-10)/(2*\k)  }, { 5 }) -- ({ 12.5 + (2*\k+1)*((30*\k+20)/(4*\k+2)-10)/(2*\k)  }, { (30*\k+20)/(4*\k+2) });
        }
        
        \draw (30/4,70/8) -- (4,6);
        \draw (70/8,150/16) -- (7,8);
        \draw (15-30/4,15-70/8) -- (15-5,15-190/28);
        \draw (15-70/8,15-150/16) -- (15-7,15-8);
        \draw (2.5,10) -- (0 , 50/6);
        \draw (1,9) -- (0,15-190/28);

        \draw (5,7)--(90/14,230/28);
        \draw (15-5,15-7)--(15-90/14,15-230/28);
        
        \filldraw[fill=white] (0,5) -- (10,10) -- (5,10) -- (0,7.5) -- (0,5);
        \filldraw[fill=white] (5,5) -- (10,7.5) -- (10,5) -- (5,5);

        \draw (0,150/16) -- (1,10);
        \draw (0,70/8) -- (2,10);
        \draw (10,150/16) -- (3,5);
        \draw (10,70/8) -- (4,5);
      
       \foreach \k in {1,...,45}{
        \draw ({10- (2*\k+2)*(7.5-7.5)/(2*\k+1)  }, { 7.5-5 })--({10- (2*\k+2)*(10-7.5)/(2*\k+1)  }, { 10-5 });
        \draw ({ 5 - (2*\k+2)*(5-7.5)/(2*\k+1)  }, { 5-5 })--({ 5 - (2*\k+2)*(7.5-7.5)/(2*\k+1)  }, { 7.5-5 });
        \draw ({ 7.5 - (2*\k+1)*(5-5)/(2*\k)  }, { 5-5 })--({ 7.5- (2*\k+1)*(10-5)/(2*\k)  }, { 10-5 });
        \draw ({ -2.5 - (2*\k+1)*(5-10)/(2*\k)  }, { 5-5 }) -- ({ -2.5 - (2*\k+1)*((30*\k+20)/(4*\k+2)-10)/(2*\k)  }, { (30*\k+20)/(4*\k+2) -5 });
        }
        
        \draw (10-30/4,70/8-5) -- (10-4,6-5);
        \draw (10-70/8,150/16-5) -- (10-7,8-5);
        \draw (-5+30/4,15-70/8-5) -- (0,15-190/28-5);
        \draw (-5+70/8,15-150/16-5) -- (-5+7,15-8-5);
        \draw (10-2.5,10-5) -- (10, 50/6-5);
        \draw (10-1,9-5) -- (10,15-190/28-5);

        \draw (10-5,7-5)--(10-90/14,230/28-5);
        \draw (0,15-7-5)--(-5+90/14,15-230/28-5);
        
        \filldraw[fill=white] (10,0) -- (0,5) -- (5,5) -- (10,2.5) -- (10,0);
        \filldraw[fill=white] (0,0) -- (5,0) -- (0,2.5) -- (0,0);

        \draw (10,150/16-5) -- (10-1,10-5);
        \draw (10,70/8-5) -- (10-2,10-5);
        \draw (0,150/16-5) -- (10-3,5-5);
        \draw (0,70/8-5) -- (10-4,5-5);

      \filldraw[fill=gray!80] (0,8) -- (10/8,10) -- (2.5,10) -- (0,5) -- (0,8);
      \filldraw[fill=gray!80] (0,9) -- (0,10) -- (10/16,10) -- (0,9);
      \filldraw[fill=gray!80] (0,0) -- (5,10) -- (5,7) -- (10/16,0) -- (0,0);
      \filldraw[fill=gray!80] (10/8,0) -- (2.5,0) -- (5,5) -- (5,6) -- (10/8,0) ;
      
      \filldraw[fill=gray!80] (5,2) -- (50/8,0) -- (7.5,0) -- (5,5) -- (5,2);
      \filldraw[fill=gray!80] (5,1) -- (5,0) -- (90/16,0) -- (5,1);
      \filldraw[fill=gray!80] (5,10) -- (10,0) -- (10,3) -- (90/16,10) -- (5,10);
      \filldraw[fill=gray!80] (50/8,10) -- (7.5,10) -- (10,5) -- (10,4) -- (50/8,10);
  
        \node at (2.5,7.5) {$A_4 \cap \sigma_1$};
  
      \node at (2.5,2.5) {$A_{3,2}^1$};
      \node at (1,3.4) {$A_{3,2}^1$};
      \node at (4,1.6) {$A_{3,2}^1$};
     
      \node at (2.7,2.9) {$A_{4,2}^1$};
      \node at (0.8,4.2) {$A_{4,2}^1$};
      \node at (4,2.1) {$A_{4,2}^1$};
      
      \node at (2.3,2.1) {$A_{4,2}^1$};
      \node at (0.9,3) {$A_{4,2}^1$};
      \node at (4,1) {$A_{4,2}^1$};
      
    \end{tikzpicture}
}
\hspace{1em}
\subfigure[][]{
\label{fig:sigmaSubpartition}
     \begin{tikzpicture}[scale=1.2]
    
    \draw (0,0) -- (0,10);
     \draw (5,0) -- (5,10);

      \filldraw[fill=gray!80] (0,8) -- (10/8,10) -- (2.5,10) -- (0,5) -- (0,8);
      \filldraw[fill=gray!80] (0,9) -- (0,10) -- (10/16,10) -- (0,9);
      \filldraw[fill=gray!80] (0,0) -- (5,10) -- (5,7) -- (10/16,0) -- (0,0);
      \filldraw[fill=gray!80] (10/8,0) -- (2.5,0) -- (5,5) -- (5,6) -- (10/8,0);

    \draw (0,0) -- (0,35/8) -- (35/21,70/21) -- (0,0);
    
   \draw (5,10) -- (10/6,50/6) -- (0,5) -- (0,30/8) -- (30/21,60/21) --  (5,10);
    
   \draw (2.5,10) -- (0,5) -- (5-10/6,10-20/6) -- (5,10);
    \draw (5,0) -- (5,5-30/8) -- (5-30/21,5-60/21) --  (2.5,0);
    
    \draw (5,5) -- (5,5-35/8) -- (5-35/21,5-70/21) -- (5,5);
 
   \draw (10/8,0) -- (5-190/89,5-215/89) -- (5-230/89,5-190/89) -- (10/16,0);
    \draw (10/16,10) -- (0,9) -- (0,8) -- (10/8,10);
 
    \draw (5,7) -- (190/89,215/89) -- (230/89,190/89) -- (5,6);
 
    \node at (10/16,2.5) {$\sigma_{1a}$};
    
    \node at (1.8,4.8) {$\sigma_{1b}$};
    
    \node at (3.5,9.5) {$\sigma_{1c}$};
    \node at (3.5,0.5) {$\sigma_{1c}$};
    
    \node at (5-10/16,2.5) {$\sigma_{1d}$};
    
    \node at (0.5,9.2) {$\sigma_{3a}$};
    \node at (1.8,1.5) {$\sigma_{3a}$};
     \node at (3.5,4.1) {$\sigma_{3b}$};
    
     \node at (2.5,2.5) {$\mathcal{R}$};
     \node [rotate=-30] at (1,3.4) {$\sigma_{1a} \cap \sigma_{1b}$};
     \node [rotate=-30] at (4,1.6) {$\sigma_{1c} \cap \sigma_{1d}$};
    \node at (2.5,7.5) {$\sigma_{1b} \cap \sigma_{1c}$};

    \end{tikzpicture}
    }

   \caption{Part (a) shows the singularity curves dividing up $\sigma_1 \cup \sigma_3$ with some key partition elements labelled. The elements $A_{3,2}^1$, $A_4 \cap \sigma$ split $\sigma_1 \cup \sigma_3$ into six subsets $\sigma_{1a},\dots,\sigma_{3b}$, any two of which are either disjoint or have intersection given by $A_4 \cap \sigma$ or one of the three subsets which make up $A_{3,2}^1$, see part (b). $\mathcal{R}$ denotes the set $\sigma_{3a} \cap \sigma_{3b}$.}
    \label{fig:caseA1Partition}
\end{figure}

Let $\Gamma$ satisfy case (A1) and suppose it has non-simple intersection with $A_{4,2}^1$. Now since $\Gamma$ has simple intersection with $A_3$, observing Figure \ref{fig:singSetFT} it is clear that $\Gamma$ must traverse $A_{3,2}^1$. Restricting $\Gamma^2 = \Gamma \cap A_2$, $H\left(\Gamma^2\right) \subset H(\Gamma)$ is a line segment which has non-simple intersection with $A_4$, i.e. (C1) is satisfied with $k=1$. Assume, then, that $\Gamma$ has simple intersection with $A_{4,2}^1$ and therefore does not traverse $A_{3,2}^1$. If $\Gamma \subset \sigma_3$ then $\Gamma$ lies entirely within one of two sets $\sigma_{3a}$, $\sigma_{3b}$ (shown in Figure \ref{fig:sigmaSubpartition}) whose union is $\sigma_3$, intersection is $\mathcal{R} = A_{3,2}^1 \cap \sigma_3$. For $\Gamma \subset \sigma_1$, simple intersection with $A_3$ implies that $\Gamma$ does not traverse $A_4 \cap \sigma_1$. This, together with the two disjoint sets which make up $A_{3,2}^1 \cap \sigma_1$, implies that $\Gamma$ lies entirely within one of four subsets $\sigma_{1a},\dots,\sigma_{1d}$, shown in Figure \ref{fig:sigmaSubpartition}. The behaviour of $H_\sigma$ over the sets $\sigma_{1a}$, $\sigma_{1b}$ is shown explicitly in Figures \ref{fig:sigma1a}, \ref{fig:sigma1b}.
 
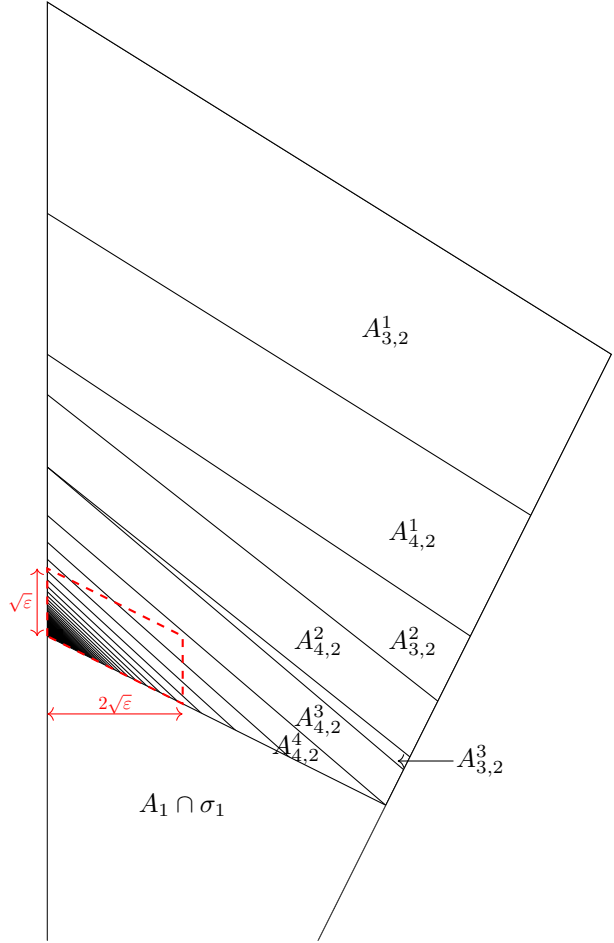
\begin{figure}
    \centering

    \begin{tikzpicture}
    \node at (0,0) {
    
    \begin{tikzpicture}[scale=4.5]
    \draw (0,1.6) -- (0,35/8) -- (35/21,70/21) -- (0.8,1.6);
    \draw[->] (1.2,2.1308) -- (1.035,2.1308);
    \node at (1.28,2.1308) {$A_{3,2}^3$};

\begin{scope}

    \clip (0,1.6) -- (0,35/8) -- (35/21,70/21) -- (0.8,1.6);
    \draw (0,5) -- (10,5);
    \draw (5,10) -- (5,0);
    
    \fill[black] (0,2.3) rectangle (0.045,2.52);
    
    \draw (0,0) rectangle (10,10);
        
        \foreach \k in {1,...,45}{
        \draw ({ (2*\k+2)*(7.5-7.5)/(2*\k+1)  }, { 7.5 })--({ (2*\k+2)*(10-7.5)/(2*\k+1)  }, { 10 });
        \draw[blue] ({ 5+ (2*\k+2)*(5-7.5)/(2*\k+1)  }, { 5 })--({ 5+ (2*\k+2)*(7.5-7.5)/(2*\k+1) }, { 7.5 });
        \draw ({ 2.5 + (2*\k+1)*(5-5)/(2*\k)  }, { 5 })--({ 2.5 + (2*\k+1)*(10-5)/(2*\k)  }, { 10 });
        \draw ({ 12.5 + (2*\k+1)*(5-10)/(2*\k)  }, { 5 }) -- ({ 12.5 + (2*\k+1)*((30*\k+20)/(4*\k+2)-10)/(2*\k)  }, { (30*\k+20)/(4*\k+2) });
        }
        
        \draw (30/4,70/8) -- (4,6);
        \draw (70/8,150/16) -- (7,8);
        \draw (15-30/4,15-70/8) -- (15-5,15-190/28);
        \draw (15-70/8,15-150/16) -- (15-7,15-8);
        \draw (2.5,10) -- (0 , 50/6);
        \draw (1,9) -- (0,15-190/28);

        \draw (5,7)--(90/14,230/28);
        \draw (15-5,15-7)--(15-90/14,15-230/28);
        
        \filldraw[fill=white] (0,5) -- (10,10) -- (5,10) -- (0,7.5) -- (0,5);
        \filldraw[fill=white] (5,5) -- (10,7.5) -- (10,5) -- (5,5);

        \draw (0,150/16) -- (1,10);
        \draw (0,70/8) -- (2,10);
        \draw (10,150/16) -- (3,5);
        \draw (10,70/8) -- (4,5);
      
       \foreach \k in {1,...,45}{
        \draw ({10- (2*\k+2)*(7.5-7.5)/(2*\k+1)  }, { 7.5-5 })--({10- (2*\k+2)*(10-7.5)/(2*\k+1)  }, { 10-5 });
        \draw ({ 5 - (2*\k+2)*(5-7.5)/(2*\k+1)  }, { 5-5 })--({ 5 - (2*\k+2)*(7.5-7.5)/(2*\k+1)  }, { 7.5-5 });
        \draw ({ 7.5 - (2*\k+1)*(5-5)/(2*\k)  }, { 5-5 })--({ 7.5- (2*\k+1)*(10-5)/(2*\k)  }, { 10-5 });
        \draw ({ -2.5 - (2*\k+1)*(5-10)/(2*\k)  }, { 5-5 }) -- ({ -2.5 - (2*\k+1)*((30*\k+20)/(4*\k+2)-10)/(2*\k)  }, { (30*\k+20)/(4*\k+2) -5 });
        }
        
        \draw (10-30/4,70/8-5) -- (10-4,6-5);
        \draw (10-70/8,150/16-5) -- (10-7,8-5);
        \draw (-5+30/4,15-70/8-5) -- (0,15-190/28-5);
        \draw (-5+70/8,15-150/16-5) -- (-5+7,15-8-5);
        \draw (10-2.5,10-5) -- (10, 50/6-5);
        \draw (10-1,9-5) -- (10,15-190/28-5);

        \draw (10-5,7-5)--(10-90/14,230/28-5);
        \draw (0,15-7-5)--(-5+90/14,15-230/28-5);
        
        \filldraw[fill=white] (10,0) -- (0,5) -- (5,5) -- (10,2.5) -- (10,0);
        \filldraw[fill=white] (0,0) -- (5,0) -- (0,2.5) -- (0,0);

        \draw (0,0) -- (5,10);
        
        \draw (10,150/16-5) -- (10-1,10-5);
        \draw (10,70/8-5) -- (10-2,10-5);
        \draw (0,150/16-5) -- (10-3,5-5);
        \draw (0,70/8-5) -- (10-4,5-5);

     \node at (1,3.4) {$A_{3,2}^1$};
     \node at (1.08,2.8) {$A_{4,2}^1$};
     \node at (1.08,2.48) {$A_{3,2}^2$};
     \node at (0.8,2.48) {$A_{4,2}^2$};
     \node at (0.8,2.25) {$A_{4,2}^3$};
     
     \node at (0.735,2.167) {$A_{4,2}^4$};
     
     \node at (0.4,2) {$A_1 \cap \sigma_1$};

\end{scope}
        \draw [red, dashed, thick] (0,2.5) -- (0,2.7) -- (0.4,2.5) -- (0.4,2.3) -- (0,2.5);
    \draw [<->, red] (0,2.27) -- (0.4,2.27);
    \draw [<->, red] (-0.03,2.5) -- (-0.03,2.7);
    \node[scale=0.7] at (-0.08,2.6) {\color{red}$\sqrt{\varepsilon}$\color{black}};
    \node[scale=0.7] at (0.2,2.3) {\color{red}$2\sqrt{\varepsilon}$\color{black}};

    \end{tikzpicture} };

    \end{tikzpicture}
    
    \caption{The singularity set of $H_\sigma$ over $\sigma_{1a}$. Unlabelled sets are given by $A_{4,2}^k$ for $k\geq 5$ which limit onto the point $(0,1/4)$ in the obvious fashion. The dashed red line is $\partial P(\varepsilon)$, useful for establishing \textbf{(KS1)} for $H_\sigma$. }
    \label{fig:sigma1a}
\end{figure}

Let $\| \cdot \|$ denote the $\| \cdot \|_\infty$ norm. Starting with $\sigma_{1a}$, $DH_\sigma$ takes values in $\mathcal{M}_{1a} = \{ M_1, M_4M_2^k, M_3M_2^l \,| \, k\in \mathbb{N}, \, l=1,2,3\}$. The unlabelled sets in Figure \ref{fig:sigma1a} are the partition elements $A_{4,2}^k$ for $k \geq 5$, limiting onto the point $(0,1/4)$ as $k\rightarrow \infty$ in the obvious fashion. We remark that any $\Gamma \subset \sigma_{1a}$ has simple intersection with all of the partition elements $A_{i,j}^k \subset \sigma_{1a}$. If $\Gamma$ is entirely contained within some partition element $A$ corresponding to $M \in \mathcal{M}_{1a}$, and is aligned with some unit vector $v\in \mathcal{C}_+$, then $\ell_v\left( H_\sigma(\Gamma) \right) = \| M v  \| \ell_v(\Gamma)$. Minimum expansion factors are straightforward to calculate. Parameterise unit vectors in $\mathcal{C}_+$ by $(v_1,1)^T$ where $1/3 \leq v_1 \leq 13/21$ and write the components of matrices $M \in \mathcal{M}$ as $\big(\begin{smallmatrix}
  a & b\\
  c & d
\end{smallmatrix}\big)$ . Then by cone invariance and the fact that vectors $(v_1,v_2)^T \in \mathcal{C}$ have norm $|v_2|$, we have that $\| M v \| = |cv_1 + d|$. This is monotone increasing in $v_1$ if $\mathrm{sgn}(c) = \mathrm{sgn}(d)$, monotone decreasing if $\mathrm{sgn}(c)\neq \mathrm{sgn}(d)$, so that $\| M v \|$ is minimal on $(1/3,1)^T$ or $(13/21,1)^T$ in these respective cases. Table \ref{tab:tab2} shows the components of matrices $M \in \mathcal{M}$ and the minimum expansion factors $K_+(M)$ which follow.


\begin{table}[ht]
    \centering
    \begin{tabular}{c|c|c|c}
        $M$ & Components & $K_+(M)$ & $K_-(M)$ \\ \hline
$M_1$ & $\left(\begin{matrix}1 & 2\\2 & 5\end{matrix}\right) $ & \( \displaystyle  \frac{17}{3} \) & \( \displaystyle  \frac{79}{21} \)  \\
$M_4$ & $ \left(\begin{matrix}1 & -2\\-2 & 5\end{matrix}\right) $ & \( \displaystyle  \frac{79}{21} \) & \( \displaystyle  \frac{17}{3} \)  \\
$M_1M_2^n$ & $(-1)^n \left(\begin{matrix}2 n + 1 & 2 n + 2\\6 n + 2 & 6 n + 5\end{matrix}\right) $ & \( \displaystyle  8 n + \frac{17}{3} \) & \( \displaystyle  \frac{16 n}{7} + \frac{79}{21} \)  \\
$M_1M_3^n$ & $(-1)^n \left(\begin{matrix}1 - 6 n & 6 n + 2\\2 - 14 n & 14 n + 5\end{matrix}\right) $ & \( \displaystyle  \frac{16 n}{3} + \frac{131}{21} \) & \( \displaystyle  \frac{56 n}{3} + \frac{13}{3} \)  \\
$M_2M_3^n$ & $(-1)^n \left(\begin{matrix}1 - 6 n & 6 n + 2\\10 n - 2 & - 10 n - 3\end{matrix}\right) $ & \( \displaystyle  \frac{80 n}{21} + \frac{89}{21} \) & \( \displaystyle  \frac{40 n}{3} + \frac{7}{3} \)  \\
$M_3M_2^n$ & $(-1)^n \left(\begin{matrix}1 - 6 n & - 6 n - 2\\2 - 10 n & - 10 n - 3\end{matrix}\right) $ & \( \displaystyle  \frac{40 n}{3} + \frac{7}{3} \) & \( \displaystyle  \frac{80 n}{21} + \frac{89}{21} \)  \\
$M_4M_2^n$ & $(-1)^n \left(\begin{matrix}1 - 6 n & - 6 n - 2\\14 n - 2 & 14 n + 5\end{matrix}\right) $ & \( \displaystyle  \frac{56 n}{3} + \frac{13}{3} \) & \( \displaystyle  \frac{16 n}{3} + \frac{131}{21} \)  \\
$M_4M_3^n$ & $(-1)^n \left(\begin{matrix}2 n + 1 & - 2 n - 2\\- 6 n - 2 & 6 n + 5\end{matrix}\right) $ & \( \displaystyle  \frac{16 n}{7} + \frac{79}{21} \) & \( \displaystyle  8 n + \frac{17}{3} \)  \\
    \end{tabular}
    \caption{Minimum expansion factors $K_\pm(M) = \inf_{v\in \mathcal{C}_\pm}\| M v \| / \| v \|$ for each $M \in \mathcal{M}$ over the cones $\mathcal{C}_\pm$.}
    \label{tab:tab2}
\end{table}

If $\Gamma$ intersects $A_{4,2}^4$ and $A_{3,2}^3$ (traversing $A_{4,2}^3$) then $H^{3}\left(\Gamma \cap A_{4,2}^3\right)$ is a line segment in $A_2' \cap A_4$, connecting the $A_3,A_4$ boundary to the $A_2,A_4$ boundary. Noting that $A_2' \cap A_4$ is made up of two quadrilaterals, see Figure \ref{fig:firstPartitions}, there are two possible ways this can occur. Firstly, it can connect points $(x,1)$ to $(2y -1,y)$ with $ 1/2 \leq x\leq 3/4$. Its image under $F$ then connects $(x,1)$ to $(1,y)$ so that then, shearing vertically by $G$, its image under $H$ connects $(x,2-2x)$ to $(1,y)$, passing through $y=0$. Since $x \leq 3/4$, we have $2-2x \geq 1/2$ so that $H^4\left(\Gamma \cap A_{4,2}^3\right)$ must have non-simple intersection with $A_2$. The second case, where $H^3\left(\Gamma \cap A_{4,2}^3\right)$ connects points $(x,1/2)$ and $(2y-1/2,y)$, is similar so that (C1) is satisfied.

Assume, then, that $\Gamma$ does not traverse $A_{4,2}^3$. Two possible cases follow; either $\Gamma$ lies entirely below the upper boundary of $A_{4,2}^3$, or $\Gamma$ lies entirely above the lower boundary of $A_{4,2}^3$. In the first case let $\Gamma_1 = \Gamma \cap A_1$. If $K_+(M_1) \, \ell_v(\Gamma_1) > \ell_v(\Gamma)$, then we may take $\Lambda = H(\Gamma_1) \subset H_\sigma(\Gamma)$ to satisfy (C2). Taking $K_+(M_1) = 17/3$ from Table \ref{tab:tab2}, this holds provided that $\ell_v(\Gamma_1)/\ell_v(\Gamma) > 3/17$. Noting that $\Gamma \subset A_1 \cup A_2$, if the above inequality does not hold, then the proportion of $\Gamma$ in $A_2$ satisfies $\ell_v(\Gamma_2)/\ell_v(\Gamma) > 14/17$. Observing Figure \ref{fig:sigma1a}, $\Gamma_2$ intersects some collection of sets $A_{4,2}^k$, indexed by a consecutive subset $\{ k_0, k_0+1, ...\}\subset \mathbb{N}$ with $k_0 \geq 3$. Assume that $\Gamma_2$ intersects just two of these sets $\Gamma_2 = \Gamma_{k_0} \cup \Gamma_{k_0+1}$. As seen in \cite{myers_hill_exponential_2022}, if 
\[\frac{1}{K_+\left(M_4M_2^{k_0}\right)} + \frac{1}{K_+\left(M_4M_2^{k_0+1}\right)} < 1 \]
then at least one of $\Gamma_k = \Gamma_{k_0}$, $\Gamma_{k_0+1}$ satisfies $\ell_v\left(H^{k+1}(\Gamma_k)\right) > \ell_v(\Gamma_2)$ and by extension if 
\[\frac{1}{K_+\left(M_4M_2^{k_0}\right)} + \frac{1}{K_+\left(M_4M_2^{k_0+1}\right)} < \frac{1}{\alpha} \]
then $\ell_v\left(H^{k+1}(\Gamma_k)\right) > \alpha \ell_v(\Gamma_2)$. Now noting that $K_+\left(M_4M_2^k\right)$ is monotonic increasing in $k$ we have
\[\sum_{k=k_0}^{k_0+1} \frac{1}{K_+\left(M_4M_2^k\right)} \leq \sum_{k=3}^4 \frac{1}{K_+\left(M_4M_2^k\right)} = \frac{3}{181} + \frac{3}{237} < \frac{14}{17} \]
so that, together with $\ell_v(\Gamma_2)/\ell_v(\Gamma) > 14/17$, for some $k$ condition (C2) follows by taking $\Lambda = H^{k+1}\left(\Gamma \cap A_{4,2}^k\right)$. The case where $\Gamma$ intersects just one of the $A_{4,2}^k$ follows as a trivial consequence.

Suppose $\Gamma \subset \sigma_{1a}$ violates the lemma, by the above we have that $\Gamma$ intersects three or more of the $A_{4,2}^k$, which by the geometry of the partition (see Figure \ref{fig:sigma1a}) implies
\begin{enumerate}[label={($\dagger$)}]
    \item $\Gamma$ traverses $A_{4,2}^k$ for some $k \geq 4$, connecting the lines $\mathcal{L}_k:$ $y=\frac{k+1-4kx}{4k+2}$ and $\mathcal{L}_{k-1}:$ $y=\frac{k-4(k-1)x}{4k-2}$.
\end{enumerate}
We will show that this leads to a contradiction through an inductive argument. If $\Gamma$ intersects $A_{4,2}^3$, it must traverse $A_{4,2}^4$. Let $y_k = (k+1)/(4k+2)$ be the sequence of points where $\mathcal{L}_k$ meets $x=0$. Since the gradients of $\mathcal{L}_k$ are monotone decreasing in $k$, a lower bound $h_4 \leq \ell_v\left(\Gamma \cap A_{4,2}^4\right)$ is given by $y_4'-y_4$ where $(x_4',y_4')$ is the intersection of the lines $y=y_4+\varphi x$ and $\mathcal{L}_3: \, y= (4-12x)/14$. Specifically
\begin{equation}
    \label{eq:h4}
    h_4 = \frac{191}{675} - \frac{5}{18} = \frac{7}{1350}.
\end{equation}
As before let $\Gamma_2 = \Gamma \cap A_2$. Observing Figure \ref{fig:sigma1a}, since $\mathcal{L}_3$ meets the boundary of $A_1$ and $A_2$ at the point $(1/10,1/5)$, the height of $\Gamma_2$ is bounded by $\ell_v(\Gamma_2) \leq L_3 = y_2 - 1/5 = 1/10$. Letting $\Lambda = H^5\left(\Gamma \cap A_{4,2}^4\right)$, we have that
\[ \ell_v(\Lambda) \geq K_+\left(M_4M_2^4\right)\, h_4 = \frac{56(4) + 13}{3} \frac{7}{1350} = \frac{553}{1350} \approx 0.4096 \]
and
$\ell_v(\Gamma) < (17/14) \, \ell_v(\Gamma_2) \leq 17/140 \approx 0.1214$, so that (C2) is satisfied. For the inductive step, assume that $\Gamma$ traverses $A_{4,2}^k$, but does not traverse $A_{4,2}^{k-1}$. Using the same method as before we calculate
\begin{equation}
    \label{eq:hk}
    h_k = \frac{21}{2\left(2k+1\right)\left(68k-47\right)}
\end{equation}
and
\begin{equation}
    \label{eq:Lk-1}
    L_{k-1} = \frac{k-1}{4k-6} - \frac{k-1}{4k-2} \frac{\left(k-1\right)}{\left(2k-3\right)\left(2k-1\right)}.
\end{equation}
Then (C2) is satisfied with $\Lambda = H^{k+1}\left(\Gamma \cap A_{4,2}^k\right)$ provided that $K_+\left(M_4M_2^k\right) \, h_k > (17/14) L_{k-1}$, i.e.
\begin{equation}
    \label{eq:sigma1aInductiveStep}
    \frac{56k+13}{3} \frac{21}{2\left(2k+1\right)\left(68k-47\right)} - \frac{17}{14} \frac{\left(k-1\right)}{\left(2k-3\right)\left(2k-1\right)} > 0,
\end{equation}
which holds for all $k>4$ as required. It follows by induction that if $\Gamma$ violates the lemma it must not traverse any $A_{4,2}^k$ for $k\geq 3$, contradicting ($\dagger$), so that the lemma must hold when $\Gamma \subset \sigma_{1a}$ lies entirely below the upper boundary of $A_{4,2}^3$. The case where $\Gamma\subset \sigma_{1a}$ lies entirely above the lower boundary of $A_{4,2}^3$ is more straightforward, with (C2) following from the inequality
\begin{equation}
    \label{eq:sigma1aSimple}
    \sum_{k=1}^3 \frac{1}{K_+\left( M_3M_2^k\right)} +  \frac{1}{K_+\left(M_4M_2^k\right)} = \sum_{k=1}^3 \frac{3}{56k+13} + \frac{3}{40k+7} \approx 0.206 < 1.
\end{equation}
The lemma holds, then, for general $\Gamma \subset \sigma_{1a}$.

\begin{figure}
    \centering
    \begin{tikzpicture}
    \node at (0,0) {
        \begin{tikzpicture}[scale=3]
    
    \draw[->] (2.1,3.965) -- (1.88,3.965);
    \node at (2.2,3.965) {$A_{3,2}^2$};
        
    \draw[<->] (70/19+0.05,5) -- (70/19+0.05,140/19);
    \node at (3.9,235/38) {$L_0$};
    
    \draw  (1,7) -- (0,5) -- (0,30/8) -- (30/21,60/21) --  (3.75,7.5);
    \clip (0,30/8) -- (30/21,60/21) -- (3.75,7.5) -- (1,7) -- (0,5) -- (0,30/8);
    \draw (0,5) -- (10,5);
    \draw (5,10) -- (5,0);
    
    \fill[black] (2.48,5) rectangle (2.54,5.04);
    
    \draw (0,0) rectangle (10,10);
        
        \foreach \k in {1,...,45}{
        \draw ({ (2*\k+2)*(7.5-7.5)/(2*\k+1)  }, { 7.5 })--({ (2*\k+2)*(10-7.5)/(2*\k+1)  }, { 10 });
        \draw ({ 5+ (2*\k+2)*(5-7.5)/(2*\k+1)  }, { 5 })--({ 5+ (2*\k+2)*(7.5-7.5)/(2*\k+1) }, { 7.5 });
        \draw ({ 2.5 + (2*\k+1)*(5-5)/(2*\k)  }, { 5 })--({ 2.5 + (2*\k+1)*(10-5)/(2*\k)  }, { 10 });
        \draw ({ 12.5 + (2*\k+1)*(5-10)/(2*\k)  }, { 5 }) -- ({ 12.5 + (2*\k+1)*((30*\k+20)/(4*\k+2)-10)/(2*\k)  }, { (30*\k+20)/(4*\k+2) });
        }
        
        \draw (30/4,70/8) -- (4,6);
        \draw (70/8,150/16) -- (7,8);
        \draw (15-30/4,15-70/8) -- (15-5,15-190/28);
        \draw (15-70/8,15-150/16) -- (15-7,15-8);
        \draw (2.5,10) -- (0 , 50/6);
        \draw (1,9) -- (0,15-190/28);

        \draw (5,7)--(90/14,230/28);
        \draw (15-5,15-7)--(15-90/14,15-230/28);
        
        \filldraw[fill=white] (0,5) -- (10,10) -- (5,10) -- (0,7.5) -- (0,5);

        \draw (0,150/16) -- (1,10);
        \draw (0,70/8) -- (2,10);
        \draw (10,150/16) -- (3,5);
        \draw (10,70/8) -- (4,5);
      
       \foreach \k in {1,...,45}{
        \draw ({10- (2*\k+2)*(7.5-7.5)/(2*\k+1)  }, { 7.5-5 })--({10- (2*\k+2)*(10-7.5)/(2*\k+1)  }, { 10-5 });
        \draw ({ 5 - (2*\k+2)*(5-7.5)/(2*\k+1)  }, { 5-5 })--({ 5 - (2*\k+2)*(7.5-7.5)/(2*\k+1)  }, { 7.5-5 });
        \draw ({ 7.5 - (2*\k+1)*(5-5)/(2*\k)  }, { 5-5 })--({ 7.5- (2*\k+1)*(10-5)/(2*\k)  }, { 10-5 });
        \draw[red] ({ -2.5 - (2*\k+1)*(5-10)/(2*\k)  }, { 5-5 }) -- ({ -2.5 - (2*\k+1)*((30*\k+20)/(4*\k+2)-10)/(2*\k)  }, { (30*\k+20)/(4*\k+2) -5 });
        }
        
        \draw (10-30/4,70/8-5) -- (10-4,6-5);
        \draw (10-70/8,150/16-5) -- (10-7,8-5);
        \draw (-5+30/4,15-70/8-5) -- (0,15-190/28-5);
        \draw (-5+70/8,15-150/16-5) -- (-5+7,15-8-5);
        \draw (10-2.5,10-5) -- (10, 50/6-5);
        \draw (10-1,9-5) -- (10,15-190/28-5);

        \draw (10-5,7-5)--(10-90/14,230/28-5);
        \draw (0,15-7-5)--(-5+90/14,15-230/28-5);
        
        \filldraw[fill=white] (10,0) -- (0,5) -- (5,5) -- (10,2.5) -- (10,0);
        \filldraw[fill=white] (0,0) -- (5,0) -- (0,2.5) -- (0,0);

        \draw (0,0) -- (5,10);
        \draw (0,5) -- (2.5,10);
        
        \draw (10,150/16-5) -- (10-1,10-5);
        \draw (10,70/8-5) -- (10-2,10-5);
        \draw (0,150/16-5) -- (10-3,5-5);
        \draw (0,70/8-5) -- (10-4,5-5);

    \node at (1.199,4.296) {$A_{4,2}^2$};
      
    \node at (1.199,3.9) {$A_{4,2}^1$};
    
    \node at (1,3.4) {$A_{3,2}^1$};
    
    \node at (1.6,4.7) {$A_1 \cap \sigma_1$};
    
    \draw[dashed] (0,70/12) -- (10,10);  
     
    \node at (1.6,5.4) {$A_{4,3}^1$};
     \node at (1.6,6.1) {$A_4 \cap \sigma_1$};
      
         \node[scale=1.5] at (2,6.78) {$\mathcal{P}$}; 
      
    \end{tikzpicture}};

    \end{tikzpicture}     
         
      \caption{The singularity set of $H_\sigma$ over the lower part of $\sigma_{1b}$ with the top portion of $A_4 \cap \sigma_1$ omitted. Unlabelled sets are given by $A_{4,3}^k$ for $k\geq 2$ which limit onto the point $(1/4,1/2)$ in the obvious fashion. The segment $\mathcal{P}$ is the preimage under $H$ of the segment joining $(1/2,3/4)$ to $(1,1)$ in $S_4$. The length $L_0$ denotes maximum height of any segment in $\sigma_{1b}$ bounded by $\mathcal{P}$ and $y=1/2$.}
    \label{fig:sigma1b}
\end{figure}
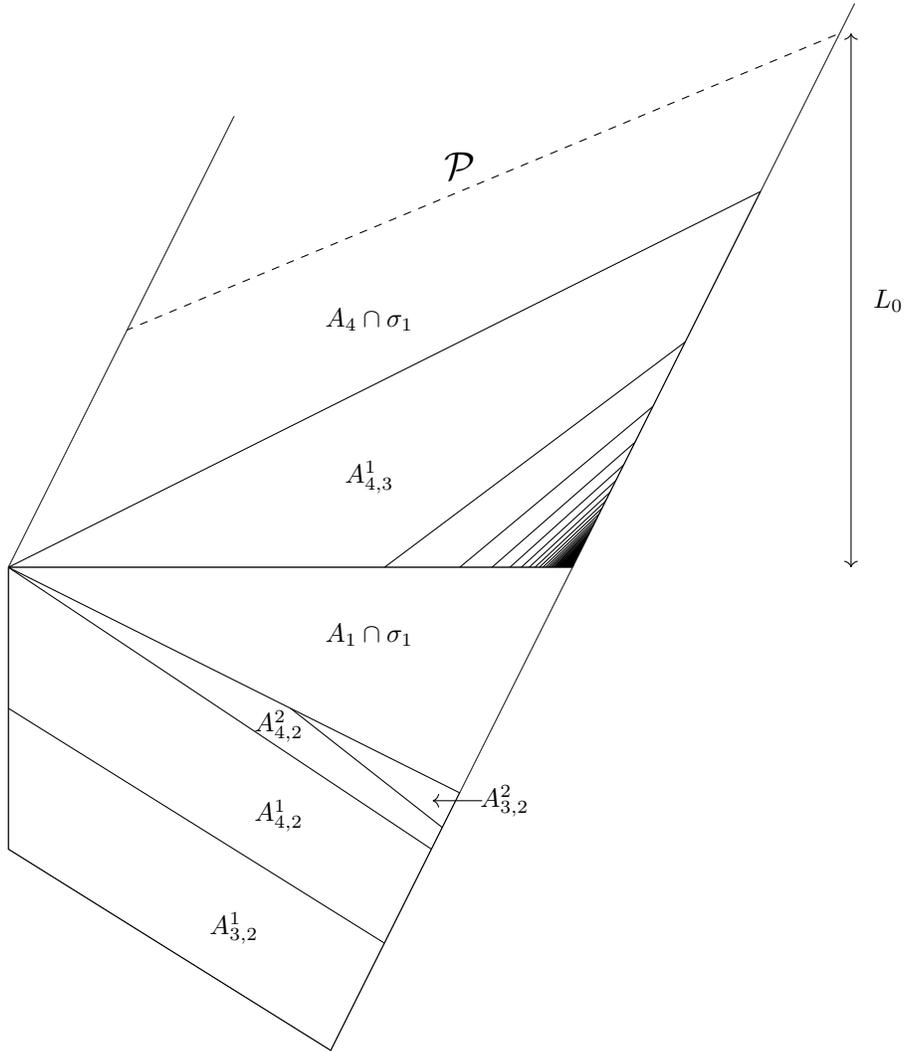

Moving onto the case $\Gamma \subset \sigma_{1b}$, write its intersections with the lower and upper regions $y\leq 1/2$ and $y \geq 1/2$ as $\Gamma_L$ and $\Gamma_U$ respectively. Observing Figure \ref{fig:sigma1b}, $\Gamma_L$ can intersect up to 5 partition elements from $\mathcal{A}_L = \{ A_{3,2}^1, \dots, A_1 \cap \sigma \}$, on which $DH_\sigma$ takes a value in $\mathcal{M}_L = \{ M_1, M_3M_2^k, M_4M_2^k \, | \, k=1,2 \}$.
Let
\[ \alpha = \sum_{M \in \mathcal{M}_L} {\frac{1}{K_+(M)}} = \frac{3}{17}+ \sum_{k=1}^{2}\left(\frac{3}{40k+7}+\ \frac{3}{56k+13}\right) \approx 0.342. \]
Dividing through by $\alpha$, for any subset $\mathcal{N} \subset \mathcal{M}_L$ (including $\varnothing$ and $\mathcal{M}_L$) we have
\[ \sum_{M \in \mathcal{N}} {\frac{1}{\alpha K_+(M)}} \leq 1. \]
Hence we may always expand from some $A \in \mathcal{A}_L$, taking $\Lambda = H_\sigma(\Gamma \cap A)$, which by the above inequality satisfies $\alpha \ell_v(\Lambda) \geq \ell_v(\Gamma_L)$. Hence (C2) is satisfied when $\ell_v(\Gamma_L) > \alpha \ell_v(\Gamma)$. It remains to show the case $\ell_v(\Gamma_L) \leq \alpha \ell_v(\Gamma)$, i.e.
\begin{equation}
    \label{eq:GammaUproportion}
    \ell_v(\Gamma_U) \geq (1-\alpha) \ell_v(\Gamma).
\end{equation}
Observing Figure \ref{fig:sigma1b}, the set of partition elements which $\Gamma_U$ can intersect is given by $\mathcal{A}_U = \{ A_4 \cap \sigma_1, A_{4,3}^k \, | \, k \geq 1 \}$, so $\mathcal{M}_U = \{ M_4, M_4M_3^k \, | \, k \geq 1 \}$. Note that any two element subset $\mathcal{N} \subset \mathcal{M}_U$ satisfies
\begin{equation}
    \begin{split}
        \sum_{M \in \mathcal{N}} \frac{1}{ K_+(M)} & \leq \frac{1}{K_+(M_4)} + \frac{1}{K_+(M_4M_3)} \\
        & = \frac{21}{79} + \frac{21}{127} = \beta \approx 0.431
    \end{split}
\end{equation}
and $\alpha + \beta <1$. It follows that if $\Gamma_U$ intersects two or fewer of the elements of $\mathcal{A}_U$, we can guarantee (C2) by the standard method, summing the reciprocals of expansion factors. Assume, then, that $\Gamma_U$ intersects three or more elements from $\mathcal{A}_U$. It follows that
\begin{enumerate}[label={($\ddag$)}]
    \item $\Gamma_U$ traverses $A_{4,3}^k$ for some $k \geq 1$, connecting the lines $\mathcal{L}_k:$ $y=\frac{(4k+2)x + k+2}{4k+4}$ and $\mathcal{L}_{k-1}:$ $y=\frac{(4k-2)x +k+1}{4k}$.
\end{enumerate}
We now follow a similar inductive argument to before, assuming that $\Gamma$ violates the lemma and aiming to contradict ($\ddag$). Let $(x_k,y_k) = \left( \frac{k+2}{4k+6}, \frac{k+2}{2k+3} \right)$ denote the intersections of the lines $\mathcal{L}_k$ with the boundary $y=2x$ of $\sigma$. Assume $\Gamma$ traverses $A_{4,3}^k$, write its restriction to this set as $\Gamma^k$. Since the gradients of the $\mathcal{L}_k$ are monotonic increasing in $k$ and vectors in $\mathcal{C}_+$ have gradients bounded above by 3, A lower bound on $\ell_v(\Gamma_k)$ is given $h_k = y_k'-y_k$, where $(x_k',y_k')$ is the intersection of the line $y-y_k = 3(x-x_k)$ and $\mathcal{L}_{k-1}$, in particular
\begin{equation}
    \label{eq:sigma1bhk}
    h_k = \frac{8k^2+18k+7}{16k^2+28k+6} - \frac{k+2}{2k+3} = \frac{3}{16k^2+28k+6}.
\end{equation}
For the base case suppose that $\Gamma_U$ traverses $A_{4,3}^1$. Let $(x_U,y_U)$ be the intersection with $y=1/2 + x/2$, the boundary between $A_{4,3}^1$ and $A_4 \cap \sigma_1$. Note that this point maps to $(1,y_U)$ under $H$ with $y_U<2/3$. Figure \ref{fig:sigma1b} shows the preimage $\mathcal{P}$ in $A_4 \cap \sigma_1$ of the segment joining $(1/2,3/4)$ to $(1,1)$ between $A_3$ and $A_4$. Specifically $\mathcal{P}$ lies on the line $y = 7/12 + 5x/12$ and $H(\mathcal{P})$ lies on $y = 1/2 + x/2$. If $\Gamma$ intersects $\mathcal{P}$, then $H(\Gamma)$ connects $(1,y_U)$ to a point on the segment joining $(1/2,3/4)$ to $(1,1)$. Since $y_U< 3/4$, it follows that $H(\Gamma)$ traverses $A_3$, making non-simple intersection with $A_4$, so that (C1) is satisfied. Assume, then, that $\Gamma_U$ does not intersect $\mathcal{P}$. This gives an upper bound $\ell_v(\Gamma_U) \leq y_0 - 1/2 =: L_0$, where $(x_0,y_0) = (7/19,14/19)$ is the intersection of $\mathcal{P}$ with the boundary of $\sigma_{1}$ on $y=2x$ (see Figure \ref{fig:sigma1b}). Noting (\ref{eq:GammaUproportion}), (C2) follows with $\Lambda = H^2\left(\Gamma^1\right)$ if the inequality $K_+(M_4M_3) h_1 > L_0/(1-\alpha)$ is satisfied. Indeed
\[ \left(\frac{16}{7} + \frac{79}{21}\right) \frac{3}{16+28+6} - \frac{9}{38(1-\alpha) } \approx 0.00277 > 0 \]
so that the base step of the induction holds. The inductive step is roughly analogous, reducing to checking the inequality 
\begin{equation}
\label{eq:sigma1bInduction}
    K_+\left(M_4M_3^k\right) h_k - \frac{L_{k-1}}{1-\alpha} >0,
\end{equation} where $L_{k-1} = y_{k-2} - 1/2$ is the height of the partition element $A_{4,3}^{k-1}$. One can verify that this inequality holds (the function is monotonic decreasing in $k\geq 2$ with limit $0$ as $k \rightarrow \infty$), establishing the lemma for $\Gamma \subset \sigma_{1b}$.

\begin{figure}
    \centering

    \begin{tikzpicture}[scale=2.3]
    \draw (190/89,215/89) -- (230/89,190/89) --  (5,6) -- (5,7) -- (190/89,215/89);
    
    \definecolor{tomato}{RGB}{255, 99, 71}
    \definecolor{teal}{RGB}{95, 158, 160} 
    
    \fill[fill=teal,opacity=0.8] (190/89, 430/178) -- (230/89,190/89) -- (10/3,10/3) -- (20/7,50/14) -- (190/89, 430/178);
    
    \node at (3.226,3.271) {$A_{4,2}^2$};
    
    \node at (4.36,5.08) {$A_{1,3}^1$};
    \node at (2.75,2.89) {$A_{4,2}^1$};
    
    \draw[->] (2.45,3.237) -- (2.673,3.237);
    \node at (2.35,3.21) {$A_{4,2}^2$};

    \draw[->] (4.4,6.62) -- (4.882,6.62);
    \node at (4.3,6.6) {$A_{2,3}^2$};
    
    \draw[->] (4.6,6.8) -- (4.92,6.8);
    \node at (4.5,6.78) {$A_{1,3}^2$};
    
    \clip (190/89,215/89) -- (230/89,190/89) --  (5,6) -- (5,7) -- (190/89,215/89);
    \draw (0,5) -- (10,5);
    \draw (5,10) -- (5,0);
    
    \draw (0,0) rectangle (10,10);
        
        \foreach \k in {1,...,45}{
        \draw ({ (2*\k+2)*(7.5-7.5)/(2*\k+1)  }, { 7.5 })--({ (2*\k+2)*(10-7.5)/(2*\k+1)  }, { 10 });
        \draw[blue] ({ 5+ (2*\k+2)*(5-7.5)/(2*\k+1)  }, { 5 })--({ 5+ (2*\k+2)*(7.5-7.5)/(2*\k+1) }, { 7.5 });
        \draw ({ 2.5 + (2*\k+1)*(5-5)/(2*\k)  }, { 5 })--({ 2.5 + (2*\k+1)*(10-5)/(2*\k)  }, { 10 });
        \draw ({ 12.5 + (2*\k+1)*(5-10)/(2*\k)  }, { 5 }) -- ({ 12.5 + (2*\k+1)*((30*\k+20)/(4*\k+2)-10)/(2*\k)  }, { (30*\k+20)/(4*\k+2) });
        }
        
        \draw (30/4,70/8) -- (4,6);
        \draw (70/8,150/16) -- (7,8);
        \draw (15-30/4,15-70/8) -- (15-5,15-190/28);
        \draw (15-70/8,15-150/16) -- (15-7,15-8);
        \draw (2.5,10) -- (0 , 50/6);
        \draw (1,9) -- (0,15-190/28);

        \draw (5,7)--(90/14,230/28);
        \draw (15-5,15-7)--(15-90/14,15-230/28);
        
        \filldraw[fill=white] (0,5) -- (10,10) -- (5,10) -- (0,7.5) -- (0,5);
        \filldraw[fill=white] (5,5) -- (10,7.5) -- (10,5) -- (5,5);

        \draw (0,150/16) -- (1,10);
        \draw (0,70/8) -- (2,10);
        \draw (10,150/16) -- (3,5);
        \draw (10,70/8) -- (4,5);
      
       \foreach \k in {1,...,45}{
        \draw ({10- (2*\k+2)*(7.5-7.5)/(2*\k+1)  }, { 7.5-5 })--({10- (2*\k+2)*(10-7.5)/(2*\k+1)  }, { 10-5 });
        \draw ({ 5 - (2*\k+2)*(5-7.5)/(2*\k+1)  }, { 5-5 })--({ 5 - (2*\k+2)*(7.5-7.5)/(2*\k+1)  }, { 7.5-5 });
        \draw ({ 7.5 - (2*\k+1)*(5-5)/(2*\k)  }, { 5-5 })--({ 7.5- (2*\k+1)*(10-5)/(2*\k)  }, { 10-5 });
        \draw[red] ({ -2.5 - (2*\k+1)*(5-10)/(2*\k)  }, { 5-5 }) -- ({ -2.5 - (2*\k+1)*((30*\k+20)/(4*\k+2)-10)/(2*\k)  }, { (30*\k+20)/(4*\k+2) -5 });
        }
        
        \draw (10-30/4,70/8-5) -- (10-4,6-5);
        \draw (10-70/8,150/16-5) -- (10-7,8-5);
        \draw (-5+30/4,15-70/8-5) -- (0,15-190/28-5);
        \draw (-5+70/8,15-150/16-5) -- (-5+7,15-8-5);
        \draw (10-2.5,10-5) -- (10, 50/6-5);
        \draw (10-1,9-5) -- (10,15-190/28-5);

        \draw (10-5,7-5)--(10-90/14,230/28-5);
        \draw (0,15-7-5)--(-5+90/14,15-230/28-5);
        
        \filldraw[fill=white] (10,0) -- (0,5) -- (5,5) -- (10,2.5) -- (10,0);
        \filldraw[fill=white] (0,0) -- (5,0) -- (0,2.5) -- (0,0);
        \draw (5,6) -- (10/8,0);
        
        \draw (10,150/16-5) -- (10-1,10-5);
        \draw (10,70/8-5) -- (10-2,10-5);
        \draw (0,150/16-5) -- (10-3,5-5);
        \draw (0,70/8-5) -- (10-4,5-5);

        \node at (2.5,2.5) {$A_{3,2}^1$};

        \node at (2.935,3.2) {$A_{3,2}^2$};

        \node at (3.6,4.245) {$A_1 \cap \sigma_3$};

        \node at (4.4,5.5) {$A_{2,3}^1$};
        \node at (4.769,6.27) {$A_{1,3}^1$};
      
    \end{tikzpicture}
    
      \caption{Behaviour of $H_\sigma$ over $\sigma_{3b} \setminus \varsigma_3$ and $H_\sigma^2$ over $\varsigma_3$, shaded in blue.}
    \label{fig:sigma3b}
\end{figure}
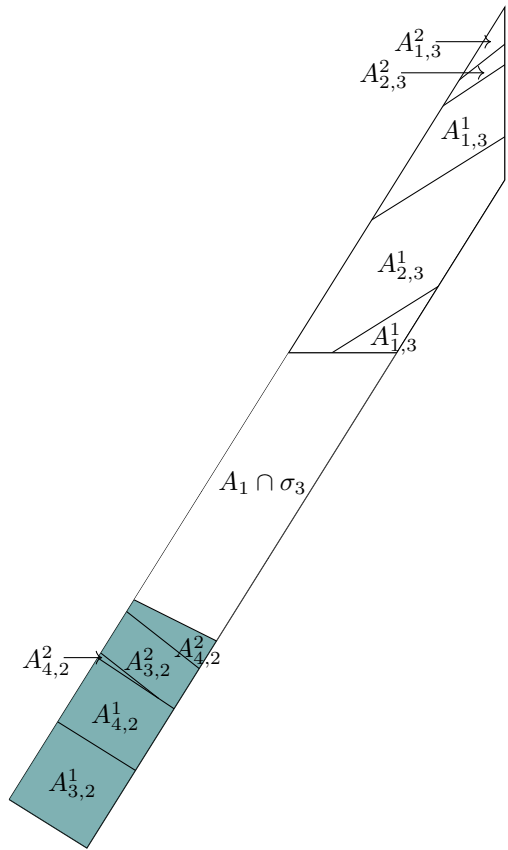
    
Next consider $\Gamma \subset \sigma_{3b}$, shown in Figure \ref{fig:sigma3b}. Note that outside of $\varsigma_3$ (shaded in blue) the Jacobian $DH_\sigma$ is some matrix from $\mathcal{M}$, but over $\varsigma_3$ we have $DH_\sigma = M_2 \notin \mathcal{M}$. Therefore if we are to expand from some subset of $\Gamma \cap \varsigma_3$, to ensure that $\Lambda$ satisfies one of (A1-2) we must map forwards using $H_\sigma^2$, whose Jacobian is always a matrix from $\mathcal{M}$ (analogous to the escape behaviour shown for $A_3$, shown in Figure \ref{fig:escapeBehaviour}). The relevant subset of matrices, then, is $\mathcal{M}_{3b} = \{M_1, M_1M_3^k, M_2M_3^k, M_3M_2^k, M_4M_2^k \,| \, k=1,2 \}$. Noting that $\Gamma$ can have non-simple intersection with the sets $A_{4,2}^2$ and $A_{1,3}^1$, the relevant inequality to verify is
\[ \left (\sum_{M \in \mathcal{M}_{3b}} \frac{1}{K_+(M)} \right ) + \frac{1}{K_+(M_1M_3)} + \frac{1}{K_+\left(M_4M_2^2\right)} < 1. \]
Indeed, the above sums to $\delta \approx 0.807 <1$, so that restricting to one of the partition elements and expanding from there (using $H_\sigma^2$ inside of $\varsigma_3$, $H_\sigma$ otherwise) will always satisfy (C2) with some $k \leq 3$. This leaves the cases $\Gamma \subset \sigma_{1c},\sigma_{1d},\sigma_{3a}$. Noting that rotating $\sigma_{1c}$ by $180^\circ$ about the point $(1/4,3/4)$ gives $\sigma_{1b}$, and $\mathcal{C}_+$ is invariant under this rotation, the argument is essentially analogous. Similarly the arguments for $\sigma_{1d},\sigma_{3a}$ are equivalent to those for $\sigma_{1a},\sigma_{3b}$ respectively. This concludes the case where $\Gamma$ satisfies (A1).

Let $\Gamma$ satisfy (A2). Define the transformation $T:\tor \rightarrow \tor$ given by $T(x,y) = (1-x,y+1/2) \text{ mod }1$. One can verify that $T \circ T = \mathrm{Id}$ and $T \circ H = H \circ T$ so that $H^n = T \circ H^n \circ T$. Now since $T(\sigma_2 \cup \sigma_4) = \sigma_1 \cup \sigma_3$ and $DT \mathcal{C}_- = \mathcal{C}_+$, the line segment $T(\Gamma)$ satisfies (A1). By our analysis above, $T(\Gamma)$ then satisfies (C1) or (C2). Noting that $T(A_j) = A_{5-j}$, if $T(\Gamma)$ satisfies (C1) then there exists $k$ such that $H^k(\Gamma) = (T \circ H^k \circ T)(\Gamma)$ has non simple intersection with $A_{5-j}$, so (C1) is satisfied. Similarly since $\ell_v(\cdot)$ is $T$-invariant, if $T(\Gamma)$ satisfies (C2) then the same holds for $\Gamma$.
\end{proof}

\subsection{Establishing the Bernoulli property}

We are now ready to establish the mixing property.

\begin{proof}[Proof of Theorem \ref{thm:mixingOTM}]
By Theorem \ref{thm:katok-strelcyn}, with \textbf{(KS1-2)} shown in \cite{myers_hill_exponential_2022} and \textbf{(KS3)} in Proposition \ref{prop:tentMap}, it suffices to show \textbf{(MR)}. By Lemmas \ref{lemma:OTMalignment}, \ref{lemma:mappedOTMalignment}, for a.e. $z$ we can find $m_0$ such that $H^{m_0}(\gamma_u(z))$ contains a line segment $\Gamma_0$ satisfying (A1) or (A2). Now iteratively apply Lemma \ref{lemma:unstableGrowth} until (C1) is satisfied, giving $m_1$ such that $H^{m_1}(\Gamma_0)$ contains a line segment $\Gamma_1$ which has non simple intersection with some $A_j$. Define a $v$-segment as any line segment traversing $S_1$, connecting its upper and lower boundaries. Similarly define a $h$-segment as any line segment in $S_1$ which connects its left and right boundaries. Consider the four parallelograms $Q_j \subset A_j$ given by $Q_1 = A_1 \cap S_2$, $Q_2 = A_2 \cap S_1$, $Q_3 = A_3 \cap S_4$, $Q_4 = A_4 \cap S_3$. We recall from \cite{myers_hill_exponential_2022}, specifically the proof of Lemma 4 in that work, that
\begin{enumerate}[label={(M\arabic*)}]
    \item If $\Gamma_1$ has non-simple intersection with some $A_{j_1}$, it traverses some $Q_{j_2}$, connecting its sloping boundaries.
    \item If $\Gamma_1$ traverses $Q_{j_2}$, $j_2 = 1,2,3,4$, then $H^k(\Gamma_1)$ traverses $Q_3$ for $k=2,1,0,3$ respectively.
    \item The image of any line segment traversing $Q_3$ contains a $v$-segment.
\end{enumerate}
The above gives $m_2 \in \{ 1,2,3,4\}$ such that $H^{m_2}(\Gamma_1)$ contains a $v$-segment $\Gamma \subset H(Q_3) \cap S_1$, with this parent set given by the quadrilateral with corners $(0,0)$, $(1/6,0)$, $(1/2,1/2)$, $(1/3,1/2)$, so that $\Gamma$ connects points $(x_1,0)$ and $(x_2,1/2)$ with $0 \leq x_1 \leq 1/6$ and $1/3 \leq x_2 \leq 1/2$. It follows that $\Gamma$ traverses $Q_2$ which, by (M2-3), implies that $H^2(\Gamma)$ contains a $v$-segment and so does $H^{2+2k}(\Gamma)$ for $k \geq 0 $ by induction. Applying $F$ to $\Gamma$ has no effect on $(x_1,0)$ and wraps $(x_2,1/2)$ horizontally around the torus so that $F(\Gamma)$ contains a segment joining $(0,y)$ to $(x_2,1/2)$ with $y<1/2$. Now $G$ has no effect on $(0,y)$ and maps $(x_2,1/2)$ to $(x_2,1/2+2x_2) \text { mod }1$. Since $1/2 + 2x_2 \geq 1/2 + 2/3 = 7/6 >1$, $H(\Gamma)$ contains a segment joining $(0,y)$ to $(x_3,1)$ with $x_3 \leq x_2 \leq 1/2$. It follows that $H(\Gamma)$ must traverse $Q_4$ which, by (M2-3), implies that $H^5(\Gamma)$ contains a $v$-segment. Using the same induction as before we have that $H^{5+2k}(\Gamma)$ contains a $v$-segment for all $k \geq 0$ which, together with the same result for $H^{2+2k}(\Gamma)$, implies that $H^k(\Gamma)$ contains a $v$-segment for all $k\geq 4$. Hence there exists $M = m_0 + m_1 + m_2 + 4$ such that $H^m(\gamma_u(z))$ contains a $v$-segment for all $m \geq M$.

Now for almost any $z'$, by Lemma \ref{lemma:OTMalignment} we can find $n_0$ such that $H^{-n_0}(\gamma_s(z'))$ contains a line segment $\Gamma'\in \sigma'$, aligned with some $v\in \mathcal{C}'$. Define the transformation $\mathcal{T}(x,y) = (1-y,1-x) \text{ mod }1$. One can verify that $\mathcal{T} \circ \mathcal{T} = \mathrm{Id}$ and $\mathcal{T} \circ H^{-1} = H \circ \mathcal{T}$ so that $H^{-k} = \mathcal{T} \circ H^k \circ \mathcal{T}$. Now since $\mathcal{T}(\sigma') = \sigma$ and $D\mathcal{T}\mathcal{C}'= \mathcal{C}$, we have that $\mathcal{T}(\Gamma')$ is a line segment in $\sigma$, aligned with some $v = D\mathcal{T}v' \in \mathcal{C}$. We now follow Lemmas \ref{lemma:mappedOTMalignment}, \ref{lemma:unstableGrowth} and the argument above to find $n_1$ such that $(H^{m}\circ \mathcal{T})(\Gamma')$ contains a $v$-segment for all $m \geq n_1$. The image of a $v$-segment under $\mathcal{T}$ is a segment joining the left and right boundaries of $S_4$. Noting Figure \ref{fig:firstPartitions}, we have that $H^{-m}(\Gamma') = (\mathcal{T} \circ H^m \circ \mathcal{T})(\Gamma')$ traverses the parallelogram $Q_2' = A_2' \cap S_4$, connecting its sloping boundaries. It was shown in \cite{myers_hill_exponential_2022} that if $\Gamma'$ traverses $Q_2'$ then $H^{-1}(\Gamma)$ contains a $h$-segment, so that $H^{-n}(\gamma_s(z'))$ contains a $h$-segment for all $n \geq N = n_0 + n_1 + 1$. Since $z$ and $z'$ were arbitrary and $h$-segments and $v$-segments must always intersect, \textbf{(MR)} holds.
\end{proof}

\begin{remark}
The $v$-segments $\Lambda$ obtained above satisfy $H^{-1}(\Lambda) \subset Q_3$, $H^{-2}(\Lambda) \subset Q_2$ so that $\Lambda \subset H(A_3 \cap H(A_2)) = \sigma_3$. Similarly the $h$-segments derived from these $v$-segments can be shown to lie in $\sigma_2'$.
\end{remark}

\section{Decay of correlations for the return map}
\label{sec:Hsigma}
As described in the introduction, we infer the polynomial decay under $H$ from exponential decay of some induced return map $H_A$, where returns to $A$ experience `strong' hyperbolic behaviour. The natural choice for $A$, following the work of section \ref{sec:Hmixing}, is the set $\sigma$. We begin by proving the Bernoulli property for $H_\sigma$.

\subsection{Bernoulli property}

\begin{prop}
\label{prop:HsigmaBernoulli}
The return map $H_\sigma$ is Bernoulli with respect to the probability measure $\mu_\sigma = \mu(\sigma)^{-1} \mu $.
\end{prop}

We will show the conditions \textbf{(KS1-3)} and \textbf{(MR)}; the result then follows from Theorem \ref{thm:katok-strelcyn}.

\begin{lemma}
\label{lemma:HsigmaKatokStrelcyn}
The return map $H_\sigma$ satisfies \textbf{(KS1-3)}.
\end{lemma}

\begin{proof}
Starting with \textbf{(KS1)} we follow a similar approach to \cite{springham_polynomial_2014}, their Lemma 4.1. We show that there exists $a,C_1>0$ s.t. $\forall \, \epsilon>0$, $\mu_\sigma(B_\varepsilon(S)) \leq C_1 \varepsilon^a$ for $S = \mathcal{S} \cap \sigma_{1a}$; the argument for the rest of $\mathcal{S}$ is similar and the result then follows by taking a larger $C_1$. Recall the line segments $\mathcal{L}_k$ from ($\dagger$) which for $k \geq 3$ terminate on the points $\left( 0 , (k+1)/(4k+2) \right)$ and $\left( 1/(4k-2) , (k-1)/(4k-2) \right)$ on the line $\mathcal{L}: y= 1/4 -x/2$. Let $P(\varepsilon)$ denote the parallelogram in $\sigma_{1a}$ of width $2\sqrt{\varepsilon}$, height $\sqrt{\varepsilon}$, with sides aligned with $x=0$ and $\mathcal{L}$ (see Figure \ref{fig:sigma1a}). For small $\varepsilon$, $P(\varepsilon)$ contains all line segments $\mathcal{L}_k$ where $2\sqrt{\varepsilon} \geq 1/(4k-2)$, i.e. $k \geq k_0 = \lceil 1/(8\sqrt{\varepsilon}) + 1/2 \rceil$, with
\[ \mu(B_\varepsilon(P(\varepsilon))) = (2\sqrt{\varepsilon} + 2\varepsilon )(\sqrt{\varepsilon} + 2\varepsilon) = 2\varepsilon + 6 \varepsilon^{3/2} + 4\varepsilon^2 < 12\varepsilon. \]
The ball $B_\varepsilon(P(\varepsilon))$ then covers all of $B_\varepsilon(S)$ except the collection $\mathcal{L}_k$, $4 \leq k \leq k_0-1$ and the seven line segments $L_j$ which terminate on $y=2x$. The measure of the ball around these latter line segments satisfies
\[ \mu(B_\varepsilon(\cup_j L_j)) \leq 14 \, \varepsilon \left(\max_j |L_j| + 2\varepsilon \right) < c_1 \varepsilon  \]
for some finite $c_1$, so it remains to estimate $ \sum_{k=4}^{k_0-1} \mu(B_\varepsilon(\mathcal{L}_k))$. We can calculate
\begin{equation}
    \label{eq:Lklength}
    |\mathcal{L}_k| = \sqrt{ \left( \frac{1}{4k-2} \right)^2 + \left( \frac{k+1}{4k+2} - \frac{k-1}{4k-2} \right)^2 } = \sqrt{\frac{8k^2 + 4k+1}{4(4k^2-1)^2 }} < \frac{1}{k} 
\end{equation}
so that
\[ \sum_{k=4}^{k_0-1} \mu(B_\varepsilon(\mathcal{L}_k)) < 2\varepsilon \sum_{k=4}^{k_0-1} \frac{1}{k} + 2\varepsilon < 4\varepsilon^2k_0 + 2\varepsilon \log k_0 < c_2 \varepsilon^a \]
for some $0<a<1$, $c_2>0$ since $k_0 < \varepsilon^{-1/2}$ and there exists finite $c$ such that $c \, \varepsilon^a> \varepsilon \log \frac{1}{\varepsilon}$ for any $0<a<1$.

Since $H_\sigma$ is piecewise linear, condition \textbf{(KS2)} follows trivially and we move onto \textbf{(KS3)}. Existence of Lyapunov exponents almost everywhere follows from Oseledets' theorem \cite{oseledets_multiplicative_1968} provided that $\max \{\log\|DH_\sigma\| ,0 \}$ is integrable. This follows from the fact that if $z \in \sigma$ has return time $R(z;H,\sigma) =k$, then the Jacobian of $H_\sigma$ at $z$ satisfies $\| DH_\sigma \| \leq c_1 k$ for some finite $c_1>0$, and that the measure of the sets $\{ z \in \sigma \, | \, R(z;H,\sigma) =k \}$ are of order $k^{-3}$. That these Lyapunov exponents are non-zero follows from Lemmas \ref{lemma:itineraries}, \ref{lemma:cone} and an argument similar to that given for $H$ in section \ref{sec:OTMhyp}.

\begin{lemma}
\label{lemma:HsigmaMR}
The return map $H_\sigma$ satisfies \textbf{\textbf{(MR)}}.
\end{lemma}

For a.e $z \in \sigma$, local manifolds $\gamma_u(z)$, $\gamma_s(z)$ under $H_\sigma$ align with those of $H$. Note that $H_\sigma$ does not immediately inherit \textbf{(MR)} from $H$ as while successive images of local manifolds under $H$ contain $h$-segments and $v$-segments, these segments may not lie in the successive images under $H_\sigma$.

\begin{figure}
    \centering
    \begin{tikzpicture}[scale=0.8]

    \definecolor{tomato}{RGB}{255, 99, 71}
    \definecolor{teal}{RGB}{95, 158, 160} 
    
    \draw (5,0) -- (5,10);

    \draw[very thick, tomato] (2,2.8) -- (3,2.2);
    
    \draw[very thick, teal] (10-2,2.8+5) -- (10-3,2.2+5);

    \filldraw[fill=gray!40] (0,8) -- (10/8,10) -- (2.5,10) -- (0,5) -- (0,8);
      \filldraw[fill=gray!40] (0,9) -- (0,10) -- (10/16,10) -- (0,9);
      \filldraw[fill=gray!40] (0,0) -- (5,10) -- (5,7) -- (10/16,0) -- (0,0);
      \filldraw[fill=gray!40] (10/8,0) -- (2.5,0) -- (5,5) -- (5,6) -- (10/8,0) ;
      \filldraw[fill=gray!40] (5,2) -- (50/8,0) -- (7.5,0) -- (5,5) -- (5,2);
      \filldraw[fill=gray!40] (5,1) -- (5,0) -- (90/16,0) -- (5,1);
      \filldraw[fill=gray!40] (5,10) -- (10,0) -- (10,3) -- (90/16,10) -- (5,10);
      \filldraw[fill=gray!40] (50/8,10) -- (7.5,10) -- (10,5) -- (10,4) -- (50/8,10);

    \filldraw[fill=gray!40] (0,150/16) -- (0,10) -- (1,10) -- (0,150/16);
      \filldraw[fill=gray!40] (0,70/8) -- (2,10) -- (5,10) -- (0,7.5);
      \filldraw[fill=gray!40] (0,5) -- (10,10) -- (10,150/16) -- (3,5) -- (0,5);
      \filldraw[fill=gray!40] (10,70/8) -- (10,7.5) -- (5,5) -- (4,5) -- (10,70/8) ;
   
     \filldraw[fill=gray!40] (10,150/16-5) -- (10,10-5) -- (10-1,10-5) -- (10,150/16-5);
      \filldraw[fill=gray!40] (10,70/8-5) -- (10-2,10-5) -- (10-5,10-5) -- (10,7.5-5);
      \filldraw[fill=gray!40] (10,5-5) -- (0,10-5) -- (0,150/16-5) -- (10-3,5-5) -- (10,5-5);
      \filldraw[fill=gray!40] (0,70/8-5) -- (0,7.5-5) -- (10-5,5-5) -- (10-4,5-5) -- (0,70/8-5) ;   

    \draw (0,5) -- (10,5);
     \draw(0,0) rectangle (10,10);

    \end{tikzpicture}
    \caption{The set $\sigma'$ superimposed on $\sigma$, their intersection left white. A $\mathfrak{h}$-segment (red) and a $\mathfrak{h}'$-segment (blue) are plotted in $\mathcal{R}=\sigma_3 \cap \sigma_2'$ and $\mathcal{R}' = \sigma_2 \cap \sigma_3'$ respectively.}
    \label{fig:frakh}
\end{figure}
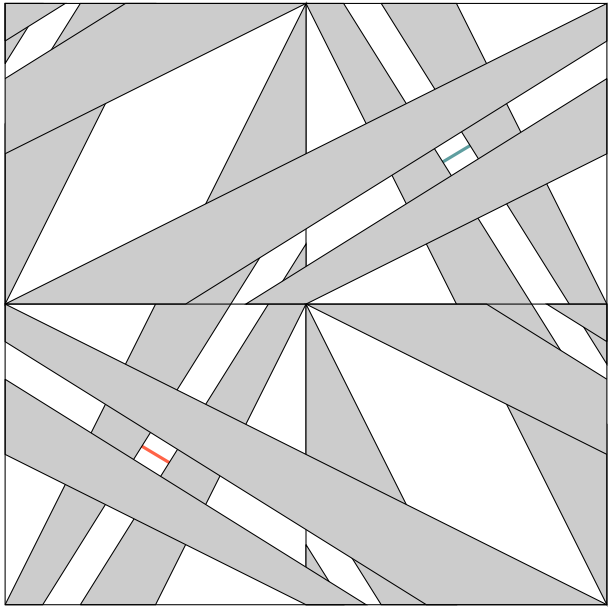

Let $\mathcal{R}$ denote the quadrilateral $\sigma_3 \cap \sigma_2'$ and $\mathcal{R}' = \sigma_2 \cap \sigma_3'$. Define a $\mathfrak{h}$-segment as a line segment spanning $\mathcal{R}$ with endpoints on $\partial \sigma_3$. Similarly define a $\mathfrak{h}'$-segment as a line segment spanning $\mathcal{R}'$ with endpoints on $\partial \sigma_2$. Examples are plotted in Figure \ref{fig:frakh}. We will show that there exists $M,N$ such that for all $m \geq M$, $n \geq N$, $H_\sigma^m(\gamma_u(z))$ intersects $H_\sigma^{-n}(\gamma_s(z'))$ in either $\mathcal{R}$ or $\mathcal{R}'$. 

By the remark after Theorem \ref{thm:mixingOTM} we can find some $n_2$ such that $H^{-n_2}(\gamma_s(z'))$ contains a $h$-segment in $\sigma_2'$, which in turn contains a $\mathfrak{h}$-segment in $\mathcal{R}$. As a line segment in $\sigma$, this $\mathfrak{h}$-segment lies in $H_\sigma^{-n_1}(\gamma_s(z'))$ for some $n_1 \leq n_2$. Note that we have a hyperbolic period 2 orbit $(1/4,1/4) \leftrightarrow (3/4,3/4)$ under $H_\sigma$, alternating between $\mathcal{R}$ and $\mathcal{R}'$. Any $\mathfrak{h}$-segment $\Lambda$ contains a point $\zeta$ on the unstable manifold through $(1/4,1/4)$ so $H_\sigma^{-2}(\Lambda)$ contains a point $\zeta'$ on the manifold closer to $(1/4,1/4)$ and extends beyond the boundaries $\partial \sigma_3$ by the expansion of $H_\sigma^{-2}$. Hence $H_\sigma^{-2}(\Lambda)$ contains a $\mathfrak{h}$-segment and by induction so does $H_\sigma^{-2k}(\Lambda)$ for all $k\geq 1$. The odd iterates $H^{-2k+1}(\Lambda)$ similarly span $\mathcal{R}'$ so that
\begin{enumerate}[label={($\vee$)}]
    \item Given arbitrary $z' \in \sigma$, there exists $N$ such that for all $n \geq N$ the image $H_\sigma^{-n}(\gamma_s(z'))$ contains a $\mathfrak{h}$-segment or a $\mathfrak{h}'$-segment.
\end{enumerate}
Define a $v'$-segment as a line segment vertically spanning $\sigma_2 \cap S_4$. Condition \textbf{(MR)} now follows from establishing
\begin{enumerate}[label={($\wedge$)}]
    \item Given arbitrary $z \in \sigma$, there exists $M$ such that for all $m \geq M$ the image $H_\sigma^{m}(\gamma_s(z'))$ contains both a $v$-segment \emph{and} a $v'$-segment.
\end{enumerate}
Recall the quadrilaterals $\varsigma_j \subset \sigma_j$ of points with return time 1 (see Figure \ref{fig:returnTime1}). It follows from the definitions of the $\sigma_j$ that $H(\varsigma_2) \subset \sigma_3$ and $H(\varsigma_3) \subset \sigma_2$. The edges of $\varsigma_3$ on the $A_1,A_2$ boundary map into the lines $x=1/2$, $x=1$, in particular onto the red dashed lines on the boundary of $S_2$ in Figure \ref{fig:HSigma} so that the image of any line segment in $\varsigma_3$ which joins these edges contains a $v'$-segment. A analogous result holds for lines segments traversing $\varsigma_2$ and since this behaviour occurs within the return set $\sigma$ we have that $v$-segments map into $v'$-segments under $H_\sigma$ and vice versa. It follows that $v$(')-segments map into $v$(')-segments under $H_\sigma^2$. It suffices to break into the odd iterates to satisfy the `and' condition ($\wedge$).

By following the steps (M1-3) in the proof of Theorem \ref{thm:mixingOTM} we can find $m_2$ such that $H^{m_2}(\gamma_u(z))$ contains a $v$-segment $\Gamma \subset \sigma_3 \cap S_1$ with $H^{-1}(\Gamma) \subset A_3$, $H^{-2}(\Gamma) \subset A_2$, $H^{-3}(\Gamma) \subset A_1$, $H^{-4}(\Gamma) \subset A_4$. In particular $H^{-1}(\Gamma)$ lies in $H(A_2 \cap H(A_1))$, i.e. outside of $\sigma_2$, so that $\Gamma$ lies in $\sigma_3 \setminus (H(\sigma_2) \cap \sigma_3)$. The set $H(\sigma_2) \cap S_1$ (shown in blue in Figure \ref{fig:HsigmaMR}) is the quadrilateral with corners $(7/68,0)$, $(3/34,0)$, $(27/68,0)$, $(7/17,0)$, which splits $\sigma_3 \cap S_1$ into left and right parts. We assume first that $\Gamma$ lies in the right part, intersecting the line $y=1/2$ at some point $(x_1,1/2)$ with $x_1 \geq 7/17$ and the $A_1,A_2$ boundary $y=1/2-x/2$ at some point $(1-2y_1,y_1)$. These intersections define a line segment $\Gamma_1 \subset \Gamma$, which lies in $A_1$ shown in Figure \ref{fig:HsigmaMR}. Applying $F$ maps $(x_1,1/2)$ to itself (wrapping horizontally around the torus) and maps $(1-2y_1,y_1)$ to $(0,y_1)$. Applying $G$ then leaves $(0,y_1)$ invariant and wraps $(x_1,1/2)$ vertically around the torus to $(x_1, 1/2 + 2x_1) \text{ mod }1 \equiv (x_1,-1/2+2x_1)$. Since $x_1 \geq 7/17$ we have that $-1/2 + 2x_1 \geq 11/34 > 10/24 \geq 1/2 - x_1/2$ so that $(x_1,-1/2+2x_1)$ lies above the line $y=1/2-x/2$. We restrict again to $A_1$, giving $\Gamma_2 \subset H(\Gamma_1)$ with endpoints on $(x_1,-1/2+2x_1)$ and some point $(1-2y_2,y_2)$ on $y=1/2-x/2$. This line meets $y=-1/2+2x$ at $(2/5,3/10)$ so that $y_2 \leq 3/10$ (see Figure \ref{fig:HsigmaMR}). Now $F(\Gamma_2)$ joins $(0,y_2)$ to $(5x_1 -2, 2x_1 - 1/2)$ and so $H(\Gamma_2)$ joins $(0,y_2)$ to $(5x_1-2,12x_1-9/2)$. The set $A_{3,2}^1$ is bounded by the parallel lines $y=7/16 -5x/8$ and $y=3/8 - 5x/8$. Since $y_2<3/8$ and $12x_1-9/2 \geq 15/34 > 109/272 \geq 7/16 - 5(5x_1 -2)/8$ we have that $H(\Gamma_2)$ contains a segment $\Gamma_3$ which traverses $A_{3,2}^1$. The image $H(\Gamma_3)$ then traverses $Q_3$ so that $H^2(\Gamma_3) \subset H^4(\Gamma)$ contains a $v$-segment in $\sigma_3$. Critically we have that $\Gamma_2,\Gamma_3 \subset H(A_1) \subset \sigma$ but $\Gamma_4 \subset H(A_2 \cap H(A_1))$ which is not in $\sigma$. Hence $H_\sigma^3(\Gamma)$ contains a $v$-segment and we can apply the $H_\sigma^2$ result above to show that $H_\sigma^k(\Gamma)$ contains a $v$-segment for all $k \geq 2$. Similar analysis can be applied to $\Gamma$ in the left portion of $\sigma_3 \setminus (H(\sigma_2) \cap \sigma_3)$. It follows that $H_\sigma^k(\Gamma)$ contains a $v'$-segment for all $k \geq 3$, establishing ($\wedge$) with $M=m_1+3$.
\end{proof}
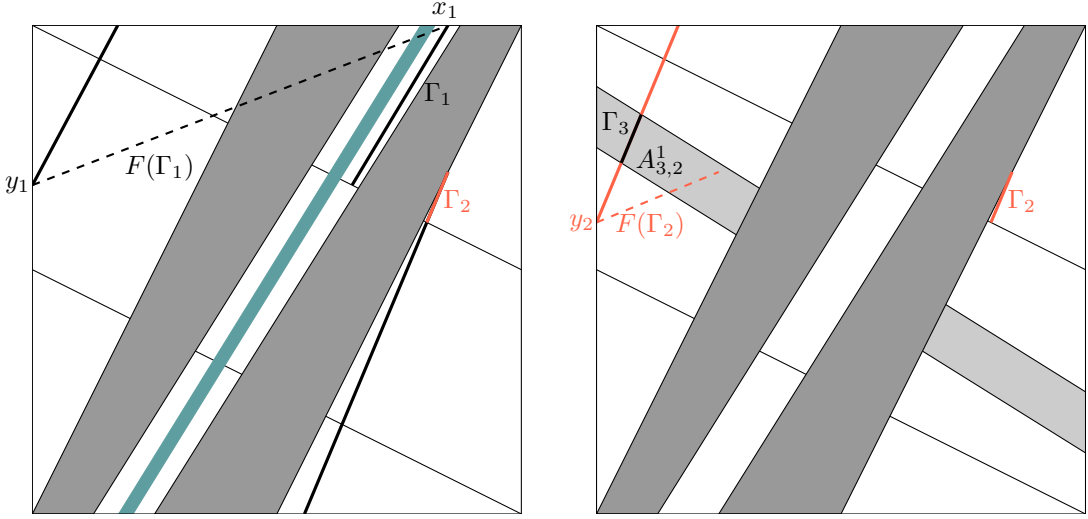
\begin{figure}
    \centering

    \begin{tikzpicture}

    \definecolor{tomato}{RGB}{255, 99, 71}
    \definecolor{teal}{RGB}{95, 158, 160} 
    
    \node at (-3.5,0) {
  
    \begin{tikzpicture}[scale=1.3]
    \node at (-0.15,350/104) {$y_1$};
    \clip (0,0) rectangle (5,5);
    
    \fill [fill=white] (0,0) -- (2.5,5) -- (5,5) -- (2.5,0) -- (0,0);
    
    \draw (0,5) -- (10,0);
    \draw (0,2.5) -- (5,0);

    \filldraw[fill=gray!80] (0,0) -- (5,10) -- (5,7) -- (10/16,0) -- (0,0);
    \filldraw[fill=gray!80] (10/8,0) -- (2.5,0) -- (5,5) -- (5,6) -- (10/8,0) ;
      
    \filldraw[teal] (70/17,5) -- (70/68,0) -- (30/34,0) -- (270/68,5) -- (70/17,5);

    \draw[very thick] (170/52,350/104) -- (4.25,5);
    
    \draw [dashed, thick] (0,350/104) -- (4.25,5);

    \draw [very thick] (0,350/104) -- (0.8724691416999097,5);
    \draw [very thick] (2.782261032261033,0) -- (4.25,3.5);
    
    \begin{scope}
    \clip (0,5) -- (5,2.5) -- (5,5) -- (0,5);
    \draw [very thick, tomato] (2.782261032261033,0) -- (4.25,3.5);
    \end{scope}
    
    \node at (4.15,4.3) {$\Gamma_1$};
    
    \node at (1.3,3.55) {$F(\Gamma_1)$};
    
    \node at (4.35,3.2) {\color{tomato}$\Gamma_2$\color{black}};

    \draw (0,0) rectangle (5,5);
    \end{tikzpicture}};

    \node at (4,0) {
    \begin{tikzpicture}[scale=1.3]
    
    \node at (-0.15,2.98323) {\color{tomato}$y_2$\color{black}};
    
    \clip (0,0) rectangle (5,5);

    \filldraw[fill=gray!40] (0,150/16-5) -- (7,0) -- (6,0) -- (0,70/8-5) --(0,150/16-5);

    \fill [fill=white] (0,0) -- (2.5,5) -- (5,5) -- (2.5,0) -- (0,0);
    
    \draw (0,5) -- (10,0);
    \draw (0,2.5) -- (5,0);

    \filldraw[fill=gray!80] (0,0) -- (5,10) -- (5,7) -- (10/16,0) -- (0,0);
    \filldraw[fill=gray!80] (10/8,0) -- (2.5,0) -- (5,5) -- (5,6) -- (10/8,0) ;
    
    \begin{scope}
    \clip (0,5) -- (5,2.5) -- (5,5) -- (0,5);
    \draw [very thick, tomato] (2.782261032261033,0) -- (4.25,3.5);
    \end{scope}
    
    \draw[dashed, thick, tomato] (0,2.98323) -- (1.25,3.5);
    \draw[very thick, tomato] (0,2.98323) -- (1.25,6);
    
    \draw [very thick] (0.25249,3.59197) -- (0.45807,4.08851);
    
    \node at (0.55,2.95) {\color{tomato}$F(\Gamma_2)$\color{black}};
    
    \node at (0.2,4) {$\Gamma_3$};
    
    \node at (0.65,3.6) {$A_{3,2}^1$};
    
    \node at (4.35,3.2) {\color{tomato}$\Gamma_2$\color{black}};

    \draw (0,0) rectangle (5,5);
    
    \end{tikzpicture}};

    \node at (-1.02,3.45) {$x_1$};
    
      \end{tikzpicture}

    \caption{Left: The upper part $\Gamma_1 \subset A_1$ of a $v$-segment in $\sigma_3 \setminus (H(\sigma_2) \cap \sigma_3)$ and its images $F(\Gamma_1)$ (dashed), $H(\Gamma_1) \cap S_1$. Right: A right part $\Gamma_2 \subset A_1$ of $H(\Gamma_1)$ and its images $F(\Gamma_2)$, $H(\Gamma_2) \cap S_1$. This image necessarily contains a line segment $\Gamma_3$ traversing $A_{3,2}^1$.}
    \label{fig:HsigmaMR}
    \end{figure}

\subsection{Invariant cones}
\label{sec:invariantCones}
We now derive specific unstable and stable cone fields for the return map $H_\sigma$, wide enough to ensure invariance \textbf{(H1.1)} yet fine enough to produce tight bounds on expansion factors, vital for verifying \textbf{(H5)}. Define the cones $\mathcal{C}_1,\dots,\mathcal{C}_4$ by
\begin{enumerate}[label={($\mathcal{C}_{\arabic*}$)}]
    \item $3|v_1| \geq |v_2| \geq 7 |v_1|/3$, $v_1v_2>0$,
    \item $5|v_1|/3 \geq |v_2| \geq \varphi |v_1|$, $v_1v_2<0$,
    \item $5|v_1|/3 \geq |v_2| \geq \varphi |v_1|$, $v_1v_2>0$,
    \item $3|v_1| \geq |v_2| \geq 7 |v_1|/3$, $v_1v_2<0$,
\end{enumerate}
and the following stable cones
\begin{enumerate}[label={($\mathcal{C}_{\arabic*}^s$)}]
    \item $|v_1| \geq |v_2|$,
    \item $9/10 \geq v_2/v_1 \geq -8/10$,
    \item $8/10 \geq v_2/v_1 \geq -9/10$,
    \item $|v_1| \geq |v_2|$.
\end{enumerate}

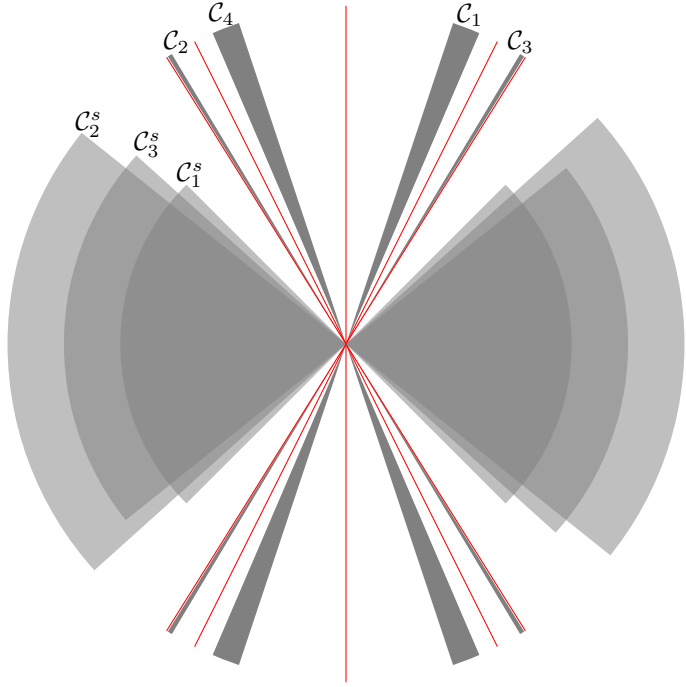
\begin{figure}
    \centering
    \begin{tikzpicture}[scale=1.5]
    
    \begin{scope}
    \clip (0,0) circle (3cm);
     \fill[gray] (0,0) -- (10,210/13) -- (10,50/3) -- (0,0);
     \fill[gray] (0,0) -- (-10,-210/13) -- (-10,-50/3) -- (0,0);
     \fill[gray] (0,0) -- (-10,210/13) -- (-10,50/3) -- (0,0);
     \fill[gray] (0,0) -- (10,-210/13) -- (10,-50/3) -- (0,0);
     
     \fill[gray] (0,0) -- (10,70/3) -- (10,30) -- (0,0);
     \fill[gray] (0,0) -- (10,-70/3) -- (10,-30) -- (0,0);
     \fill[gray] (0,0) -- (-10,70/3) -- (-10,30) -- (0,0);
     \fill[gray] (0,0) -- (-10,-70/3) -- (-10,-30) -- (0,0);
     
     \fill[gray,opacity=0.5] (0,0) -- (-10,8) -- (-10,-9) -- (0,0);
     \fill[gray,opacity=0.5] (0,0) -- (10,-8) -- (10,9) -- (0,0);
     
     \draw[red] (-10,-80/5) -- (10,80/5);
     \draw[red] (10,-80/5) -- (-10,80/5);
     
     \draw[red] (-10,-20) -- (10,20);
     \draw[red] (10,-20) -- (-10,20);
     
     \draw[red] (0,-10) -- (0,10);

    \end{scope}
   
    \begin{scope}
    \clip (0,0) circle (2.5cm);
    \fill[gray,opacity=0.5] (0,0) -- (-10,9) -- (-10,-8) -- (0,0);
     \fill[gray,opacity=0.5] (0,0) -- (10,-9) -- (10,8) -- (0,0);
    
    \end{scope}

   \begin{scope}
    \clip (0,0) circle (2cm);
    \fill[gray,opacity=0.5] (0,0) -- (-10,10) -- (-10,-10) -- (0,0);
     \fill[gray,opacity=0.5] (0,0) -- (10,-10) -- (10,10) -- (0,0);
    \end{scope}

  \node at (-1.38,1.4+0.1) {$\mathcal{C}_1^s$};
  \node at (-1.76,1.67+0.1) {$\mathcal{C}_3^s$};
  \node at (-2.28,1.85+0.1) {$\mathcal{C}_2^s$};
  
  \node at (-1.5,2.68) {$\mathcal{C}_2$};
    \node at (1.55,2.68) {$\mathcal{C}_3$};
    
    \node at (-1.1,2.93) {$\mathcal{C}_4$};
    \node at (1.1,2.92) {$\mathcal{C}_1$};

    \end{tikzpicture}
    \caption{Unstable and stable cone fields $\mathcal{C}_j$, $\mathcal{C}_j^s$ over the subsets $\sigma_j \subset \sigma$ for the return map $H_\sigma$. Also shown in red are the gradients of the line segments which make up the boundary $\partial \sigma$ which lie outside of all cone fields.}
    \label{fig:HsigmaCones}
\end{figure}
In the notation of section \ref{sec:outline}, for general $z \in \sigma$ we take $C_z^u = \mathcal{C}_j$ and $C_z^s = \mathcal{C}_j^s$ for $z \in \sigma_j$. These cone fields are plotted in Figure \ref{fig:HsigmaCones}.

\begin{table}[h]
    \centering
    \begin{tabular}{c||c|c|c|c}
         $DH_\sigma$ & $\sigma_1$ & $\sigma_2$ & $\sigma_3$ & $\sigma_4$  \\
         \hline \hline
         $\sigma_1$  & $M_1$ & - & $M_3M_2^k$ & \makecell{$M_4$ \\ $M_4M_2^k$ \\ $M_4M_3^k$} \\
         \hline
         $\sigma_2$  & - & - & \makecell{$M_3$ \\ $M_3M_2$ \\ $M_3M_2^2$} & \makecell{$M_4$ \\ $M_4M_2$ \\ $M_4M_2^2$} \\
         \hline
         $\sigma_3$  & \makecell{$M_1$ \\ $M_1M_3$ \\ $M_1M_3^2$ }& \makecell{$M_2$ \\ $M_2M_3$ \\  $M_2M_3^2$} & - & - \\
         \hline
         $\sigma_4$  & \makecell{$M_1$ \\$M_1M_2^k$ \\ $M_1M_3^k$} & $M_2M_3^k$ & - & $M_4$ \\
    \end{tabular}
    \caption{Possible values of the Jacobian $DH_\sigma$ at $z$ if $z \in \sigma_i$ (rows) and $z' = H_\sigma(z) \in \sigma_j$ (columns). Exponent $k$ takes values in $\mathbb{N}$, dashes are shown if no transition is possible, e.g $H_\sigma(\sigma_1) \cap \sigma_2 = \varnothing$.}
    \label{tab:sigmaTransitions}
\end{table}

\begin{lemma}
\label{lemma:finalReducedCones}
The above cones satisfy $DH_\sigma \, C_z^u \subset C_{z'}^u$ and $DH_\sigma \, C_z^s \supset C_{z'}^s$ for all $z\in \sigma$ where $DH_\sigma$ exists, $z'=H_\sigma(z)$. 
\end{lemma}

\begin{proof}
We begin with the unstable cones. Table \ref{tab:sigmaTransitions} shows the possible values of $DH_\sigma$ at $z$ if $z \in \sigma_i$ and $z' \in \sigma_j$. The calculations for $z' \in \sigma_1,\sigma_4$ are similar to those made in the proof of Lemma \ref{lemma:mappedOTMalignment}, noting that each $\mathcal{C}_j$ is contained within $\mathcal{C}$ and, for example, $M_1M_j^k \mathcal{C} \subset \mathcal{C}_1$ for $j=2,3$, and $k \geq 0$. For $z' \in \sigma_2$ we verify that $M_2M_3^k (-1,3)^T =(-1)^k (-24k+5,-40k-7)^T \in \mathcal{C}_2$ and $M_2M_3^k (-3,7)^T =(-1)^k (60k+11,-100k-15)^T \in \mathcal{C}_2$ for all $k \geq 1$ so that $DH_\sigma \, C_z^u \subset C_{z'}^u$ for $z \in \sigma_4$. For $z \in \sigma_3$ we have $M_2 (3,5)^T = (13,-21)^T \in \mathcal{C}_2$ and $M_2 (13,21)^T = (55,-89)^T \in \mathcal{C}_2$, ensuring invariance in this particular case also, despite $M_2$ being non-hyperbolic. Entirely symmetric calculations can be made for $z' \in \sigma_3$, verifying the result for all unstable cones.

For the stable cones, we remark that taking $C_z^s = \mathcal{C}_1^s$ for all $z \in \sigma $ would satisfy $DH_\sigma \, C_z^s \supset C_{z'}^s$ but since $M_2^{-1} (1,-1)^T = (-1,1)^T$ we would be unable to derive sufficient uniform bounds on expansion factors \textbf{(H1.2)}. The matrix $M_3^{-1}$ exhibits a similar problem so we must slim down the cones $C_z^s$ when $DH_\sigma^{-1} \in \{ M_2^{-1}, M_3^{-1} \}$ which, observing Table \ref{tab:sigmaTransitions}, is for $z \in \sigma_2,\sigma_3$. To remedy this, for such $z$ we slim down the cones $C_z^s$ to $\mathcal{C}_2^s,\mathcal{C}_3^s$ above. As these cones lie in the wider invariant cone $|v_1| \geq |v_2|$, the lemma follows from checking that $DH_\sigma \, C_z^s \supset C_{z'}^s$ for $z \in \sigma_2,\sigma_3$. This can be verified via direct calculations.
\end{proof}

\subsection{Structure of the singularity set}
\label{sec:singSetStructure}
Using the notation of \textbf{(H2)} in section \ref{sec:outline}, let $\mathcal{S}_0 = \cup_j \partial \sigma_j$, the union of $\partial \sigma$ and the red dashed lines in Figure \ref{fig:HSigma}. The set $M = \Omega \setminus \mathcal{S}_0$ is clearly dense in $\Omega$ and $H_\sigma$ is a $C^2$ diffeomorphism from $M\setminus \mathcal{S}_1$ onto $M \setminus \mathcal{S}_{-1}$, being linear on each component. 

The set $\mathcal{S}_0 \cup \mathcal{S}_1$ is the countable union of bounded line segments with the endpoints of each segment terminating on another segment, giving \textbf{(H2.2)}. 

The gradients of the segments in $\mathcal{S}_0$ take values in $\{ \pm 8/5, \pm 2, \infty \}$ which avoid unstable and stable cones $C_z^u$, $C_z^s$ (see Figure \ref{fig:HsigmaCones}). The gradients of singularity curves in $\sigma_1$ and $\sigma_4$ are bounded between -1 and 1 (approaching these limits as we approach the accumulation points) so lie in $\mathcal{C}_1^s,\mathcal{C}_4^s$. The gradients of singularity curves in $\sigma_2$ and $\sigma_3$ are bounded between -11/14 and 11/14 so lie in $\mathcal{C}_2^s,\mathcal{C}_3^s$ since $11/14 < 8/10$. Similar calculations show that the gradients of segments in $\mathcal{S}_{-1}$ lie in unstable cones.

We conclude this section with showing \textbf{(H2.4)}. Condition (\ref{eq:DxbCondition}) can only fail when $\| DH_\sigma \|$ becomes unbounded, i.e. at points $z$ approaching the accumulation points. We consider the case with $z \in A_{4,2}^k$ near $(0,1/4)$, the other cases are similar. Recall Figure \ref{fig:sigma1a} and the lines $\mathcal{L}_k$ from ($\dagger$). We note that $d(z,\mathcal{S}_1)$ is bounded above by the length of the segment joining $z=(x,y)$ to $(x,y_k(x))$ on $\mathcal{L}_k$, which in turn is bounded above by the height of the segment joining $(x,y_k(x))$ to $(x,y_{k-1}(x))$ on $\mathcal{L}_{k-1}$. This height is
\[ y_{k-1}(x) - y_k(x) \leq y_{k-1}\left( \frac{1}{4k-2}\right) - y_k\left( \frac{1}{4k-2} \right) = \frac{1}{2(2k-1)^2} \leq c_1 /k^2  \]
for some constant $c_1>0$. The operator norm of $DH_\sigma$ over $A_{4,2}^k$ satisfies $\| DH_\sigma \| \leq c_2 k$ for some $c_2>0$ so that (\ref{eq:DxbCondition}) holds for some $c>0$ whenever we choose $b>1/2$.

\subsection{One-step expansion}
\label{sec:twoStepExpansionOTM}
We will verify (\ref{eq:oneStep}) for the map $f=H_\sigma^2$, $q=1/2$. We begin with a basic statement on expansion over unstable curves.

\begin{lemma}
\label{lemma:piecewiseApprox}
Let $M$ be the constant Jacobian of $f$ over $W_i$, $V_i =f(W_i)$. Then
\[ \lambda^{-} := \inf_{v\in \mathcal{C}} \frac{\| M v \|}{\|v \|} \leq \frac{|V_i|}{|W_i|} \leq  \sup_{v\in \mathcal{C}} \frac{\| M v \|}{\|v \|} =: \lambda^+. \]
\end{lemma}

\begin{proof}
Given any $\varepsilon>0$, consider a piecewise linear approximation $\hat{W_i}$ to $W_i$ such that
\[ \left| \frac{|V_i|}{|W_i|} - \frac{|\hat{V_i}|}{|\hat{W_i}|}  \right| < \varepsilon \]
where $\hat{V_i} = f(\hat{W_i})$ gives a piecewise approximation for $V_i$. Each of the piecewise components will be line segments aligned with vectors in $\mathcal{C}$ so that their expansion factors will be bounded by $\lambda^\pm$, giving the result.
\end{proof}

We also derive basic inequalities on the length of a given $W_i$.

\begin{lemma}
\label{lemma:componentLength}
Let $\mathcal{L}_0,\mathcal{L}_1$ be the singularity curves on which $W_i$ terminates, write these intersections as $(x_0,y_0)$ and $(x_1,y_1)$. Then
\[ \sqrt{(x_1-x_0)^2 + (y_1-y_0)^2} \leq |W_i| \leq |y_1-y_0|\sqrt{1+\frac{1}{g^2}}\]
where $g = \inf|v_2/v_1|$ over $(v_1,v_2)^T \in \mathcal{C}$. 
\end{lemma}

\begin{proof}
Noting that the lower bound is trivial, we focus on the upper bound. Since $g>0$ for all unstable cones $\mathcal{C}$, the projection of $W_i$ to the $y$-axis is injective. Without loss of generality suppose $y_1>y_0$, then we can parameterise $W_i$ as a curve $(x(y),y)$ for $y_0\leq y \leq y_1$. Now
\begin{equation*}
\begin{split}
    |W_i| &= \int_{y_0}^{y_1} \sqrt{ \left(\frac{\mathrm{d}x}{\mathrm{d}y} \right)^2 + \left(\frac{\mathrm{d}y}{\mathrm{d}y}\right)^2 } \d y \\
    & \leq (y_1-y_0) \sup_{y_0\leq y \leq y_1} \sqrt{1+ \left(\frac{\mathrm{d}x}{\mathrm{d}y}\right)^2}\\
    & \leq (y_1-y_0)  \sqrt{1+ \frac{1}{g^2}}
    \end{split}
\end{equation*}
as tangent vectors $(x'(y),1)^T$ to $W_i$ lie in $\mathcal{C}$. 
\end{proof}

Let $P_1 = \{ (0,1/4), (1/2,1/4), (1/2,3/4), (1,3/4) \}$ denote the accumulation points similar to that of $\sigma_{1a}$, $P_2 = \{ (1/4,1/2), (1/4,1), (3/4,0), (3/4,1/2) \}$ the accumulation points similar to that of $\sigma_{1b}$. Let $\varepsilon$ be small. Given a set $P$, let $B_\varepsilon(P)$ denote the union of the balls $B_\varepsilon(p) \cap \sigma$, centred at $p\in P$ of radius $\varepsilon$. The following describes the images of balls about $P_1 \cup P_2$ under $H_\sigma$.

\begin{lemma}
\label{lemma:accumulationMapping}
Given small $\varepsilon>0$, there exists some $\varepsilon'>0$ such that $H_\sigma(B_\varepsilon(P_1 \cup P_2))$ covers $B_{\varepsilon'}(P_1 \cup P_2)$.
\end{lemma}

\begin{figure}
    \centering
        \subfigure[][]{
    \label{fig:twoStepA}
    \begin{tikzpicture}[scale=0.85*7]
    
    \clip (1.95,4.5) rectangle (3.2,6);
    \foreach \k in {1,...,60}{
    \draw[very thin] ({ 5+ (2*\k+2)*(5-7.5)/(2*\k+1)  }, { 5 })--({ 5+ (2*\k+2)*(7.5-7.5)/(2*\k+1) }, { 7.5 });
    }

    \fill[black] (2.47,5) -- (2.53,5.06) -- (2.53,5) -- (2.47,5) ;
    \draw (0,5) -- (5,5);
    \filldraw[fill=white](0,0) -- (5,10) -- (5,0) -- (0,0);
    \draw [blue, thick] (1713/700,4.9) --(2.49,5)-- (2.8,5.9);
    \draw [red, thick] (6639/2515,5.4348) -- (739/275,613/110);
    
    \filldraw[red,fill=red,opacity=0.2] (2,5) -- (30/14,5) -- (50/18,50/9) -- (40/14,40/7) -- (2,5);
    \node at (2.13, 5.05) {$A_{4,3}^k$};
    \node at (6639/2515-0.05,5.4348+0.05) {$W_k$};
    \node at (2.4,4.95) {$W_\star$};
    \node at (2.75,5.9) {\color{blue}$W$\color{black}};

    \node[scale=0.5] at (2,5) {$\bullet$};
    \node at (2,4.95) {$r_2$};
    
    \node[scale=0.5] at (30/14,5) {$\bullet$};
     \node at (30/14,4.95) {$r_3$};
    \node[scale=0.5] at  (50/18,50/9) {$\bullet$};
     \node at (50/18+0.05,50/9) {$r_4$};
     \node[scale=0.5] at (40/14,40/7) {$\bullet$};
      \node at (40/14+0.05,40/7) {$r_1$};
    \node at (30/14,4.7) {$A_1 \cap \sigma$};

    \end{tikzpicture}
    }
    \subfigure[][]{
    \label{fig:twoStepB}
    \begin{tikzpicture}
    
    \node at (0,0) {
    \begin{tikzpicture}[scale=0.85*12]
    \clip (4.96,6.6) rectangle (5.35,7.7);
     \foreach \k in {1,...,80}{
        \draw({ 2.5 + (2*\k+1)*(5-5)/(2*\k)  }, { 5 })--({ 2.5 + (2*\k+1)*(10-5)/(2*\k)  }, { 10 });
    }
    \fill[fill=white] (5,10) -- (5,7.5) -- (10,10) -- (5,10);
    
    \node at (5.20,7.36) {$A_{1,3}^l$};
    
    \draw (5,7.5) -- (10,10);
   
    \filldraw[fill=white] (5,10) rectangle (4.8,5);
    
    \filldraw[red,fill=red,opacity=0.3] ({5+10/(4*8+2)},{7.5 - 50/(8*8+4) }) -- ({5+10/(4*8-2)},{7.5 - 50/(8*8-4) }) -- (5, {7.5 + 10/(8*8+4) } ) -- (5, {7.5 + 10/(8*8+12) } ) -- ({5+10/(4*8+2)},{7.5 - 50/(8*8+4) });
    
    
    \node[scale=0.5] at ({5+10/(4*8+2)},{7.5 - 50/(8*8+4) }) {$\bullet$};
     \node at ({5+10/(4*8+2) - 0.02},{7.5 - 50/(8*8+4) }) {$r_1'$};
    
    \node[scale=0.5] at ({5+10/(4*8-2)},{7.5 - 50/(8*8-4) }) {$\bullet$};
    \node at ({5+10/(4*8-2)-0.02},{7.5 - 50/(8*8-4) }) {$r_2'$};
    \node[scale=0.5] at  (5, {7.5 + 10/(8*8+4) } ) {$\bullet$};
    \node at (5+0.02, {7.5 + 10/(8*8+4) } ) {$r_3'$};
     \node[scale=0.5] at (5, {7.5 + 10/(8*8+12) } ) {$\bullet$};
    \node at (5-0.02, {7.5 + 10/(8*8+12) } ) {$r_4'$};
    
    \fill[fill=black] (5,7.48) -- (5,7.5) -- (5+0.03,7.5+0.015) -- (5,7.48);
    
    \end{tikzpicture}
    };

    \draw[dashed] (-0.4,1.4) circle (0.6cm);
    \draw[dashed] ({-0.4+0.6*cos(315)},{1.4+0.6*sin(315)}) -- ({2.6+2*cos(160)},{0.1+2*sin(160)});
    
    \filldraw[dashed, fill=white] (2.6,0.1) circle (2cm);
    
    \node[scale=1.5] at (2.6+1.05,-1.2) {\color{red}$U_{k}$\color{black}};
    
    \begin{scope}
    \clip (2.6,0.1) circle (2cm);
    \draw[thick]  ({2.6+2*cos(170)},{0.1+2*sin(170)}) -- ({2.6+2*cos(80)},{0.1+2*sin(80)});
    \draw[thick]  ({2.6+2*cos(250)},{0.1+2*sin(250)}) -- ({2.6+2*cos(-5)},{0.1+2*sin(-5)});
     
    \filldraw[red,fill=red,opacity=0.2] (1,2) -- (2.5,2.6) -- (4.5,-2.4) -- (3,-3) -- (1,2);
     
    \draw[red,very thick] ({2.6+2*cos(110)},{0.1+2*sin(110)}) -- ({2.6+2*cos(295)},{0.1+2*sin(295)});
    \draw[blue,very thick] ({2.6+1.48*cos(109)},{0.1+1.48*sin(109)}) -- ({2.6+1.24*cos(296.5)},{0.1+1.24*sin(296.5)});
    \node[scale=1.5] at (2.6+0.6,0.1+0.2) {\color{blue}$U_{k,l}$\color{black}};
    \end{scope}

    \end{tikzpicture}
    }

    \caption{Part (a) shows an unstable curve $W$ passing near to the accumulation point $(1/4,1/2)$, split into $W_\star$ below $y=1/2$ and the collection $W_k \subset A_{4,3}^k$. Part (b) shows the image $U_k = H_\sigma(W_k) \subset H_\sigma \left(A_{4,3}^k \right)$, which for odd $k$ lies near the accumulation point $(1/2,3/4)$ and contains subcurves $U_{k,l} \subset A_{1,3}^l$.}
    \label{fig:oneStep}
\end{figure}
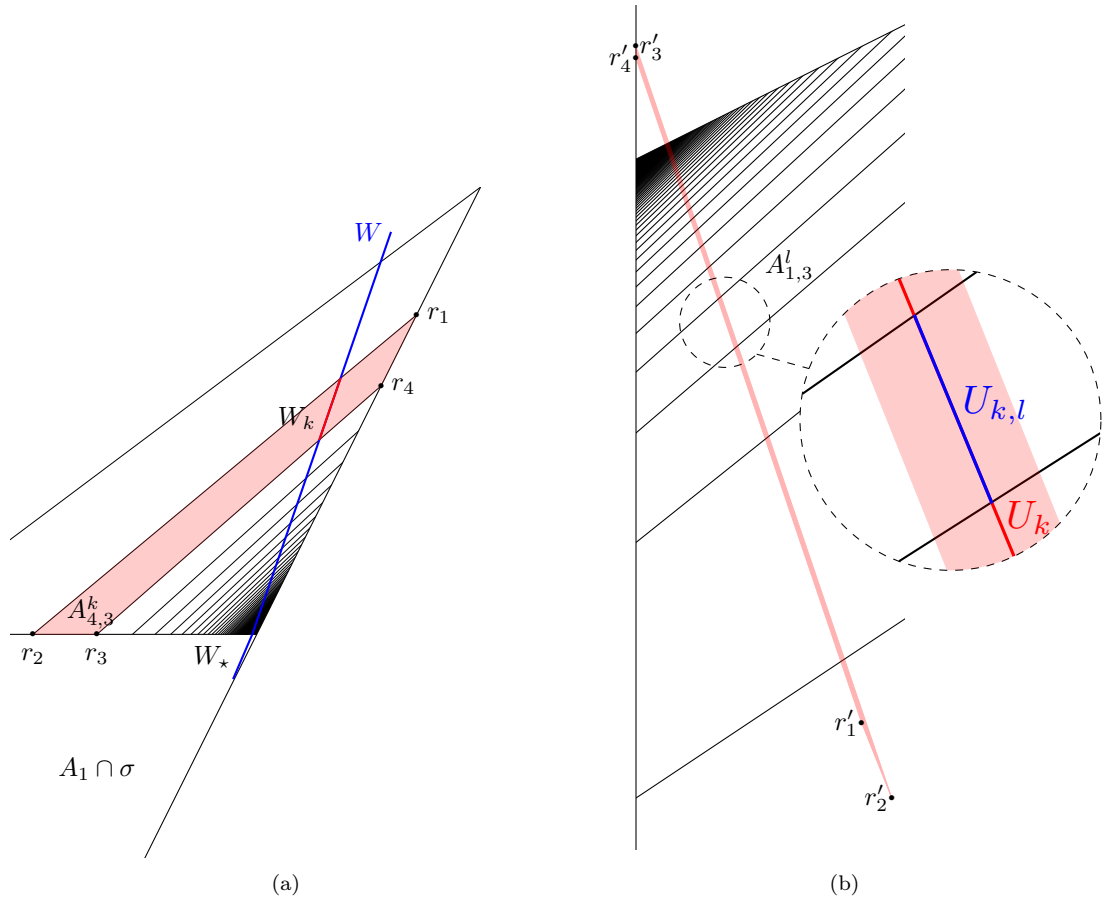

\begin{proof}

We describe the covering of $B_{\varepsilon'}((1/2,3/4))$, analysis for the other points in $P_1\cup P_2$ is analogous. For any $\varepsilon>0$, $B_\varepsilon(P_1 \cup P_2)$ contains the sets $A_{4,3}^k$ for all $k \geq k_0$ where $k_0 \in \mathbb{N}$ depends on $\varepsilon$. Each $A_{4,3}^k$ consists of two quadrilaterals, one in the ball around $(1/4,1)$ and the other in the ball around $(1/4,1/2)$. Figure \ref{fig:twoStepA} shows this latter quadrilateral, with corners on the points
\[  r_1 = \left( \frac{k+1}{4k+2}, \frac{k+1}{2k+1}  \right), \quad r_2 = \left(\frac{k-1}{4k-2}, \frac{1}{2}  \right), \quad r_3 = \left(\frac{k}{4k+2}, \frac{1}{2}  \right), \quad r_4 = \left(\frac{k+2}{4k+6}, \frac{k+2}{2k+3}  \right).  \]
Since $DH_\sigma$ is constant on $A_{4,3}^k$, given by the integer valued matrix $M_4M_3^k = (-1)^k \big(\begin{smallmatrix}
  2k+1 & -2k-2\\
  -6k-2 & 6k+5
\end{smallmatrix}\big)$,
its image $H_\sigma\left(A_{4,3}^k\right)$ is given by the quadrilateral with corners given by $M_4M_3^k \, r_j^T \text{ mod }1$. For odd $k$ we can calculate these corners as
\[  r_1'(k) = \left( \frac{1}{2} + \frac{1}{4k+2} , \frac{3}{4} -\frac{5}{8k+4}  \right), \quad r_2'(k) = \left( \frac{1}{2} + \frac{1}{4k-2} , \frac{3}{4} - \frac{5}{8k-4}  \right),\] \[r_3'(k) = \left( \frac{1}{2} , \frac{3}{4} + \frac{1}{8k+4}  \right), \quad r_4'(k) = \left( \frac{1}{2} , \frac{3}{4} + \frac{1}{8k+12}  \right), \]
shown in Figure \ref{fig:twoStepB}. For even $k$ the corners of $A_{4,3}^k$ in the ball around $(1/4,1)$ map into the $r_j'$. Writing this quadrilateral as $Q(k)$, since $r_2'(k+1)$ = $r_1'(k)$ and $r_3'(k+1) = r_4'(k)$ we have that $\cup_{k \geq k_0} Q(k)$ is the polygon with corners $r_2'(k_0)$, $r_3'(k_0)$, and $\lim_{k\to \infty} r_1'(k) = \lim_{k \to \infty} r_4'(k) = (1/2,3/4)$. Noting that $r_3'(k_0) > 3/4$ and $r_2'(k_0)$ lies on the line $y-\frac{3}{4} = -\frac{5}{2}(x-\frac{1}{2})$, there exists $\varepsilon'$ such that $\cup_{k \geq k_0} Q(k)$ covers all points $(x,y) \in B_{\varepsilon'}((1/2,3/4))$ with $y \geq \frac{3}{4} - \frac{5}{2}(x-\frac{1}{2})$. The image $H_\sigma( B_\varepsilon((1,3/4)) \cap A_4 )$ fills the remaining portion of $B_{\varepsilon'}((1/2,3/4))$, since $H_\sigma = H$ on $A_4$, $H(1,3/4) = (1/2,3/4)$, $DH \,(-2,1)^T = (0,-1)^T$, and $DH \,(0,-1)^T = (2,-5)^T$.
\end{proof}

\begin{prop}
\label{prop:simple2step}
Condition (\ref{eq:oneStep}) holds for $H_\sigma^2$ when there exists $\varepsilon>0$ such that $W \cap B_\varepsilon(P_1 \cup P_2) = \varnothing$.
\end{prop}

\begin{proof}

We claim that an unstable curve $W$ of vanishing length, bounded away from the accumulation points, is split into at most 9 components $W_i$ by the singularity set for $H_\sigma^2$. The upper bound follows from analysis of the original singularity set for $H_\sigma$. Let $P_F$ denote the set of fixed points under $H$, $P_F = \{ (0,1/2), (1/2,0), (1/2,1/2), (1,1) \}$. Observing Figure \ref{fig:sigmaSingSet}, if $W \cap B_\varepsilon(P_F) \neq \varnothing$ then $W$ is split by $\mathcal{S}$ into at most 5 components $W_j$, and if $W \cap B_\varepsilon(P_F) = \varnothing$ then the upper bound is 3. We consider these cases separately. 

Take, for example, $W \cap B_\varepsilon((0,1/2)) \neq \varnothing$. Observing Figure \ref{fig:sigma1b}, four of the components $W_j$ map into $A_4'$ under $H_\sigma$, and their images lie in some sector $B_{\varepsilon'}((1,1/2)) \cap A_4'$. We can take $\varepsilon$ small enough that this sector lies entirely in $A_4$, so that no further splitting occurs during the next iterate of $H_\sigma$. The other component $W \cap A_1$ maps into some sector $B_{\varepsilon'}((0,1/2)) \cap A_1'$ and is split into at most 5 components, giving at most $N=9$ components in total. The other cases $W \cap B_\varepsilon(p) \neq \varnothing$, $p\in P_F$, are analogous. Now suppose $W \cap B_\varepsilon(P_F) = \varnothing$. $\mathcal{S}$ splits $W$ into at most 3 components $W_j$ and, by Lemma \ref{lemma:accumulationMapping} and the above, each $H_\sigma(W_j)$ is bounded away from the accumulation points $P_1 \cup P_2$ and the fixed points $P_F$. Hence each $H_\sigma(W_j)$ is split into at most 3 components during the next iterate of $H_\sigma$, again giving at most $N=9$ components in total.

The weakest expansion of $DH_\sigma^2$ over cones $\mathcal{C}_j$ on $\sigma_j$ using the euclidean norm is that of $M_1M_4 = \big(\begin{smallmatrix}
  -3 & 8\\
  -8 & 21
\end{smallmatrix}\big)$ on $\sigma_1$ (or equivalently $M_4M_1$ on $\sigma_4$), and is given by
\[ c = \frac{\| M_1M_4 (3,7)^T \|}{\| (3,7)^T \|} \sqrt{ \frac{(-9+56)^2+(-24+147)^2}{3^2+7^2}} = \sqrt{\frac{8669}{29}} \approx 17.29\]
so that, by Lemma \ref{lemma:piecewiseApprox}, $|V_i| \geq c \,|W_i|$ for each component $W_i$. 
Now for $q=1/2$ we have
\begin{equation*}
  \begin{split}
      \sum_i \left( \frac{|W|}{|V_i|}\right)^q \frac{|W_i|}{|W|} & = \sum_i \sqrt{\frac{|W|}{|V_i|}} \frac{|W_i|}{|W|} \\
      &= \sum_i \sqrt{ \frac{|W_i|}{|V_i|} } \sqrt{ \frac{|W_i|}{|W|}}\\
      &\leq \frac{1}{\sqrt{c}}  \sum_{i=1}^N \sqrt{ \frac{|W_i|}{|W|}}.
  \end{split} 
\end{equation*}
Letting $x_i = |W_i|/|W|$ and taking vectors $u = \left(\sqrt{x_1},\dots,\sqrt{x_N} \right)$, $v = (1,\dots,1)$ we have that $\sum_{i=1}^N x_i = 1$ and so
$\left( \sum_{i=1}^N \sqrt{x_i} \right)^2 = \left( u \cdot v \right)^2 \leq (u \cdot u) \, (v \cdot v) = \left( \sum_{i=1}^N x_i \right) N = N$
by the Cauchy-Schwarz inequality. Hence
\[\sum_i \left( \frac{|W|}{|V_i|}\right)^q \frac{|W_i|}{|W|} \leq \frac{\sqrt{N}}{\sqrt{c}} <1\]
since $N\leq9<c$. 
\end{proof}

\begin{prop}
\label{prop:upperTwoStep}
Condition (\ref{eq:oneStep}) holds for $H_\sigma^2$ when $W \cap B_\varepsilon(P_2) \neq \varnothing$ for all $\varepsilon>0$.
\end{prop}

\begin{proof}
We begin with the case $W \cap B_\varepsilon((1/4,1/2)) \neq \varnothing$ and let $\varepsilon \to 0$. We may choose $\delta$ sufficiently small so that $W$ intersects $A_1 \cap \sigma$ and some collection of sets $A_{4,3}^k$, $k_0 \leq k \leq k_1$, where $k_1 \to \infty$ as $\varepsilon \to 0$, $k_0 \to \infty$ as $\delta \to 0$. Therefore, $\mathcal{S}$ splits $W$ into a lower component $W_\star \subset A_1 \cap \sigma$ and upper components $W_k \subset A_{4,3}^k$, illustrated in Figure \ref{fig:twoStepA}. We study how the images of these components under $H_\sigma$ are split up by $\mathcal{S}$.

Recall the corners $r_j(k)$ which define $A_{4,3}^k$ near $(1/4,1/2)$. The curve $W_k$ has endpoints on $r_1r_2$ and $r_3r_4$, and all tangent vectors to $W$ lie in $\mathcal{C}_1$. For odd $k$ the image $U_k = H_\sigma(W_k)$ is a curve joining $r_1'r_2'$ to $r_3'r_4'$, with tangent vectors aligned in $M_4M_3^k \, \mathcal{C}_1$. This curve is split by $\mathcal{S}$ into an upper portion $U_{k,\star} \subset A_4 \cap \sigma$, and a collection $U_{k,l} \subset A_{1,3}^l$ for some consecutive range $l_0 \leq l \leq l_1$ which depends on $k$. Each $A_{1,3}^l$ is bounded by the lines
\begin{equation}
    \label{eq:linesLl}
    \mathcal{L}_l: y- \frac{1}{2} = \frac{2l}{2l+1} (x-1/4)
\end{equation}
and $\mathcal{L}_{l-1}$, hence a lower bound on $l_0(k)$ is given by the largest $l$ such that $r_2'(k)$ lies on or above $\mathcal{L}_{l-1}$. One can verify that $r_2'(k)$ lies on $\mathcal{L}_{l-1}$ when $k = 7l-4$ and approaches $(1/2,3/4)$ monotonically in $x$ and $y$ so that $l_0(k) \geq \lfloor \frac{k+4}{7} \rfloor$. To determine an upper bound on $l_1$, note that $r_4'r_1'$ lies on the line
\begin{equation}
    \label{eq:cellLinesk}
    y - \frac{3}{4} -\frac{1}{8k+12} = - \frac{6k+8}{2k+3} \left(x-\frac{1}{2} \right),
\end{equation}
meeting the $A_4$ boundary $\mathcal{L} : y=1/2 + x/2$ at the point
\begin{equation}
\label{eq:xkykonL}
    (x_k,y_k) =  \left( \frac{7k+10}{14k+19}, \frac{21k+29}{28k+38} \right).
\end{equation}
We similarly calculate that the line $\mathcal{L}_l$ meets $y=1/2 + x/2$ at the point
\[ (X_l,Y_l) = \left( \frac{l}{2l-1}, \frac{1-3l}{2-4l}  \right). \]
The intersection of $U_k$ with $y=1/2 + x/2$ must be some point $(x,1/2+x/2)$ with $x\geq x_k$ so that an upper bound on $l_1(k)$ is the smallest $l$ such that $x_k \geq X_l$, which reduces to $l \geq 7k+10$, hence $l_1(k) \leq \lceil 7k+10 \rceil = 7k+10$. For even $k$ the splitting behaviour is entirely analogous, with $H_\sigma(W_k)$ intersecting $\mathcal{S}$ in the neighbourhood of $(1,1/4)$.

For the lower component $W_\star$, the image $U_\star$ = $H_\sigma(W_\star) = H(W_\star)$ lies in a neighbourhood of $H(1/4,1/2) = (1/4,1)$ and is split by $\mathcal{S}$ into a collection $U_{\star,j} \subset A_{4,3}^j$, $j_0 \leq j \leq j_1$, where $j_1 \to \infty$ as $\varepsilon \to 0$, $j_0 \to \infty$ as $\delta \to 0$. Write $W_{\star,j} = H_\sigma^{-1}(U_{\star,j})$, $W_{k,\star} = H_\sigma^{-1}(U_{k,\star})$, $W_{k,l} = H_\sigma^{-1}(U_{k,l})$, then $W$ splits into components
\begin{equation}
    \label{eq:Wsplitting}
    W = \left( \bigcup_{j \geq j_0} W_{\star,j} \right) \cup \left( \bigcup_{k \geq k_0}  W_{k,\star} \right) \cup \left(  \bigcup_{k \geq k_0} \bigcup_{l=l_0}^{l_1} W_{k,l}  \right)
\end{equation}
on which $DH_\sigma^2$ is constant. Let $V_{i} = H_\sigma(U_i) = H_\sigma^2(W_i)$, then for $q=1/2$:
\begin{equation*}
\begin{split}
    \sum_i \left( \frac{|W|}{|V_i|}\right)^q \frac{|W_i|}{|W|} & = \sum_i \sqrt{ \frac{|W_i|}{|V_i|} } \sqrt{ \frac{|W_i|}{|W|}}\\
    & = \sum_{j \geq j_0}  \sqrt{ \frac{|W_{\star,j}|}{|V_{\star,j}|} } \sqrt{ \frac{|W_{\star,j}|}{|W|}} +  \sum_{k \geq k_0}         \sqrt{ \frac{|W_{k,\star}|}{|V_{k,\star}|} } \sqrt{ \frac{|W_{k,\star}|}{|W|}}   + \sum_{k \geq k_0} \sum_{l=l_0}^{l_1} \sqrt{ \frac{|W_{k,l}|}{|V_{k,l}|} } \sqrt{ \frac{|W_{k,l}|}{|W|}} \\
    & \leq \sum_{j \geq j_0}  \sqrt{ \frac{1}{\Lambda_{\star,j} }} \sqrt{ \frac{|W_{\star,j}|}{|W|}} +  \sum_{k \geq k_0}  \sqrt{\frac{1}{ \Lambda_{k,\star}}} \sqrt{ \frac{|W_{k,\star}|}{|W|}}   + \sum_{k \geq k_0} \sum_{l=l_0}^{l_1} \sqrt{\frac{1}{ \Lambda_{k,l}}} \sqrt{ \frac{|W_{k,l}|}{|W|}}
    \end{split}
\end{equation*}
by Lemma \ref{lemma:piecewiseApprox}, where $\Lambda_i$ is the minimum expansion factor of $DH_\sigma^2$ on $W_i$ over the cone $\mathcal{C}_1$. Define $W_\diamond = W \setminus W_\star$ and let $0\leq p \leq 1$ denote the proportion $|W_\diamond| = p\,|W|$, then
\begin{equation*}
\begin{split}
  \liminf_{\delta \to 0} \sup_{W: |W|< \delta}  \sum_i \left( \frac{|W|}{|V_i|}\right)^q \frac{|W_i|}{|W|} \leq \sup_{0 \leq p \leq 1} \Bigg( & \lim_{j_0 \to \infty}  \sum_{j \geq j_0}  \sqrt{ \frac{1}{\Lambda_{\star,j} }} \sqrt{ \frac{(1-p)|W_{\star,j}|}{|W_\star|}} \\
  & + \lim_{k_0 \to \infty} \sum_{k \geq k_0}  \sqrt{\frac{1}{ \Lambda_{k,\star}}} \sqrt{ \frac{p\,|W_{k,\star}|}{|W_\diamond|}} \\
  & + \lim_{k_0 \to \infty} \sum_{k \geq k_0} \sum_{l=l_0}^{l_1} \sqrt{\frac{1}{ \Lambda_{k,l}}} \sqrt{ \frac{p\,|W_{k,l}|}{|W_\diamond|}}\, \Bigg).
    \end{split}
\end{equation*}
We put upper bounds on each of these sums using lower bounds on the expansion factors $\Lambda_i$ and geometric bounds on the curves $U_i$ terminating on $\mathcal{S}$. We use asymptotic notation $f \sim g$ for functions $f,g$ if $f/g \to 1$, and write $f \lesssim g$ if there is some function $h$ such that $f \leq h \sim g$.

Starting with the first sum, $DH_\sigma^2$ is given by $M_4M_3^jM_1 = (-1)^j \left( \begin{smallmatrix} - 2 j - 3 & - 6 j - 8\\6 j + 8 & 18 j + 21  \end{smallmatrix} \right)$ on each component $W_{\star,j}$ with minimum expansion factors given by
\begin{equation}
    \label{eq:Lambdastarj}
    \begin{split}
         \Lambda_{\star,j} & = \inf_{7/3 \leq m \leq 3} \sqrt{ \frac{ \left( -2j-3 -6jm-8m \right)^2 + (6j+8 +18jm + 21)^2  }{1+m^2}}\\
         & \sim \inf_{7/3 \leq m \leq 3} \sqrt{ \frac{ \left( 2+6m \right)^2 + (6 +18m)^2}{1+m^2}} \,j = \frac{48 \sqrt{145}}{29} \,j := c_\star \,j.
    \end{split}
\end{equation}
Each curve $U_{\star,j}$ has tangent vectors in $M_1\mathcal{C}_1$ satisfying $41/17 \leq |v_2|/|v_1| \leq 17/7$, $v_1v_2 \geq 0$. For each $j>j_0$, $U_{\star,j}$ traverses $A_{4,3}^j$ so that (making a similar calculation to equation \ref{eq:sigma1bhk}) Lemma \ref{lemma:componentLength} gives
\begin{equation}
    \label{eq:Ustarj}
    \frac{ a_\star }{j^2} \lesssim |U_{\star,j}| \lesssim \frac{b_\star}{j^2}
\end{equation}
for $a_\star = \frac{13}{80} \sqrt{2}$, $b_\star = \frac{41}{192} \sqrt{1+ 17^2/41^2}$ (calculated in the appendix, section \ref{sec:upperCalculationsAppendix}). The upper bound also trivially holds for $j=j_0$. Let $\Lambda_1^+$, $\Lambda_1^-$ denote the maximum and minimum expansion factors of $M_1$ over $\mathcal{C}_1$, then
\[|W_\star| = \sum_{j\geq j_0} |W_{\star,j}| \geq \sum_{j\geq j_0} \frac{|U_{\star,j}|}{\Lambda_1^+} \gtrsim \frac{a_\star}{\Lambda_1^+} \sum_{j \geq j_0+1} \frac{1}{j^2} \geq \frac{a_\star}{\Lambda_1^+ (j_0+1)} \sim \frac{a_\star}{\Lambda_1^+ j_0}\]
where we have used the fact that
\[ \frac{1}{j^2} \geq \frac{1}{j(j+1)} = \frac{1}{j} - \frac{1}{j+1}\]
and considered the telescoping sum. Similarly $|W_{\star,j}| \lesssim b_\star / \left(\Lambda_1^- j^2 \right)$ so that
\[ \frac{|W_{\star,j}|}{|W_\star|} \lesssim \frac{b_\star \Lambda_1^+ j_0}{a_\star\Lambda_1^-} \, \frac{1}{j^2}.\]
Hence
\begin{equation}
\label{eq:Wstar}
    \begin{split}
        \sum_{j \geq j_0}  \sqrt{ \frac{1}{\Lambda_{\star,j} }} \sqrt{ \frac{(1-p)|W_{\star,j}|}{|W_\star|}} & \lesssim  \sum_{j \geq j_0} \sqrt{\frac{1}{c_\star j}} \sqrt{ \frac{(1-p)b_\star \Lambda_1^+ j_0}{a_\star\Lambda_1^-} \, \frac{1}{j^2} }\\
        & = \sqrt \frac{(1-p)b_\star \Lambda_1^+ j_0}{c_\star a_\star\Lambda_1^-} \sum_{j \geq j_0} j^{-3/2}\\
        & \leq \sqrt \frac{(1-p)b_\star \Lambda_1^+}{c_\star a_\star\Lambda_1^-} \, 2 \sqrt{\frac{j_0}{j_0-1}} \\
        & \to 2\sqrt \frac{(1-p)b_\star \Lambda_1^+}{c_\star a_\star\Lambda_1^-}
    \end{split}
\end{equation}
as $j_0 \to \infty$, where we have used $\sum_{j \geq j_0} j^{-3/2} \leq \int_{j_0-1}^\infty x^{-3/2} \d x$.

Moving onto the next summation, $\Lambda_{k,\star}$ is determined by $M_4^2M_3^k = (-1)^k \left( \begin{smallmatrix}  14 k + 5 & - 14 k - 12\\- 34 k - 12 & 34 k + 29 \end{smallmatrix} \right)$ and satisfies
\[ \Lambda_{k,\star} \sim \inf_{7/3 \leq m \leq 3} \sqrt{\frac{(14-14m)^2+(34-34m)^2}{1+m^2}} \, k = \frac{104}{\sqrt{29}}\, k =: c_\diamond \, k.\]
For odd $k$ the curve $U_{k,\star}$ has endpoints $(1/2,y_0)$ on $r_3'r_4'$ and $(2y_1-1,y_1)$ on $\mathcal{L}$, where $y_1$ is bounded by the intersections of $r_4'r_1'$ and $r_3'r_2'$ with $\mathcal{L}$ (see Figure \ref{fig:twoStepB}). An upper bound on $|y_1-y_0|$ is given by taking $(1/2,y_0) = r_3'$ and $(2y_1-1,y_1) = r_4'r_1' \cap \mathcal{L}$. Noting (\ref{eq:xkykonL}), this gives
\[ |y_1-y_0| \leq \frac{3}{4} + \frac{1}{8k+4} - \frac{21k+29}{28k+38} = \frac{6k+9}{56k^2+104k+38} \sim \frac{6}{56}k^{-1}. \]
$U_{k,\star}$ has tangent vectors $(v_1,v_2)^T$ in $M_{4,3}^k \, \mathcal{C}_1 \subset \mathcal{C}_4$, so that $|v_2/v_1| \geq 7/3$. By Lemma \ref{lemma:componentLength} we then have $|U_{k,\star}| \lesssim b_\diamond /k$, where $b_\diamond = \frac{6}{56} \sqrt{1 + \frac{9}{49}} \approx 0.117$. The minimum expansion factor of $M_4M_3^k = (-1)^k \left( \begin{smallmatrix}
  2k+1 & -2k-2\\
  -6k-2 & 6k+5
\end{smallmatrix}\right)$ over $\mathcal{C}_1$ is given by
\[\inf_{7/3 \leq m \leq 3} \sqrt{\frac{(2-2m)^2+(6-6m)^2}{1+m^2}} \, k = \frac{8\sqrt{145}}{29}\, k =: \gamma \, k\]
which gives $|W_{k,\star}| \lesssim b_\diamond/\left(\gamma k^2 \right)$. The analysis for even $k$ is analogous and gives the same upper bound. We next require a lower bound on $|W_\diamond|$. For $k>k_0$, $W_k$ is a curve with tangent vectors in $\mathcal{C}_1$ which traverses $A_{4,3}^k$. Making the same calculation as (\ref{eq:sigma1bhk}), $|W_k|$ is bounded below by the shortest path across $A_{4,3}^k$, the line segment passing through $r_4$ with gradient $3$. That's
\begin{equation}
\label{eq:Wklength}
    |W_k| \geq \sqrt{ \left( \frac{1}{16k^2+28k+6} \right)^2 + \left( \frac{3}{16k^2+28k+6} \right)^2 } \sim \frac{\sqrt{10}}{16k^2}
\end{equation}
so that 
\begin{equation}
\label{eq:Wdiamond}
    |W_\diamond| \geq \sum_{k\geq k_0+1} |W_k| \gtrsim \frac{a}{k_0}
\end{equation}
with $a:= \sqrt{10}/16$. Hence
\[\sum_{k \geq k_0}  \sqrt{\frac{1}{ \Lambda_{k,\star}}} \sqrt{ \frac{p\,|W_{k,\star}|}{|W_\diamond|}} \lesssim  \sum_{k \geq k_0} \sqrt{\frac{1}{c_\diamond k}} \sqrt{ \frac{p b_\diamond k_0}{a \gamma} \, \frac{1}{k^2} } \to 2\sqrt{\frac{pb_\diamond}{c_\diamond a \gamma}}\]
as $k_0 \to \infty$, by a similar argument to (\ref{eq:Wstar}).

For the third summation, $\Lambda_{k,l}$ is determined by the matrix \[M_1M_3^lM_4M_3^k = (-1)^{k+l} \begin{pmatrix} -48 k l - 10 k - 18 l - 3 & 48 k l + 10 k + 42 l + 8 \\ -112 k l - 26 k - 42 l - 8 & 112 k l + 26 k + 98 l + 21 \end{pmatrix}  \]
and satisfies (for large $k,l$)
\[ \Lambda_{k,l} \sim \inf_{7/3 \leq m \leq 3} \sqrt{\frac{(48-48m)^2+(112-112m)^2}{1+m^2}} \, kl = 64 \,kl =: c \,kl.\]
We can show an upper bound $|U_{k,l}| \lesssim b/l^2$ where $b = \frac{3}{32}\sqrt{1 + \frac{9}{49}} \approx 0.102$ (see section \ref{sec:upperCalculationsAppendix}) so that $|W_{k,l}| \lesssim b/(\gamma k l^2)$. Now by (\ref{eq:Wdiamond}),
\begin{equation*}
    \begin{split}
        \sum_{l=l_0}^{l_1} \sqrt{\frac{1}{ \Lambda_{k,l}}} \sqrt{ \frac{p|W_{k,l}|}{|W_\diamond|}} & \lesssim \sum_{l=l_0}^{l_1} \sqrt{\frac{1}{c\,kl}} \sqrt{\frac{p b k_0}{a \gamma k l^2}} \\
        & \leq \sqrt{\frac{1}{c\,k}} \sqrt{\frac{p b k_0}{a \gamma k}} \sum_{l=\lfloor \frac{k+4}{7} \rfloor}^{7k+10} l^{-3/2} \\
        & \leq 2 \sqrt{\frac{1}{c\,k}} \sqrt{\frac{p b k_0}{a \gamma k}} \left( \frac{1}{\sqrt{\lfloor \frac{k+4}{7} \rfloor -1}} - \frac{1}{\sqrt{7k+10}}   \right) \\
        & \sim 2 \sqrt{\frac{1}{c\,k}} \sqrt{\frac{p b k_0}{a \gamma k}} \left( \sqrt{7}- \frac{1}{\sqrt{7}} \right) \frac{1}{\sqrt{k}}.
    \end{split}
\end{equation*}
Letting $h = (\sqrt{7} - 1/\sqrt{7})^2 = 36/7$, we have that
\begin{equation*}
    \begin{split}
        \sum_{k \geq k_0} \sum_{l=l_0}^{l_1} \sqrt{\frac{1}{ \Lambda_{k,l}}} \sqrt{ \frac{p|W_{k,l}|}{|W_\diamond|}} & \lesssim 2 \sqrt{\frac{p b h k_0}{ca \gamma}}\sum_{k \geq k_0} k^{-3/2} \\
        & \to 4 \sqrt{\frac{p b h}{ca \gamma}}
    \end{split}
\end{equation*}
as $k_0 \to \infty$. Hence for $q=1/2$
\begin{equation}
    \label{eq:justConstants}
 \liminf_{\delta \to 0} \sup_{W: |W|< \delta}  \sum_i \left( \frac{|W|}{|V_i|}\right)^q \frac{|W_i|}{|W|} \leq \sup_{0 \leq p \leq 1} \left( 2\sqrt \frac{(1-p)b_\star \Lambda_1^+}{c_\star a_\star\Lambda_1^-} +    2\sqrt{\frac{pb_\diamond}{c_\diamond a \gamma}} +  4 \sqrt{\frac{p b h}{ca \gamma}} \right).
\end{equation}
It is simple to show that for $s,t>0$ the function $f(p) = s\sqrt{1-p} + t\sqrt{p}$ always attains its maximum value at $p = t^2/(s^2+t^2)$. Hence letting
\[ s =  2\sqrt \frac{b_\star \Lambda_1^+}{c_\star a_\star\Lambda_1^-} \approx 0.450, \quad t = 2\sqrt{\frac{b_\diamond}{c_\diamond a \gamma}} +  4 \sqrt{\frac{ b h}{ca \gamma}} \approx 0.639 \]
gives
\[ \liminf_{\delta \to 0} \sup_{W: |W|< \delta}  \sum_i \left( \frac{|W|}{|V_i|}\right)^q \frac{|W_i|}{|W|} \leq s\sqrt{\frac{s^2}{s^2+t^2}} + t\sqrt{\frac{t^2}{s^2+t^2}} \approx 0.781 < 1  \]
as required. The analysis is analogous for $W$ near $(1/4,1)$ and extends to $W$ near $(3/4,0)$ and $(3/4,1/2)$ using the symmetry $T(x,y) = (1-x,y+1/2)$ which commutes with $H_\sigma$ (as seen in the proof of Lemma \ref{lemma:unstableGrowth}).
\end{proof}

Equivalent analysis verifies the two step expansion for curves near the other accumulation points $p \in P_1$. We provide the relevant calculations to the appendix, Proposition \ref{prop:lowerTwoStep}. We are now ready to apply Theorem \ref{thm:chernovZhang}.

\subsection{Decay of correlations}

\begin{thm}
\label{thm:HsigmaExp}
The return map $H_\sigma: \sigma \to \sigma$ enjoys exponential decay of correlations. In particular it admits a Young tower with base $\Delta_0$ satisfying the exponential tail bound
\begin{equation}
    \label{eq:delta0tailbound}
    \mu\left( \{ z \in \sigma \, | \, R(z,H_\sigma,\Delta_0) > n \} \right) \leq \mathrm{const} \, \theta^n
\end{equation}
for all $n \geq 1$ where $\theta<1$ is some constant. 
\end{thm}

\begin{proof}

We run through the conditions for applying Theorem \ref{thm:chernovZhang}. Invariance of the unstable and stable cone fields $C_z^u$, $C_z^s$ was the subject of section \ref{sec:invariantCones}, satisfying \textbf{(H1.1)}. Condition \textbf{(H1.2)} follows by taking $\Lambda = \sqrt{85/41}$, with this lower bound attained by considering the expansion of $M_2^{-1}$ over the cone boundary of $\mathcal{C}_2^s$ with gradient $-8/10$. Noting Remark \ref{remark:firstRemark}, we next show \textbf{(H1.3')}. The cone fields are continuous over the components $\sigma_j$ of $\Omega \setminus \mathcal{S}_0$, indeed they are constant. Noting that all of the stable cones lie within $\mathcal{C}_1^s$ and all of the unstable cones lie in $\mathcal{C}$, a positive angle between stable and unstable cone fields follows from $\mathcal{C}_1^s \cap \mathcal{C} = \varnothing$. \textbf{(H2)} was the subject of section \ref{sec:singSetStructure}. Unstable manifolds provide a class of $H_\sigma$ invariant unstable curves which satisfy the regularity conditions listed in \cite{chernov_statistical_2009}. Piecewise linearity of the map trivially implies their bounded curvature and bounds on distortion; absolute continuity follows from Lemma \ref{lemma:HsigmaKatokStrelcyn}. \textbf{(H4)} follows from Proposition \ref{prop:HsigmaBernoulli}, where we showed that $H_\sigma$ is Bernoulli with respect to the normalised Lebesgue measure on $\sigma$. Noting Remark \ref{remark:firstRemark}, \textbf{(H5)} follows for the map $H_\sigma^2$ by Propositions \ref{prop:simple2step}, \ref{prop:upperTwoStep}, \ref{prop:lowerTwoStep}.
\end{proof}

\section{Decay of correlations for the OTM}
\label{sec:polyMixingRate}
We now turn to the upper bound on correlations for $H$. We follow the approach outlined in \cite{chernov_improved_2008} to infer the polynomial mixing rate of $H$ from the exponential mixing rate of $H_\sigma$. For reference, their $\mathcal{M}$ and $M \subset \mathcal{M}$ are our $\tor$ and $\sigma \subset \tor$, their $\mathcal{F}: \mathcal{M} \to \mathcal{M}$ and $F:M \to M$ are our $H$ and $H_\sigma$ respectively. With $\Delta_0$ above, define
\[ A_n = \{ z \in \tor \, | \, R(z,H,\Delta_0 > n) \}. \]
We will show
\begin{prop}
\label{prop:AnPolynomial}
$\mu(A_n) = \mathcal{O}(n^{-1})$. 
\end{prop}
Theorem \ref{thm:mainTheoremMixingRate} then follows from the work of \cite{young_recurrence_1999}. Proving Proposition \ref{prop:AnPolynomial} involves treating separately a set of infrequently returning points, a method due to \cite{markarian_billiards_2004}. For each $z \in \tor$ and $n \geq 1$ define
\[ r(z;n,\sigma) = \# \{ 1 \leq i \leq n \, | \, H^i(z) \in \sigma  \},  \]
counting the number of times the orbit of $z$ hits $\sigma$ over $n$ iterates of $H$. Define 
\[B_{n,b} = \{ z \in \tor \, | \, r(z;n,\sigma) > b \ln n \} \]
where $b$ is a constant to be chosen shortly.
\begin{lemma}
\label{lemma:AncapBn}
$\mu(A_n \cap B_{n,b}) = \mathcal{O}(n^{-1})$.
\end{lemma}
\begin{proof}
This follows from (\ref{eq:delta0tailbound}), choosing $b$ large enough so that $n \, \theta^{b \ln n} < n^{-1}$. See \cite{chernov_improved_2008}, or \cite{springham_polynomial_2014} for a detailed proof.
\end{proof}
Proposition \ref{prop:AnPolynomial} then follows from similarly establishing
\begin{lemma}
\label{lemma:AnlessBn}
$\mu(A_n \setminus B_{n,b}) = \mathcal{O}(n^{-1})$.
\end{lemma}
Analysis of the set $A_n \setminus B_{n,b}$ is the focus of \cite{chernov_improved_2008}. It consists, for large $n$, of points which take many iterates to hit $\Delta_0$ and hit $\sigma$ infrequently during these iterates. \citeauthor{chernov_improved_2008} define $m$-cells
\[ M_m = \{ z \in \sigma \, | \, R(z;H,\sigma) = m+1 \}  \]
for $m \geq 0$. For the OTM the coloured regions of Figure \ref{fig:returnTime1} form $M_0$ and for $m >0$ each set $M_m$ is the union $\cup_{i,j} A_{i,j}^m$. For these latter sets, the authors assume that their measures decrease polynomially
\begin{equation}
\label{eq:muMm}
    \mu(M_m) \leq C_1 / m^r,
\end{equation}
where $r \geq 3$. Further they assume that if $z \in M_m$ then $F(z) \in M_k$ with
\begin{equation}
    \label{eq:betaCondition}
    \beta^{-1} m -C_2 \leq k \leq \beta m +C_3
\end{equation}
for some $\beta >1$ and unimportant constants $C_i>0$. It is straightforward to verify that (\ref{eq:muMm}) holds for $r=3$: each $A_{i,j}^m$ has length similar to $|\mathcal{L}_m| = \mathcal{O}(m^{-1})$, see (\ref{eq:Lklength}), and width similar to $|W_m| = \mathcal{O}(m^{-2})$, see (\ref{eq:Wklength}). Recalling the bounds $l_0$ and $l_1$ found in the proof of Proposition \ref{prop:upperTwoStep}, we see that for our map $H_\sigma$, condition (\ref{eq:betaCondition}) holds with $\beta=7$. In §4 of \cite{chernov_improved_2008} the authors describe an `ideal situation' under which the action of $F$ on the cells $M_m$ is equivalent to a discrete Markov chain. This requires:
\begin{enumerate}[label={(I\arabic*)}]
    \item The components of each $M_m$ and their images under $F$ are exact trapezoids which shrink homotetically as $m$ grows,
    \item The measure $\mu$ has constant density,
    \item $F$ is linear over each component,
    \item Condition (\ref{eq:betaCondition}) holds with $C_2 = C_3 =0$ (no irregular intersections).
\end{enumerate}
These conditions, together with:
\begin{equation}
\label{eq:conditionalProp}
   \frac{ \mu(F(M_m) \cap M_k)}{\mu(F(M_m))} = C_4 \frac{m}{k^2} + \mathcal{O}\left(\frac{1}{m^2}\right)
\end{equation}
for some $C_4>0$ and $k$ satisfying (\ref{eq:betaCondition}), are sufficient to establish the lemma. The authors show that their cells admit good linear approximations and the irregular intersections are of relative measure $\mathcal{O}(1/k)$ so that (I1) and (I4) are essentially satisfied, removing some portion of negligible measure from each cell. They then go on to estimate the effect of nonlinearity and nonuniform density of $\mu$ to address (I2) and (I3), requiring a more sophisticated approach. For our system (I2) and (I3) are already satisfied by $H_\sigma$, so it remains to verify (\ref{eq:conditionalProp}) for $H_\sigma$, show that our $A_{i,j}^k$ are well approximated by exact trapezoids, and calculate the relative measure of the irregular intersections. For (\ref{eq:conditionalProp}), areas of the regular intersections can be calculated using the shoelace formula on the corner coordinates $p_{k,l}$, $\overline{p}_{k,l}$, given explicitly in the appendix, section \ref{sec:upperCalculationsAppendix} and Proposition \ref{prop:lowerTwoStep}. Our cells are near exact trapezoids; unlike the billiards systems considered in \cite{chernov_improved_2008} whose cell sides are curvilinear, our cell boundaries are linear with the sides (e.g. $\mathcal{L}_m$, $\mathcal{L}_{m-1}$) near parallel for large $m$. For the irregular intersections, (\ref{eq:conditionalProp}) still gives an upper bound on their measure and there are some constant $C_2+C_3$ of them. The total number of intersections scales with $m$ by (\ref{eq:betaCondition}) so they have negligible measure compared to $\mu(H_\sigma(M_m))$.

\section{Final remarks}
\label{sec:discussion}
We provide some comments on the growth conditions which constituted the majority of our analysis in sections \ref{sec:Hmixing} and \ref{sec:Hsigma}. In the simplest cases of Lemma \ref{lemma:unstableGrowth}, growth was established in an analogous fashion to the old one-step expansion condition (\ref{eq:oldOneStepExpansion}), finding the relevant Jacobians $M_j$ and checking that their expansion factors $K(M_j)$ satisfy
\begin{equation}
    \label{eq:discussionOneStep}
    \sum_j \frac{1}{K(M_j)} <1.
\end{equation}
For the more complicated cases, the inductive method used to establish growth near the accumulation points in Lemma \ref{lemma:unstableGrowth} and the weakened one-step expansion condition (\ref{eq:oneStep}) both address the same fundamental issue: the splitting of unstable curves by singularities into an unbounded number of small components. They circumvent this obstacle in rather different ways, however. While (\ref{eq:oneStep}) generalises (\ref{eq:discussionOneStep}) to ensure an growth of unstable curves `on average' (see \cite{chernov_statistical_2009} for a precise statement), our inductive method is a more direct adaptation of (\ref{eq:discussionOneStep}), using it to generate contradictory geometric conditions which a hypothetical non-growing unstable curve must satisfy. It may be possible to prove Theorem \ref{sec:Hmixing} using (\ref{eq:oneStep}) as the basis for growth. Since we required (\ref{eq:oneStep}) anyway for proving Theorem \ref{thm:HsigmaExp}, this could potentially condense our analysis, but only to a minor extent. A convenience of the method used in section \ref{sec:Hmixing} is that, by way of the `simple intersection' property, it naturally gives geometric information on the images of manifolds, useful for proving the property \textbf{(M)} of Theorem \ref{thm:katok-strelcyn}.

We expect that essentially analogous analysis can be applied to establish mixing properties in a wide class of piecewise linear non-uniformly hyperbolic maps, including those (like the OTM) which sit on the boundary of ergodicity and beyond. While we have relied on the precise partition structure of $H_\sigma$, its fundamental feature (self-similar sequences of elements $A^k$, sharing boundaries with its neighbours $A^{k-1},A^{k+1}$ and accumulating onto some point $p$) is quite typical to return map systems. See, for example, those of various stadium billiards \cite{chernov_chaotic_2006,chernov_improved_2008,chernov_statistical_2009} and LTMs \cite{springham_polynomial_2014}. Indeed, the same method can be used to prove the Bernoulli property for non-monotonic LTMs \cite{myers_hill_mixing_2022}, where monotonicity of the manifold images cannot be assumed and the classical argument \cite{sturman_mathematical_2006} fails. The OTM is the pointwise limit of these maps as the boundary shrinks to null measure. It further has utility in proving growth conditions for maps which are uniformly hyperbolic but possess regions $A_j$ where the hyperbolicity is very weak, signified by $K(M_j) \approx 1$, so that (\ref{eq:discussionOneStep}) fails. Typically this leads to suboptimal bounds on mixing windows, see e.g. \cite{wojtkowski_model_1981,przytycki_ergodicity_1983,myers_hill_family_2022}. The map $H_{(\eta,\eta)}$ for $\eta \approx 1/2$ is another example, possessing weak hyperbolicity over $A_2, A_3$. Letting $\varepsilon = |\eta-1/2|>0$, there is an upper bound $N = N(\varepsilon)$ on escape times from the intersections $A_2\cap \sigma, A_3 \cap \sigma$. The growth lemma then follows by applying the inductive step roughly $N$ times and can be established for arbitrarily small $\varepsilon$, opening the door to establishing optimal mixing windows.

\begin{figure}
    \centering
    \includegraphics[width=0.24 \linewidth]{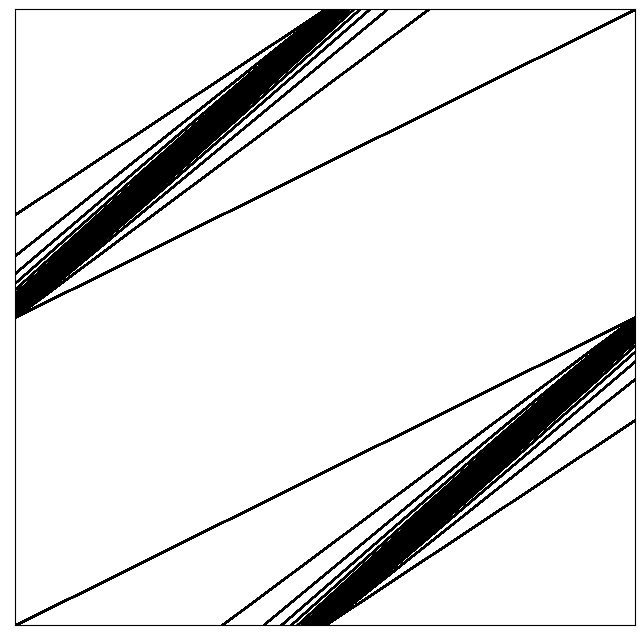}
    \includegraphics[width=0.24 \linewidth]{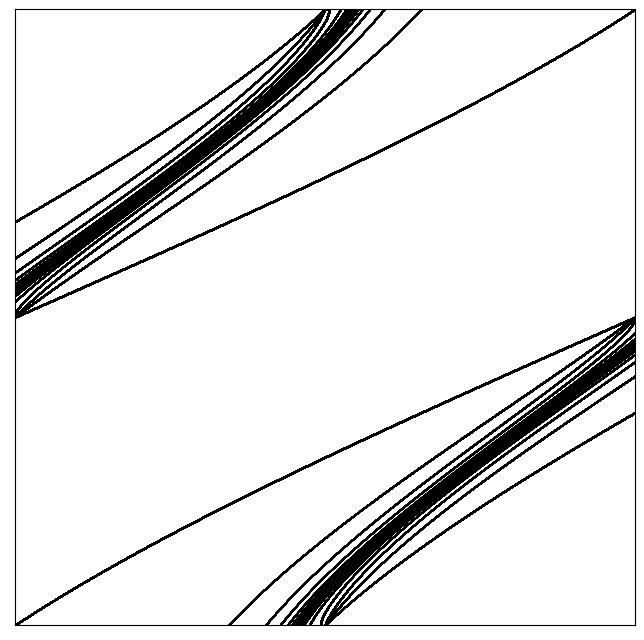}
    \includegraphics[width=0.24 \linewidth]{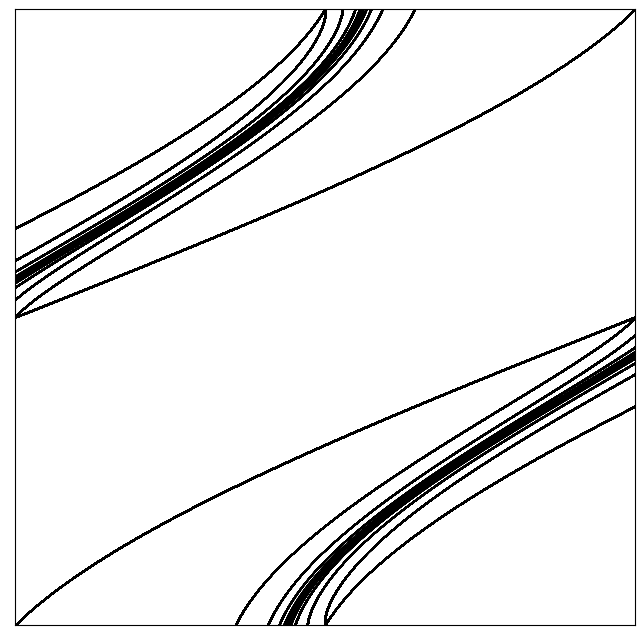}
    \includegraphics[width=0.24 \linewidth]{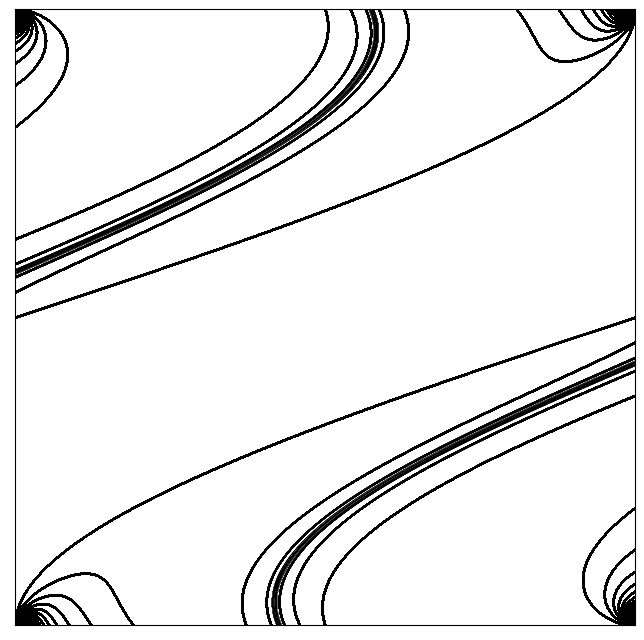}
    \caption{Partition of escape times from $S_3$ under the mapping $F \circ G$ for $a= 0,4,8,16$. }
    \label{fig:perturbations}
\end{figure}

The above gives two examples of piecewise linear perturbations to $H$ where mixing with respect to Lebesgue is preserved and our methods can be applied. Nonlinear perturbations to the shear profiles complicate the analysis in several ways. Firstly as the map's Jacobians takes on a broader range of values, cone invariance becomes an increasingly harder condition to establish. Cones must be widened, giving looser bounds on expansion factors, which may already be weak due to new regions of weaker stretching. This, together with the change from polygonal to curvilinear return time partition elements and nonlinear local manifolds, adds some complexity to showing growth conditions. This does not rule out certain (small) nonlinear perturbations however. There is some leeway in the inequalities which govern cone invariance and growth of local manifolds, the latter of which is not too dissimilar from the piecewise linear setting (see Lemmas \ref{lemma:piecewiseApprox}, \ref{lemma:componentLength}). Certain small perturbations would not alter the \emph{topological} structure of the return time partition, i.e. which elements share boundaries, the key information needed for setting up the induction. Finally while the partition elements would no longer be polygonal, only coarse geometric information is required for verifying each inductive step. Following the above, a potential perturbation could be to replace the linear portions of each shear by a cubic, perturbing the tent profile
\[  f(t) = \begin{cases} 2t & 0 \leq t \leq 1/2, \\ 2(1-t) & 1/2 \leq t \leq 1 ,\end{cases} \]
of the OTM shears to
\[  f_a(t) = \begin{cases} \frac{1}{8} t \left(16 - a + 6at - 8at^{2} \right) & 0 \leq t \leq 1/2, \\ \frac{1}{8}\left(1-t\right)\left( 16 - a + 6a\left(1-t\right) - 8a\left(1-t\right)^{2}\right)  & 1/2 \leq t \leq 1, \end{cases}   \]
for $a>0$. For small enough $a$ the gradient range $f'(t)$ is restricted to small neighbourhoods of $\{ 2, -2\}$ and the escape time partition retains a similar structure. We illustrate this in Figure \ref{fig:perturbations}, showing escapes from the square $S_3$ under the map $G \circ F$, equivalent to escapes from the perturbed $A_3$ under the $G \circ F$, but with a cleaner geometry for comparison. When $a$ is too large the analogy to the OTM breaks down. At $a=16$ the map is twice differentiable everywhere and features a new source of slowed mixing, the Jacobian is the identity at the corner points $x,y \in \{  0, 1/2 \}$ giving locally parabolic behaviour (visible in the escape time partition). 

\section{Appendix}
\subsection{Calculations for Proposition \ref{prop:upperTwoStep}}
\label{sec:upperCalculationsAppendix}
We begin by showing the bounds (\ref{eq:Ustarj}). Note that $U_{\star,j}$ is a curve traversing $A_{4,3}^j$ near $(1/4,1)$ with tangent vectors $(v_1,v_2)^T$ in the cone $M_1 \mathcal{C}_1$ satisfying $41/17 \leq |v_2|/|v_1| \leq 17/7$. Noting that the geometry of $A_{4,3}^j$ near $(1/4,1)$ is a $180^\circ$ rotation of $A_{4,3}^j$ near $(1/4,1/2)$ and the cone is invariant under this rotation, we can follow an analogous argument to (\ref{eq:Wdiamond}) to calculate $a_\star$. In particular $|U_{\star,j}|$ is bounded below by the length of the segment passing through $r_4$ with gradient $17/7$, which gives $a_\star = 13\sqrt{2}/80$ as required. For the upper bound, the height of $|U_{\star,j}|$ is bounded above by the height of the line segment with endpoints on $r_3(j)$ and $\mathcal{L}_{j-1}$ with gradient $41/17$. In particular
\[ \ell_v(U_{\star,j}) \leq \frac{48j^2 + 41j+29}{(2j+1)(48j+17)} - \frac{1}{2} = \frac{41}{2(96j^2+82j+17)} \sim \frac{41}{192\, j^2} \]
so that, by Lemma \ref{lemma:componentLength} and $|v_2/v_1| \geq 41/17$, $|U_{\star,j}| \geq b_\star/j^2$ with $b_\star = \frac{41}{192} \sqrt{1+17^2/41^2}$ as required.

We move onto calculating $b$ such that $|U_{k,l}| \lesssim b/l^2$. Define a $(k,l)$-cell as the intersection $H_\sigma\left(A_{4,3}^k \right) \cap A_{1,3}^l$ near the accumulation point $(1/2,3/4)$, shown as the magnified region in Figure \ref{fig:twoStepB}, the quadrilateral bounded by the lines $\mathcal{L}_l$, $\mathcal{L}_{l-1}$ (as defined in equation \ref{eq:linesLl}) on $\partial A_{1,3}^l$ and $\mathscr{L}_k$, $\mathscr{L}_{k-1}$ on $\partial H_\sigma \left(A_{4,3}^k \right)$. The explicit equation for $\mathscr{L}_k$ is given in (\ref{eq:cellLinesk}), letting us calculate the corner coordinates $p_{k,l} \in \mathscr{L}_k \cap \mathcal{L}_l$ as
\begin{equation}
    \label{eq:Pklcoords}
     p_{k,l} = \left( x_{k,l},y_{k,l} \right) =  \left(     \frac{16 k l + 7 k + 23 l + 10}{32 k l + 12 k + 44 l + 16} , \frac{12 k l + 3 k + 17 l + 4}{16 k l + 6 k + 22 l + 8}    \right).
\end{equation}
The curve $U_{k,l}$ traverses the $(k,l)$-cell with endpoints on the segments $p_{k,l}p_{k-1,l}$ and $p_{k,l-1}p_{k-1,l-1}$ and has tangent vectors in the cone $M_4M_3^k \mathcal{C}_1$. Roughly speaking, for large $k$ the vectors in this cone are essentially parallel to the cell boundaries $\mathscr{L}_k$, $\mathscr{L}_{k-1}$ with gradient approaching -3, so that $\ell_v(U_{k,l})$ is given to leading order by $y_{k,l} - y_{k,l-1} \sim \frac{3}{32} l^{-2} $. Noting that $M_4M_3^k \mathcal{C}_1 \subset \mathcal{C}_4$ for any $k$, we can bound the gradient of vectors as $|v_2/v_1| \geq 7/3$ so that by Lemma \ref{lemma:componentLength} we have $|U_{k,l}| \lesssim b/l^2$ with $b = \frac{3}{32}\sqrt{1+\frac{9}{49}}$.
A more careful calculation similar to that of $b_\star$ above gives the same bound to leading order.

\subsection{Two-step expansion near $P_1$}
\label{sec:lowerCalculationsAppendix}
We will follow similar analysis to the proof of Proposition \ref{prop:upperTwoStep} to show:
\begin{prop}
\label{prop:lowerTwoStep}
Condition (\ref{eq:oneStep}) holds for $H_\sigma^2$ when $W \cap B_\varepsilon(P_1) \neq \varnothing$ for all $\varepsilon>0$.
\end{prop}

\begin{proof}

We consider the case where $W$ lies near the accumulation point $(0,1/4)$, split by $\mathcal{S}$ into subcurves $\overline{W}_\star = W \cap A_1$ and $\overline{W}_k = W \cap A_{4,2}^k$. The image of the lower subcurve $\overline{U}_\star = H_\sigma \left(\overline{W}_\star \right)$ lies near the accumulation point $(1/2,1/4) = H(0,1/4)$ and is split by $\mathcal{S}$ into curves $\overline{U}_{\star,j} \subset A_{4,2}^j$. The image of each upper subcurve $\overline{U}_k = \overline{W}_k$ maps close to $(3/4,1/2)$ for $k$ odd, $(3/4,0)$ for $k$ even. Analysis for both of these cases is analogous, as before we take $k$ to be odd and consider the geometry of $\mathcal{S}$ near the accumulation point $(3/4,1/2)$. We calculate the corners of $A_{4,2}^k$ near $(0,1/4)$ as
\[  \overline{r}_1 = \left( 0, \frac{k}{4k-2}  \right), \quad \overline{r}_2 = \left(\frac{1}{4k-6}, \frac{k-2}{4k-6}  \right), \quad \overline{r}_3 = \left(\frac{1}{4k-2}, \frac{k-1}{4k-2}  \right), \quad \overline{r}_4 = \left(0, \frac{k+1}{4k+2}  \right).  \]
so that, using the integer valued matrix $M_4M_2^k = (-1)^k \big(\begin{smallmatrix}
  1-6k & -6k-2\\
  14k-2 & 14k+5
\end{smallmatrix}\big)$, its image $H_\sigma \left(A_{4,2}^k \right)$ is the quadrilateral with corners
\[  \overline{r}_1' = \left( \frac{3k+1}{4k-2}, \frac{2k-7}{4k-1}  \right), \quad \overline{r}_2' = \left(\frac{3k-2}{4k-6}, \frac{2k-9}{4k-6}  \right),\] \[ \overline{r}_3' = \left(\frac{3k-2}{4k-2}, \frac{k}{2k-1}  \right), \quad \overline{r}_4' = \left(\frac{3k+1}{4k+2}, \frac{k+1}{2k+1}  \right).  \]
The curve $\overline{U}_{k}$ has endpoints on the segments $\overline{r}_1'\overline{r}_2'$ and $\overline{r}_3'\overline{r}_4'$ and is split by $\mathcal{S}$ into an upper portion $\overline{U}_{k,\star}$ in $A_4$ above $y=1/2$ and subcurves $\overline{U}_{k,l} \subset A_{1,2}^l$ where $l_0 \leq l \leq l_1$. Comparison of the point $\overline{r}_2'$ with the lines $\overline{\mathcal{L}}_l: y -1/4 = -\frac{2l+1}{2l+2} (x-1)$ and $\overline{\mathcal{L}}_{l-1}$ on $\partial A_{1,2}^l$ yields $l_0(k) \geq \lfloor \frac{k-4}{7} \rfloor $, intersecting $\overline{r}_1'\overline{r}_4'$ with $y=1/2$ yields $l_1(k) \leq 7k+2$. Let $W_i = H_\sigma^{-1}\left(\overline{U}_i \right)$ then $W$ splits in an analogous fashion to (\ref{eq:Wsplitting}) with $DH_\sigma^2$ constant on each component. It follows that for $q=1/2$, 
\[ \liminf_{\delta \to 0} \sup_{W: |W|< \delta}  \sum_i \left( \frac{|W|}{|V_i|}\right)^q \frac{|W_i|}{|W|} \leq \sup_{0 \leq p \leq 1} \left( 2\sqrt \frac{(1-p)\overline{b}_\star \Lambda_1^+}{\overline{c}_\star \overline{a}_\star\Lambda_1^-} +    2\sqrt{\frac{p\overline{b}_\diamond}{\overline{c}_\diamond \overline{a} \overline{\gamma}}} +  4 \sqrt{\frac{p \overline{b} h}{\overline{c}\overline{a} \overline{\gamma}}} \right) \]
where the new constants satisfy (letting $K(M)$ denote the minimum expansion of $M$ over $\mathcal{C}_1$)
\begin{itemize}
    \item $K\left(M_4M_2^jM_1\right) \sim \overline{c}_\star j$
    \item $K\left(M_4^2M_2^k\right) \sim \overline{c}_\diamond k$
    \item $K\left(M_1M_2^lM_4M_2^k\right) \sim \overline{c} kl$
    \item $K\left(M_4M_2^k\right) \sim \overline{\gamma} k$
    \item $\overline{a}_\star /j^2 \lesssim |\overline{U}_{\star,j}| \leq \overline{b}_\star /j^2$
    \item $|\overline{U}_{k,\star}| \lesssim \overline{b}_\diamond / k$
    \item $|\overline{W}_k| \gtrsim \overline{a}/k^2$
    \item $|\overline{U}_{k,l}| \lesssim \overline{b}/l^2$
\end{itemize}
and $\Lambda_1^{\pm}$, $h$ are unchanged from (\ref{eq:justConstants}). The expansion factors can be calculated in the same fashion as (\ref{eq:Lambdastarj}), in particular
\[ \overline{c}_\star = \frac{48\sqrt{145}}{5}, \quad \overline{c}_\diamond = 8\sqrt{197},  \quad \overline{c} =64, \quad \overline{\gamma} = \frac{8\sqrt{145}}{5}.   \]
The constant $\overline{a}_\star$ is obtained by considering the shortest path across $A_{4,2}^j$ with tangent vectors aligned in the cone $M_1\mathcal{C}_1$, bounded by the length of the segment with endpoints on $r_4(j)$ and $\mathcal{L}_{j-1}$ (as defined in ($\dagger$), proof of Lemma \ref{lemma:unstableGrowth}) with gradient 41/17. The constant $\overline{b}_\star$ is obtained by considering the maximum height of a segment joining $\mathcal{L}_{j-1}$ to $\mathcal{L}_j$, given by the segment passing through $r_3(j) \in \mathcal{L}_j$ with gradient 17/17, and applying Lemma \ref{lemma:componentLength}. In particular $\overline{a}_\star = \sqrt{1970}/464$ and $\overline{b}_\star = \frac{17}{192} \sqrt{1+\frac{17^2}{41^2}}$. Similar analysis to the calculation of $\overline{a}_\star$ but using the wider cone $\mathcal{C}_1$ yields $\overline{a} = \sqrt{55}/80$. We again apply Lemma \ref{lemma:componentLength} to find $\overline{b}_\diamond$, with $\ell_v \left(\overline{U}_{k,\star} \right)$ bounded above by the height $1/(4k-2) \sim 1/(4k)$ of $r_3'(k)$ above $y=1/2$. Tangent vectors of $\overline{U}_{k,\star}$ lie in the cone $M_4M_2^k \mathcal{C}_1 \subset \mathcal{C}_4$ so that $\overline{b}_\diamond = \frac{1}{4}\sqrt{1 + 9/49}$ provides the upper bound. Finally we calculate $\overline{b}$, following a similar approach to section \ref{sec:upperCalculationsAppendix}. For each $k$ the segments $\overline{r}_1'\overline{r}_4'$ and $\overline{r}_2'\overline{r}_3'$ lie on the lines $\overline{\mathscr{L}}_k$ and $\overline{\mathscr{L}}_{k-1}$ respectively, with 
\[ \overline{\mathscr{L}}_k: y - \frac{k+1}{2k+1} = -\frac{14k+5}{6k+2} \left( x- \frac{3k+1}{4k+2} \right).\]
Define a $\overline{(k,l)}$ cell as the intersection of $H_\sigma(A_{4,2}^k) \cap A_{1,2}^l$, given the by quadrilateral bounded by the lines $\overline{\mathscr{L}}_k$, $\overline{\mathscr{L}}_{k-1}$, $\overline{\mathcal{L}}_l$, $\overline{\mathcal{L}}_{l-1}$. Its corners $\overline{p}_{k,l},\dots,\overline{p}_{k-1,l-1}$ are given by
\[ \overline{p}_{k,l} = \left( \frac{(3k+1)(2l+3)}{8kl+11k+3l+4}  ,\frac{16kl+15k+7l+6}{4(8kl+11k+3l+4)} \right)  \]
with (as before, to leading order terms for $k$ large) $\ell_v\left(\overline{U}_{k,l}\right)$ bounded above by the height of the segment joining $\overline{p}_{k-1,l}$ to $\overline{p}_{k-1,l-1}$, $\ell_v \left(\overline{p}_{k-1,l}\overline{p}_{k-1,l-1} \right) \sim 7/32 l^{-2}$. Again, tangent vectors to $\overline{U}_{k,l}$ lie in $\mathcal{C}_4$ so that $\overline{b} = \frac{7}{32}\sqrt{1 + 9/49}$ gives an upper bound $|\overline{U}_{k,l}| \lesssim \overline{b}/l^2$ by Lemma \ref{lemma:componentLength}.

As before we take
\[ \overline{s} =  2\sqrt \frac{\overline{b}_\star \Lambda_1^+}{\overline{c}_\star \overline{a}_\star\Lambda_1^-} \approx 0.186, \quad \overline{t} = 2\sqrt{\frac{\overline{b}_\diamond}{\overline{c}_\diamond \overline{a} \overline{\gamma}}} +  4 \sqrt{\frac{ \overline{b} h}{\overline{c}\overline{a} \overline{\gamma}}} \approx 0.488, \]
giving
\[ \liminf_{\delta \to 0} \sup_{W: |W|< \delta}  \sum_i \left( \frac{|W|}{|V_i|}\right)^q \frac{|W_i|}{|W|} \leq \overline{s}\sqrt{\frac{\overline{s}^2}{\overline{s}^2+\overline{t}^2}} + \overline{t}\sqrt{\frac{\overline{t}^2}{\overline{s}^2+\overline{t}^2}} \approx 0.522 < 1  \]
as required.

\end{proof}

\printbibliography


\end{document}